\documentclass{article}
\usepackage{geometry} 
\usepackage[utf8]{inputenc}
\usepackage{appendix}
\usepackage{units}
\usepackage {xfrac} 
\usepackage{amsfonts,amssymb}
\usepackage{microtype}
\usepackage{subfigure}
\usepackage{apacite}
\usepackage{booktabs} 
\graphicspath{{picture/}}
\usepackage{algorithm}
\usepackage{algorithmic}
\usepackage{multirow}
\usepackage{graphicx}
\usepackage{threeparttable} 
\usepackage{array}
\renewcommand{\arraystretch}{1.5}
\usepackage{bm}
\usepackage{tikz}
\usetikzlibrary{arrows.meta}
\usepackage{amstext}
\usepackage{amssymb} 
\usepackage{amsmath}
\usepackage{hyperref}
\usepackage{amstext}
\usepackage{hyperref}
\usepackage{float}
\usepackage{color}
\usepackage{indentfirst}
\usepackage[round]{natbib}
\usepackage{graphicx} 
\usepackage{amsthm}
\usepackage{titlesec}
\usepackage{sectsty}
\sectionfont{\fontsize{12}{14}\selectfont} 
\subsectionfont{\fontsize{10}{12}\selectfont}
\newtheorem{theorem}{Theorem}
\newtheorem{proposition}[theorem]{Proposition}
\newtheorem{lemma}{Lemma}
\newtheorem{remark}{Remark}

\newtheorem{definition}{Definition}
\newtheorem{assumption}[theorem]{Assumption}
\usepackage{rotating}      
\usepackage{tabularx}      
\usepackage{xltabular}   
\usepackage{makecell}    
\usepackage[table]{xcolor} 
\usepackage{longtable}
\usepackage{verbatim}
\usepackage{pdflscape}

\newcommand{\E}{\mathbb E}
\newcommand{\Var}{\operatorname{Var}}

\newcommand{\rhostar}{\rho^\ast}
\newcommand{\sigstar}{\sigma^\ast}
\newcommand{\F}{\mathcal F}
\newcommand{\G}{\mathcal G}

\begin{document}
\title{A Bayesian oriented exploration and reinforcement learning in continuous-time mean-variance portfolio optimization}
\author{%
  Kexin Chen\thanks{Department of Applied Mathematics, The Hong Kong Polytechnic
  University. Email: kexin-neme.chen@polyu.edu.hk}
  \qquad Xiang Yu\thanks{Department of Applied Mathematics, The Hong Kong
  Polytechnic University. Email: xiang.yu@polyu.edu.hk}
  \qquad Shanshan Zhao\thanks{Department of Applied Mathematics, The Hong Kong
  Polytechnic University. Email: shan2.zhao@connect.polyu.hk}
}

\date{\vspace{-0.3in}}
\maketitle
\begin{abstract}
We propose a new reinforcement learning (RL) framework for the continuous-time mean-variance portfolio optimization with unknown market parameters. To better reflect the impact of improved information of the unknown drift on the degree of exploration, we consider a novel randomization based on continuous-time stochastic filtering. Specifically, we construct a continuous-time analog to Thompson Sampling by treating the unknown drift as a random draw from a generalized sampling distribution parameterized by the filter's posterior mean and variance. The new learning objective measures the tradeoff between the wealth performance and the weighted KL divergence between the policy and the true Bayesian posterior, which transforms the classical feedback control into an information guided stochastic policy. As a result, the degree of exploration endogenously depends on the Reward-Biased Maximum Likelihood Estimation that decays with the shrinking filter's variance. This allows us to devise an efficient model-based learning algorithm in the mean-variance portfolio optimization with theoretical guarantee and superior learning performance.

\vspace{1em} 
\end{abstract}

\section{Introduction}
Mean--variance portfolio selection provides an intuitive way to balance the expected terminal wealth against risk,  which has been widely used in practice since the seminal work of \citet{markowitz1952portfolio} and its dynamic extensions in \citet{li2000optimal} for multi-period model and in  \citet{zhou2000continuous} for the continuous-time model. 
The endogenous linear-quadratic structure brings the favorable tractability in the dynamic control problem, however, relying crucially on the precise information of model parameters. In practice, the asset drift is unobservable and particularly difficult to estimate using the history data \citep{merton1980estimating}, leaving possibly large errors in implementing the theoretical optimal solution using the inaccurate drift \citep{best1991sensitivity}. A common response in the literature is to adopt the identify-then-optimize procedure used in adaptive control \citep{bielecki2019adaptive,kumar2015stochastic,faradonbeh2020adaptive}. Once the resulting estimates are treated as the true parameters, the remaining parameter uncertainty is hidden. Consequently, when the estimated model is substituted directly into the complete-information feedback policy, estimation errors might be further amplified \citep{best1991sensitivity,kan2007optimal}.

These limitations motivate the development of reinforcement learning (RL) approach that allows the agent to learn the optimal decision via the trial-and-error interactions with unknown environment. For continuous-time diffusion models, \citet{wang2020reinforcement} propose the exploratory formulation of stochastic policies via the entropy regularization. Within this framework, \citet{wang2020continuous} study the precommitted mean--variance problem and show that the optimal stochastic policy is Gaussian, which leads to a structured and implementable learning method. Along this direction, there is a fast growth in the development of entropy regularized reinforcement learning for continuous-time stochastic control problems and financial applications, see among \cite{jia2022policy}, \cite{JZhou22b}, \cite{GXZ22}, \cite{JZhou23}, \cite{Dong2023}, \cite{szpruch2024optimal}, \cite{wu2024reinforcement}, \cite{GRZ}, \cite{WYY24}, \cite{BHY25}, \cite{WY25}, \cite{RL-Merton}, \cite{FGLS},  \cite{jia2026accuracy}, \cite{huang2026meanvariance}, \cite{GLZ26}, \cite{BHYZ26}, \cite{HYZ26}, and references therein. 

In all aforementioned studies, although the entropy regularization provides a tractable way for the policy randomization, its induced exogenous noise in the randomized policy cannot reflect the agent's improving uncertainty on unknown model parameters along the learning. Intuitively, when the agent continues to learn, he would gain more information of the unknown environment and employ less randomness for the purpose of exploration. In addition, how to choose the temperature parameter in the standard entropy regularization to balance the tradeoff between exploration and exploitation in RL algorithms is still an open problem. On the other hand, Bayesian filtering quantifies that reduced unknownness through a posterior updated from observed prices, providing a statistical basis for calibrating exploration. The use of posterior beliefs to inform investment decisions is well established in portfolio optimization with partial information \citep{de2019bayesian,xiong2007mean}. \citet{lakner1995utility} derives optimal investment and consumption decisions when only prices are observed, while \citet{sass2004optimizing} use a filter for a hidden Markov-chain drift to characterize the terminal-wealth optimizer. \citet{bismuth2019portfolio} extend Bayesian learning to portfolio-choice, liquidation, and transition problems with unknown expected returns. In these models, the filter carries information of both the estimated parameter and its remaining uncertainty along the investment horizon. The economic relevance of retaining that uncertainty is illustrated by \citet{rogers2001relaxed}, which reveals that the welfare cost can exceed the cost of infrequent portfolio review.

Some recent studies simply combine this approach with the standard entropy-regularization in continuous-time RL. In a regime-switching market, \citet{wu2024reinforcement} incorporate the filtered regime belief into an actor--critic algorithm in the classical exploratory formulation. For an unknown constant drift, \citet{au2026entropy} combines Bayesian filtering with action-space entropy regularization and shows that the mean allocation remains the deterministic Bayesian feedback, while entropy changes the allocation variance. In these formulations, the randomness in the stochastic policy still comes exogenously from the entropy regularization. The filtered belief can affect the optimal allocation variance, but a positive entropy weight can sustain the randomization in the stochastic policy even when parameter uncertainty vanishes, as illustrated by the known-drift limit in \citet{au2026entropy}. For a small investor whose trades do not affect market prices, randomized policy does not change the return observations available to the filter and therefore does not provide additional information to improve the knowledge of the unknown drift. Such randomization can generate different wealth outcomes for policy evaluation, but it also incurs the cost of trying unnecessary suboptimal policies \citep{szpruch2024optimal}. 

\begin{figure}[!t]
\centering
\begin{tikzpicture}[
    >=Latex,
    box/.style={draw=blue!70!black, rounded corners=2pt, fill=blue!7,
        align=center, minimum height=12mm, text width=2.15cm,
        inner sep=3pt, font=\small, text=black},
    widebox/.style={draw=blue!70!black, rounded corners=2pt, fill=blue!10,
        align=center, minimum height=15mm, text width=3cm,
        inner sep=3pt, font=\small, text=black},
    flow/.style={->, thick, draw=blue!70!black},
    feedback/.style={->, thick, dashed, draw=blue!70!black},
    edge label/.style={font=\scriptsize, align=center, fill=white,
        inner sep=1pt, text=black}
]
\node[box] (market) at (0,0) {Market returns};
\node[box] (posterior) at (2.85,0) {Bayesian posterior};
\node[widebox] (exploration) at (6.10,0) {Posterior-guided,\\reward-biased\\Sharpe-ratio exploration};
\node[box] (action) at (9.55,0) {Mean--variance\\portfolio action};
\node[box] (outcome) at (12.40,0) {Wealth and reward};
\draw[flow] (market) -- (posterior);
\draw[flow] (posterior) -- (exploration);
\draw[flow] (exploration) -- (action);
\draw[flow] (action) -- (outcome);
\coordinate (loop-right) at (12.40,-1.55);
\coordinate (loop-left) at (6.10,-1.55);
\draw[feedback] (outcome.south) -- (loop-right) --
    node[edge label, below] {Reward feedback biases future Sharpe-ratio exploration toward better performance}
    (loop-left) -- (exploration.south);
\end{tikzpicture}
\caption{Posterior-guided reward-biased exploration.}
\label{fig:posterior-guided-framework}
\end{figure}
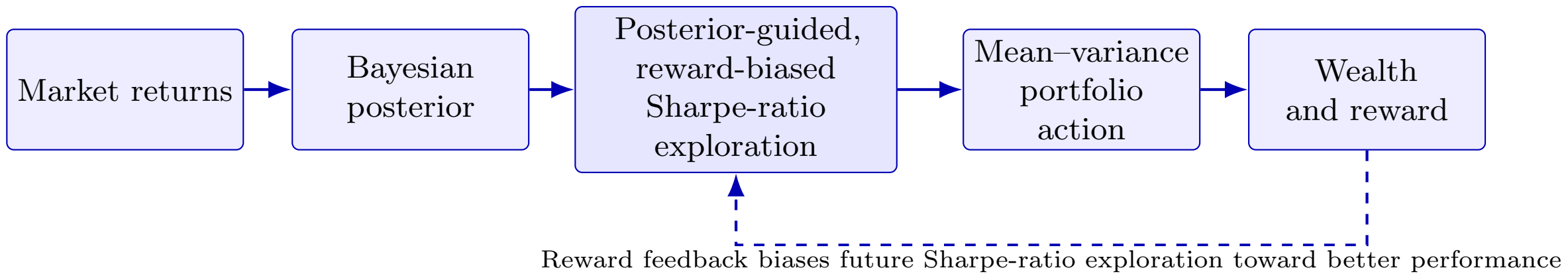

This work aims to provide a new perspective to address the issue of  over-randomization in standard entropy regularization as well as a proper choice of temperature parameter during the learning procedure when the agent gets more information on the unknown parameters by filtering. We contribute a novel continuous-time RL formulation that uses the posterior to guide the exploration search over candidate Sharpe ratios based on the classical mean--variance feedback form. At each time, the agent selects a sampling density $q_t$, and each sampled ratio generates a portfolio allocation through that feedback rule. The Bayesian filter supplies the current Gaussian posterior $p_t$ as a reference distribution. A Kullback--Leibler (KL) penalty on deviations of the sampling density from the posterior connects the parameter search to the information accumulated from returns. As the posterior variance decreases, it provides a shrinking effect for exploration. This construction is related to posterior sampling in Bayesian RL \citep{strens2000bayesian}, with the sampling law also allowed to respond to investment performance. Following the principle of reward-biased maximum likelihood estimation \citep{mete2021reward}, reward feedback can shift probability toward candidate parameters whose allocations yield better outcomes.

Our new formulation is built on the balance between the weighted posterior-relative KL penalty and the current policy value. This value scaling yields a multiplicative effect and preserves the quadratic structure inherited from the mean--variance objective, facilitating the explicit characterization of the optimal stochastic policy. For the algorithm aspect, the randomized allocations generate wealth outcomes for evaluating and improving portfolio policies. The Bayesian posterior continues to be updated from exogenous market returns independently of the investor's position, so randomization does not accelerate posterior contraction. In deriving the exploratory dynamics, we retain the common return innovation shared by wealth and the posterior mean. Figure~\ref{fig:posterior-guided-framework} summarizes the flow chart of how the posterior, parameter search, portfolio allocation, and reward feedback are implemented.

The main contributions of the present paper are summarized as follows:
\begin{itemize}
    \item First, we formulate a posterior-guided and reward-biased exploratory control framework, new to the literature, to learn the optimal policy in the context of mean-variance portfolio optimization. Its main novelty is a multiplicative criterion in which the posterior-relative KL divergence is weighted by both the exploration coefficient and the current policy value. This value scaling preserves the quadratic homogeneity inherited from the mean--variance objective.
    
    \item Second, we solve this exploratory control problem analytically. Notably, one state variable is an evolving probability density and the associated equation is nonlinear. However, by strategically choosing the temperature parameter in a proper form, we can still exploit the mean--variance structure and Gaussian posterior to reduce the optimizer to a unique Gaussian law and derive closed-form expressions for the value function and the stochastic policy, including the exceptional case in which the general solution is not directly applicable. The solution explicitly separates the role of posterior uncertainty in controlling the exploration scale from that of reward feedback, and shows that the degree of exploration declines automatically as the filtering information accumulates.
    \item Third, we turn the analytical result into an implementable learning method and establish some learning performance guarantees. To overcome the unknown oracle term in the closed-form policy and the instability caused by vanishing exploration variance, we develop a structure-preserving episodic actor--critic algorithm with specially normalized and centered martingale updates. We rigorously prove convergence of the filter, actor, sampling policy as well as Lagrange multiplier without requiring the critic proxy itself to identify the true Sharpe ratio. We further show that first-order policy errors cancel at the limiting optimum, leading to a nearly inverse-linear last-iterate performance rate in the stationary market.

    \item Fourth, to illustrate the learning performance of the proposed algorithm, we conduct both simulation experiments and empirical study by comparing benchmarks from the rolling maximum likelihood and entropy-regularized Bayesian-filtering. The results complement the theoretical findings and showcase the improved risk-adjusted performance and greater stability of learning.
\end{itemize}

The remainder of the paper is organized as follows. Section~2 introduces the classical benchmark, the Bayesian filter, and the exploratory state dynamics. Section~3 formulates the posterior-guided KL objective and derives the explicit characterization of optimal Gaussian stochastic policy. Section~4 presents the actor--critic model-based RL algorithm and establishes its convergence analysis. Sections~5 and~6 report the simulation experiments and empirical studies of the RL algorithm, respectively. The concluding section summarizes the findings, and the Appendix collects the proofs and auxiliary results.

\section{Bayesian Oriented Exploratory Formulation}
In this section, we first identify the dependence of the mean-variance feedback control on the market Sharpe ratio, and then introduce a filtering state and derive exploratory dynamics that preserve the common return innovation between wealth and beliefs.

\subsection{Classical Benchmark}
Let \(T>0\) be the finite horizon, and let \(\{W_t, 0 \leq t \leq T\}\) be a standard one-dimensional Brownian motion defined on a filtered probability space \((\Omega, \mathcal{F}, \{\mathcal{F}_t\}_{0 \leq t \leq T}, \mathbb{P})\) that satisfies the usual conditions. The price process of the risky asset is a geometric Brownian motion governed by

\begin{equation}\label{eq:St}
dS_t = S_t \left( \mu dt + \sigma dW_t \right), \quad  0 \leq t \leq T,
\end{equation}
with \(S_0 = s_0 > 0\) being the initial price, and \(\mu \in \mathbb{R}\), \(\sigma > 0\) being the mean and volatility parameters, respectively. The riskless asset has a constant interest rate \(r > 0\). The Sharpe ratio of the risky asset is defined by 
\begin{align}
    \rho = \dfrac{\mu - r}{\sigma}.
\end{align}

Denote by $\{ X^{u}_t, \, 0 \le t \le T \}$ the discounted wealth process under the portfolio strategy $u = \{ u_t, \, 0 \le t \le T \}$. Here, $u_t$ is the discounted dollar value allocated in the risky asset at time $t$, which satisfies the standard self-financing assumption and admissibility  conditions that will be given later. It follows from \eqref{eq:St} that the discounted wealth process satisfies
\begin{equation}\label{eq:X}
dX_t^u =\sigma u_t(\rho dt + dW_t),  \quad 0 \leq t \leq T,
\end{equation}
with initial wealth \( X^u_0 = x_0 \in \mathbb{R} \).

When the parameter $\bm{\theta}^*:=(\mu^*,\sigma^*)$ is known, the classical continuous-time Mean-Variance (MV) model aims to solve the optimization problem \citep{wang2020continuous}:
\begin{equation}\label{eq:V_cl}
V^{\mathrm{cl}}(s,y;w):= \min_{u\in {\cal A}^{\mathrm{cl}}(s,y)} \mathbb{E}\left[(X_T^u-w)^2\vert X_s^u = y\right]-(w-z)^2,
\end{equation}
for $(s,y)\in[0,T)\times\mathbb{R}$, where \( w \in \mathbb{R} \) is a fixed Lagrange multiplier. Once this problem is solved, $w$ can be determined by the constraint $\mathbb{E}[X_T^*] =z$, with $\{X^*_t\}_{t\in[0,T]}$ being the optimal wealth process under optimal portfolio $u^*$.

In this case, it is well known that the optimal control $u^*$ of \eqref{eq:V_cl} is given in a feedback form by 
\begin{align}\label{eq:u*}
    u^*= \psi(t,X_t^*,w;\bm{\theta}^*) = -\frac{\rho^*}{\sigma^*}(X_t^*-w),
\end{align}
and the resulting optimal wealth process is the unique strong solution to the SDE
\begin{align}
    d X_t^* = - (\rho^*)^2(X_t^*-w) dt - \rho^*(X_t^*-w) dW_t, \quad X_0^*=x_0.
\end{align}
If the parameter $\bm{\theta}^*$ is unknown to the agent, the agent searches for the optimal control to solve the MV model while learning the system, i.e., the parameter $\bm{\theta}^*$. However, if the agent is a small investor, the current portfolio selection problem has a distinct feature: the instantaneous returns of the risky and risk-free assets can be observed directly from the market, therefore, the information used for learning is unaffected by the portfolio strategies the agent applies. That said, the agent is always able to use an estimate $\hat{\bm{\theta}}$ (by observing market's time evolution) to solve the optimization problem while disregarding the modeling error, and use the resulting controller $\psi(\bm{\theta})$ on the true system with $\bm{\theta}^*$. That is, at time $t>0$, by observing $\{S_s\}_{s\in[0,t]}$, a maximum likelihood estimator (MLE) can be derived from maximizing log-likelihood: $\hat{\bm{\theta}}_t$. Then, a feedback control $\psi(\cdot;w,\hat{\bm{\theta}}_t)$ is exercised, and the state process evolves according to the dynamics \eqref{eq:X} controlled by policy  $\psi(\cdot;w,\hat{\bm{\theta}}_t)$:
\begin{align}\label{eq:Xphys_psi}
    dX_t =\sigma^* \psi(t,X_t ;w,\hat{\bm{\theta}}_t)(\rho^* dt + dW_t),  \quad 0 \leq t \leq T,
\end{align}
The agent then inherits the sampling error of $\hat{\bm{\theta}}_t$. Because these parameters are estimated from random samples, they are themselves random and introduce modeling errors that may propagate into undesirable behaviors over long time horizons. To address this, the next subsection retains the estimation uncertainty as an explicit state variable rather than treating the estimate as exact.

\subsection{Filtering and Exploratory Dynamics}
When $\sigma^*$ is known, the stochastic filtering theory characterizes the conditional distribution process $\{p_t:t\ge 0\}$ of $\rho^*$ \footnote{This is equivalent to filter $\mu^*$, but the filter on $\rho^*$ simplifies the formulation as $u^*$ is directly linked to $\rho$.}  given ${\cal S}_t$, where ${\cal S}_t$ is the (completed) $\sigma$-algebra generated by the observation (stock price process  $\{S_s\}_{s\in[0,t]}$) up to time $t$. We recall that in this setting, the conditional distribution $p_t$ of $\rho^*$ given ${\cal S}_t$ is a normal distribution, and hence completely characterized by its mean $m_t = \E[\rho^*|{\cal S}_t]$ and the variance $R_t =\E[(\rho^* -m_t)^2 |{\cal S}_t]$, which satisfy the following equation
\begin{align}\label{eq:SDE_mt}
    &d m_t = \frac{R_t}{\sigma^*}\Big(\frac{d S_t}{S_t} -r dt  - \sigma^* m_t dt\Big) =  R_t\big((\rho^* - m_t)dt + dW_t\big),\\
    & d R_t = -R_t^2 dt, \label{eq:SDE_Rt}
\end{align}
with the prior distribution of $\rho^*$ given by $N(m_0,R_0)$ at time $0$. It follows that posterior variance $R_t$ is deterministic and has closed-form: $R_t = R_0/(1+R_0 t)$, which means the uncertainty decays at rate $ O(1/t)$ regardless of the observations or policy. In addition, $m_t$ is a ${\cal S}_t$-martingale, converges to $\rho^*$ a.s. 

When \(\sigma^*\) is unknown, we estimate it from observed returns and substitute the estimate into the filtering equations, obtaining a plug-in filtering state \((m_t,R_t)\). Although this procedure provides an estimate of the unknown Sharpe ratio, uncertainty about that parameter remains and carries over to the portfolio allocation determined by the classical feedback rule. We therefore use the filtering state to guide a stochastic portfolio policy, denoted by \(\pi_t\), retaining parameter uncertainty in the policy construction rather than relying solely on a point estimate.

Recall that the classical optimal feedback control, $u^*=\psi(\cdot;\bm{\theta})$ as given by \eqref{eq:u*}, is explicitly determined by the true parameter. To construct the exploratory policy, we treat $\rho$ as a random variable drawn from a sampling distribution $q_t(\rho \mid m_t, R_t)$. This mapping naturally transforms the classical feedback control into a stochastic continuous-time policy $\pi_t(u \mid m_t, R_t)$. 
Crucially, $q_t$ need not be the exact Bayesian posterior $p_t$.
While this approach serves as a continuous-time analog of Thompson Sampling \citep{daniel2018tutorial}, relaxing the strict Bayesian assumption leaves mathematical room to integrate robust, generalized learning frameworks, such as the Reward-Biased Maximum Likelihood Estimation (RBMLE) approach which will be introduced in Section~\ref{sec:RBMLE}.

To rigorously analyze the macroscopic behavior of the state process under this randomized policy, we adopt the exploratory stochastic differential equation (SDE) framework introduced by \cite{wang2020continuous}. We consider $N$ independent identical agents, each experiencing an independent realization path of the physical Brownian motion $W_t^i$, $i=1,\cdots,N$, and updates their Bayesian posterior $(m^i_t,R_t)$, $i=1,\cdots,N$ identically. They sample independently their own control $u_t^i$, $i=1,\cdots,N$, from random policy $\pi_t$, resulting in $X_t^i$, $i=1,\cdots,N$ being copies of state process.
Then the increments of these state process copies are, for $i=1,\cdots,N$,
\begin{align}
    \Delta X_t^i =\sigma^*\rho^* u_t^i \Delta t + \sigma^* u_t^i (W^i_{t+\Delta t}-W^i_t), \quad \Delta m_t^i = R_t\big((\rho^*-  m_t^i)\Delta t + (W^i_{t+\Delta t}-W^i_t)\big), \quad t\ge 0.
\end{align}

Because the randomized policy $\pi_t$ is parameterized by the current Bayesian posterior $(m_t, R_t)$, it is inherently path-dependent. However, the randomization of the action $u_t^i$ is conditionally independent of the physical noise $\Delta W_t$. Therefore,  by the law of large numbers, as $N\rightarrow\infty$, the empirical averages of the state increments converge almost surely  to expectation taken over action space, while leaving the physical filtration of $W_t$ intact:
\begin{align*}
    &\frac{1}{N}\sum_{i=1}^N \Delta X_t^i \xrightarrow{a.s.} \mathbb{E}_{u \sim \pi_t}[ \Delta X_t^i ], \quad
    \frac{1}{N}\sum_{i=1}^N (\Delta X_t^i)^2 \xrightarrow{a.s.} \mathbb{E}_{u \sim \pi_t}[ (\Delta X_t^i)^2 ];\\
    &
    \frac{1}{N}\sum_{i=1}^N (\Delta X_t^i \Delta m_t^i) \xrightarrow{a.s.} \mathbb{E}_{u \sim \pi_t}[ \Delta X_t^i \Delta m_t^i ].
\end{align*}
This heuristic convergence motivates the definition of {\it exploratory state dynamics}, denoted by $X_t^{\pi}$. 
Each pair $(X^i, m^i)$, $i=1,\cdots,N$, can be viewed as an independent sample from the structural exploratory state dynamics $(X^{\pi},m)$. 
That is, each $(X^i,m^i)$ is generated according to the classical dynamic \eqref{eq:Xphys_psi} and \eqref{eq:SDE_mt}, with the corresponding $u$ sampled independently under policy $\pi$.
Under this framework, the effective drift, squared volatility, and cross-variations of $X_t^\pi$ are obtained by integrating the classical physical coefficients over the control distribution $\pi_t$.  

In our structural exploration method, the control is mapped from a parameterized distribution over the standardized drift, $\rho \sim q_t(\cdot | m_t, R_t)$, with mean $M_t$ and variance $V_t$. Substituting the action mapping $u_t = -{\rho}(X_t^\pi - w)/\sigma^*$, we explicitly evaluate these exploratory coefficients:
\begin{align*}
    &\frac{1}{\Delta t}\mathbb{E}_{u \sim \pi_t}[ \Delta X_t^i ] \xrightarrow{\Delta t \downarrow 0} \mathbb{E}_{\rho \sim q_t} [\sigma^* \rho^* u_t] 
    = \sigma^* \rho^* \mathbb{E}_{\rho \sim q_t} \bigg[ -\frac{\rho}{\sigma^*}(X_t^\pi - w) \bigg] 
    = - \rho^*(X_t^\pi - w) M_t.\\
    &\frac{1}{\Delta t}\mathbb{E}_{u \sim \pi_t}[ (\Delta X_t^i)^2 ]  \xrightarrow{\Delta t \downarrow 0}\mathbb{E}_{\rho \sim q_t} [(\sigma^* u_t)^2] 
    = (\sigma^*)^2 \mathbb{E}_{\rho \sim q_t} \bigg[ \frac{\rho^2}{(\sigma^*)^2}(X_t^\pi - w)^2 \bigg] 
    = (X_t^\pi - w)^2 (M_t^2 + V_t).\\
    &\frac{1}{\Delta t} \mathbb{E}_{u \sim \pi_t}[ \Delta X_t^i \Delta m_t ]\xrightarrow{\Delta t \downarrow 0}  \mathbb{E}_{\rho \sim q_t} [ (\sigma^* u_t) R_t ] 
    = \sigma^* R_t \mathbb{E}_{\rho \sim q_t} \bigg[ -\frac{\rho}{\sigma^*}(X_t^\pi - w) \bigg] 
    = - M_t R_t (X_t^\pi - w).
\end{align*}
Synthesizing these limits, we arrive at the complete macroscopic SDE system. Let $B$ be a standard Brownian motion independent of $W$. The exploratory state dynamics $X_t^\pi$ and the path-dependent learning process $m_t$ are governed by the following coupled system:
\begin{align}\label{eq:macroSDE}
    d X_t^\pi &= - {\rho^*} (X_t^\pi - w) M_t dt - (X_t^\pi - w)(M_tdW_t + \sqrt{V_t} dB_t), \cr
    d m_t &= R_t \big( (\rho^* - m_t) dt + dW_t \big), \quad dR_t = - R_t^2 dt, 
\end{align}
with the explicit cross-variation $d\langle X^\pi, m \rangle_t = - M_t R_t (X_t^\pi - w) dt$. 

\section{RBMLE-based Reinforcement Learning Objective}\label{sec:RBMLE}
Having established the exploratory state dynamics, we next formulate the objective of learning. 
In classical adaptive control, certainty-equivalent policies rely solely on an MLE estimator or purely Bayesian estimates $m_t$. 
The statistical estimator is inherently biased in favor of $\mathbb{P}^{\hat{\bm{\theta}}}$'s lower optimal rewards (or higher optimal cost) than under $\mathbb{P}^{\bm{\theta}^*}$. Therefore to extricate oneself from this bind, one must necessarily tilt the balance toward exploring parameters with higher optimal rewards (or lower optimal cost). Motivated by this,
 Reward Biased MLE (RBMLE) is proposed to counteract this statistical trap by introducing an optimistic bias, tilting the parameter estimation toward regions that promise higher optimal rewards (or lower optimal cost) \citep{mete2021reward}.

We employ this RBMLE principle into our continuous-time structural exploration framework. Rather than acting on a biased point estimate, the agent controls the trajectory of the system by selecting a dynamic distribution $q_t(\rho)$ over the Sharpe ratio $\rho$.
Mathematically, the agent seeks to minimize the expected terminal tracking cost $\E[(X_T^{\pi}-w)^2]$ regularized by the cumulative Kullback-Leibler (KL) divergence between the chosen exploratory policy $q_t$ and the true Bayesian posterior $p_t = N(m_t,R_t)$:
\begin{align*}
 J(t,x,m) = \inf_{q\in\mathcal Q(t,x,m)} \E\bigg[ (X_T^{\pi}-w)^2 + \int_t^T \lambda_s \eta(s,X_s^\pi,m_s) D_{\mathrm{KL}}\Big(q_s \parallel N(m_s, R_s)\Big) ds \Bigg| X_t^\pi = x, m_t = m \bigg],
\end{align*}
subject to system \eqref{eq:macroSDE}. In this formulation, the state-dependent weight $\eta(s,X_s^\pi,m_s)$ is chosen as $J(s,X_s^{\pi},m_s)$ to align with the magnitude of optimal objective function, and $\lambda_s>0$ controls the scale-invariant exploration-exploitation trade-off. 

It is worth noting the connection between our RBMLE objective and classical entropy-regularized reinforcement learning. Specifically, if we were to remove the Bayesian learning mechanism—ignoring the continuous-time filtration of observed data encapsulated by $(m_t, R_t)$—and instead replace the reference distribution $\mathcal{N}(m_s, R_s)$ with a static, uninformative prior, the KL divergence penalty would mathematically degenerate into a standard differential entropy regularizer ${\cal H}(q_t)=-\int_{\mathbb{R}}q_t(\rho) \log q_t(\rho) d\rho$. In this light, our formulation dynamically incorporates the agent's learning from observed data. This differs from \citet{au2026entropy}, where entropy regularizes the allocation distribution directly and the mean allocation remains the deterministic Bayesian feedback.

The optimization problem above cannot be posed over arbitrary density-valued processes. Progressive measurability alone does not guarantee that the policy moments are finite, that the controlled system \eqref{eq:macroSDE} is well posed, or that the terminal loss and cumulative KL penalty have finite expectation. These properties are also needed when applying It\^o's formula, taking conditional expectations, and invoking dynamic programming. Hence, Definition~\ref{def:admissible-distributional-control} introduces the admissible distributional controls $\mathcal Q(t,x,m)$ and thereby specifies the effective domain of the value function.

\begin{definition}\label{def:admissible-distributional-control}
Let $\mathbb F=(\mathcal F_s)_{s\in[t,T]}$ be the usual augmentation of the filtration generated by $(W,B)$. A distributional control $q=\{q_s\}_{s\in[t,T]}$ belongs to $\mathcal Q(t,x,m)$ if:
\begin{enumerate}
\renewcommand{\labelenumi}{(\roman{enumi})}
\item $q_s(\cdot)$ is a probability density on $\mathbb R$ almost surely for almost every $s$;
\item for every $A\in\mathcal B(\mathbb R)$, the process $\{\int_Aq_s(\rho)d\rho\}_{s\in[t,T]}$ is $\mathbb F$-progressively measurable;
\item its moments $M_s=\int_{\mathbb R}\rho q_s(\rho)d\rho$, 
$V_s=\int_{\mathbb R}(\rho-M_s)^2q_s(\rho)d\rho$ 
satisfy 
\[
\E\int_t^T(M_s^2+V_s)ds<\infty;
\]
\item the system \eqref{eq:macroSDE} has a unique strong solution satisfying
\[
\E\left[\sup_{t\le s\le T}
\bigl(|X_s^\pi|^2+|m_s|^2\bigr)\right]<\infty;
\]
\item Define $\Lambda_s^q
:=\exp\left\{\int_t^s\lambda_vD_{\mathrm{KL}}
\bigl(q_v\Vert N(m_v,R_v)\bigr)dv\right\}$, $t\le s\le T$.
The terminal integrating factor and the corresponding weighted terminal
loss satisfy
\[
\Lambda_T^q<\infty\quad\text{almost surely},
\quad
\E\left[\Lambda_T^q|X_T^\pi-w|^2\right]<\infty;
\]
\item For the positive classical HJB candidate $v$ used in
Theorem~\ref{thm:recursive-verification}, the family
\[
\left\{
\exp\left(
\int_t^\tau\lambda_sD_{\mathrm{KL}}
\bigl(q_s\Vert N(m_s,R_s)\bigr)ds
\right)
v(\tau,X_\tau^\pi,m_\tau):
\tau\text{ is an }[t,T]\text{-valued stopping time}
\right\}
\]
is uniformly integrable.
\end{enumerate} 
\end{definition}
Following the linear BSDE formulation of recursive utility in
\citet[Section~2.3]{dai2026merton}, the policy
value satisfies
\begin{align}
dJ_s^q={}&-\lambda_sD_{\mathrm{KL}}
\bigl(q_s\Vert N(m_s,R_s)\bigr)J_s^qds
+\mathbf Z_s^q\cdot d(W_s,B_s)^\top,
\quad
J_T^q=(X_T^\pi-w)^2.
\label{eq:recursive-policy-bsde}
\end{align}
Under condition~(v), its unique nonnegative solution for which
$\{\Lambda_\tau^qJ_\tau^q:\tau\le T\}$ is uniformly integrable is
\begin{equation}
J_s^q
=\frac{1}{\Lambda_s^q}
\E\left[
\Lambda_T^q(X_T^\pi-w)^2\mid\mathcal F_s
\right].
\label{eq:l1-policy-value-representation}
\end{equation}
Indeed, $\Lambda^qJ^q$ is the uniformly integrable martingale with terminal
value $\Lambda_T^q(X_T^\pi-w)^2$: condition~(v) makes this terminal
value integrable, so multiplying
\eqref{eq:l1-policy-value-representation} by $\Lambda_s^q$ identifies
$\Lambda_s^qJ_s^q$ with a uniformly integrable martingale. The Brownian
martingale representation theorem and the product rule then yield
\eqref{eq:recursive-policy-bsde}. Conversely, if another
solution satisfies the stated uniform-integrability condition, the product
rule makes $\Lambda^qJ^q$ a local martingale, while uniform integrability
upgrades it to a martingale. Its terminal value therefore forces the same
conditional-expectation representation. This $L^1$ solution concept and the
role of uniform integrability in uniqueness are standard; see
\citet[Theorems~6.2--6.3]{briand2003lp}. Tonelli's theorem and the tower
property show that \eqref{eq:l1-policy-value-representation} is equivalent
to the recursive relation \eqref{eq:policy-value-q}. Condition~(vi) is the single additional condition needed to
pass from localized It\^o identities to the terminal time in the verification
argument.

A feedback distributional policy is a Borel kernel $q(\rho\mid s,x,m)$ whose induced process $q_s(\rho)=q(\rho\mid s,X_s^\pi,m_s)$ belongs to this class for every initial state. The process $\{q_s\}$ is the corresponding open-loop control, whereas the deterministic kernel is the feedback (closed-loop) control. For an admissible feedback distributional policy $q$, let $\{q_s\}$ be
its induced open-loop control and let $X^\pi$ be the corresponding wealth
process. The policy-specific value function is defined recursively by
\begin{align}\label{eq:policy-value-q}
J^q(t,x,m)
:={}&\E\bigg[
(X_T^\pi-w)^2
+\int_t^T\lambda_sJ^q(s,X_s^\pi,m_s)
D_{\mathrm{KL}}\Big(q_s\Vert N(m_s,R_s)\Big)ds
\,\bigg|\,X_t^\pi=x,\ m_t=m
\bigg].
\end{align}
Consequently, the optimal value function satisfies $J(t,x,m)=\inf_{q\in\mathcal Q(t,x,m)}J^q(t,x,m)$.

For the Gaussian policies used below, $M_s$ and $V_s$ depend on $(s,m_s)$ but not on wealth. Since the filtering equation has a unique strong solution, the wealth equation is then linear and, whenever the displayed integrals are finite, has the unique solution
\begin{align}\label{eq:wealth-explicit-admissibility}
X_s^\pi-w=(x-w)\exp\bigg\{
&-\int_t^s\left(\rho^*M_v+\frac12(M_v^2+V_v)\right)dv
-\int_t^sM_vdW_v
-\int_t^s\sqrt{V_v}\,dB_v\bigg\}.
\end{align}
Thus, existence and pathwise uniqueness follow directly for the policy class considered in this paper. Meanwhile, the moment and cost bounds remain substantive admissibility conditions.

\subsection{Exploratory HJB Equation and Gaussian Optimal Policy}
Having specified a well-defined learning objective for the agent, we now proceed to solve for the optimal control. Since $R_s$ is deterministic for a fixed initial posterior variance, the controlled Markov state is $(X_s^\pi,m_s)$, with $R_s$ suppressed from the arguments of $J$. Motivated by the dynamic programming principle for relaxed controls
\citep{elkaroui1987compactification}, we consider the following Bellman relation 
\begin{align*}
 J(t,x,m) = \inf_{q\in\mathcal Q(t,x,m)} \E\bigg[ J(s,X_s^{\pi},m_s)  + \int_t^s \lambda_v J(v,X_v^{\pi},m_v) D_{\mathrm{KL}}\Big(q_v \parallel N(m_v, R_v)\Big) dv \Big| X_t^\pi = x, m_t = m \bigg].
\end{align*}
For fixed $(t,x,m)$ and fixed derivatives of $J$, the part of the controlled generator that depends on the sampled Sharpe ratio $\rho$ is
\begin{align*}
g_t(\rho)
:={}&-\rho^*\rho(x-w)J_x
+\frac12\rho^2(x-w)^2J_{xx}
-\rho R_t(x-w)J_{xm},
\end{align*}
which is quadratic in the sampled Sharpe ratio $\rho$. Thus the pointwise optimization over densities is
\begin{align}\label{eq:gibbs-variational-problem}
\inf_q\left\{
\int_{\mathbb R}g_t(\rho)q(\rho)d\rho
+\lambda_tJ\,D_{\mathrm{KL}}(q\Vert p_t)
\right\}.
\end{align}
For $x\ne w$, let 
\[
N_t:=\int_{\mathbb R}\exp\left\{-\frac{g_t(\rho)}{\lambda_tJ(t,x,m)}\right\}
p_t(\rho)d\rho.
\]
By the Donsker--Varadhan variational representation of relative entropy
\citep[Lemma~2.1]{donsker1975asymptotic}, or directly by the following KL
decomposition, for every admissible density $q$,
\begin{align}\label{eq:gibbs-kl-decomposition}
&\int_{\mathbb R}g_t(\rho)q(\rho)d\rho
+\lambda_tJ(t,x,m)D_{\mathrm{KL}}(q\Vert p_t)=-\lambda_tJ(t,x,m)\log N_t
+\lambda_tJ(t,x,m)D_{\mathrm{KL}}(q\Vert q_t^*).
\end{align}
where $$q_t^*(\rho)
={p_t(\rho)\exp\{-g_t(\rho)/(\lambda_tJ(t,x,m))\}}/{N_t}.$$ 
As the relative entropy is nonnegative, \eqref{eq:gibbs-kl-decomposition}
shows that $q_t^*$ uniquely minimizes \eqref{eq:gibbs-variational-problem},
provided $N_t<\infty$. Because $p_t$ is Gaussian and $g_t$ is quadratic,
the normalized Gibbs tilt $q_t^*$ is again Gaussian. Hence the
infinite-dimensional optimization may be written, without loss, in terms of
its mean $M$ and variance $V>0$.
Proceeding with standard arguments, we deduce that $J$ satisfies the Hamilton-Jacobi-Bellman (HJB) equation
\begin{align}\label{eq:HJB_J}
 J_t +& \inf_{M,V}\bigg\{ -{\rho^*}M (x-w) J_x + \frac{1}{2}(x-w)^2 (M^2+V) J_{xx} \cr
&+ R_t(\rho^*-m) J_m + \frac{1}{2} R_t^2 J_{mm} - M R_t(x-w)J_{xm} + \lambda_t J \times D_{\mathrm{KL}}(N(M,V)|N(m,R_t))  \bigg\} =0,
\end{align}
subject to $J(T, x,m) = (x-w)^2$. In \eqref{eq:HJB_J},  the KL divergence between two Gaussians with respective means and variances $M$, $V$, $m$, and $R$ can be written as
\begin{align*}
    D_{\mathrm{KL}}(q||p) = \frac{1}{2}\Big( \log R -\log V -1 + \frac{V+(M-m)^2}{R} \Big).
\end{align*}
The verification argument below uses the linear-BSDE representation of recursive utility.
\begin{theorem}
\label{thm:recursive-verification}
Let $v$ be continuous on $[0,T]\times\mathbb R^2$, belong to
$C^{1,2,2}$ on
$[0,T)\times(\mathbb R\setminus\{w\})\times\mathbb R$, be positive for
$x\ne w$, and satisfy \eqref{eq:HJB_J} and
$v(T,x,m)=(x-w)^2$. Suppose that its measurable Hamiltonian minimizer
$q^*$ is admissible in the sense of
Definition~\ref{def:admissible-distributional-control}.
Then
\[
v(t,x,m)\le J^q(t,x,m)
\quad\text{for every }q\in\mathcal Q(t,x,m),
\]
with equality for $q=q^*$. Consequently,
$v=J$ and $q^*$ is optimal.
\end{theorem}

The next subsection constructs the closed-form candidate. Its
candidate-specific admissibility, moment bounds and the uniform-integrability
condition in Definition~\ref{def:admissible-distributional-control} are
verified after the explicit solution is obtained.

\subsection{Closed-Form Solution}
The wealth dynamics are multiplicative in $X^\pi-w$ and the terminal loss is quadratic, so it is natural to seek a solution of $J$ in  the form $J(t,x,m) = A(t,m)(x-w)^2$, subject to $A(T,m) = 1$. 
Substituting this ansatz into \eqref{eq:HJB_J}, we have
\begin{align*}
    &A_t (x-w)^2 + \inf_{M,V}\bigg\{ -2A {\rho^*} M (x-w)^2 + A(M^2+V)(x-w)^2 + R_t(\rho^*-m)A_m (x-w)^2\\
    &\quad\quad\quad+ \frac{1}{2}R_t^2 A_{mm} (x-w)^2 -2A_m R_t M(x-w)^2 + \frac{\lambda_t A}{2} (x-w)^2 \Big( \log \frac {R_t} V -1 + \frac{V+(M-m)^2}{R_t} \Big) \bigg\} =0.
\end{align*}
The above equation is homogeneous in $(x-w)^2$, therefore the infimum over $M$ and $V$ completely decouples from $x$, the first order condition gives the optimal pair $(M^{\ast},V^{\ast})$ as follows: 
\begin{align*}
&A + \frac{\lambda_t A}{2} \left( -\frac{1}{V^{\ast}} + \frac{1}{R_t} \right) = 0
\Rightarrow \frac{1}{V^{\ast}}= \frac{1}{R_t} + \frac{2}{\lambda_t}\\
& 2AM^{\ast}- 2A{\rho^*} - 2A_m R_t + \lambda_t A \frac{M^{\ast} - m}{R_t} = 0 \Rightarrow  \frac{M^{\ast}}{V^{\ast}}= \frac{m}{R_t} + \frac{2}{\lambda_t} \Big({\rho^*} + \frac{A_m(t,m)}{A(t,m)} R_t \Big).
\end{align*}
By substituting the optimal structural controls $M^{\ast}$, $V^{\ast}$ back into the HJB \eqref{eq:HJB_J}, we obtain the PDE that governs $A(t,m)$ as follows,
\begin{align*}
A_t + \frac{1}{2} R_t^2 A_{mm} - \Gamma(t) \frac{A_m^2}{A} + (\alpha_A(t)m + \beta_A(t)) A_m + (\kappa_A(t)m^2 + \gamma_A(t)m + \delta_A(t)) A = 0,\; A(T,m)=1,
\end{align*}
with time-dependent coefficients given by 
\begin{align*}
 & \Gamma(t)=\frac{2 V^* R_t^2}{\lambda_t}, \quad \alpha_A(t)=-R_t-2 V^*, \\ 
 & \beta_A(t)=R_t \rho^*\left(1-\frac{4 V^*}{\lambda_t }\right), \quad
 \kappa_A(t)=\frac{\lambda_t}{2 R_t}-\frac{V^* \lambda_t}{2 R_t^2}=\frac{V^*}{R_t}\;(\text {after simplification}), \\ 
 & \gamma_A(t)=-\frac{2 V^* \rho^*}{R_t}, \quad \delta_A(t)=\frac{\lambda_t}{2} \log \left(\frac{R_t}{V^*}\right)-\frac{2 V^*\left(\rho^*\right)^2}{\lambda_t}.
\end{align*}
The PDE is semilinear because of the term $(A_m)^2/A$. 
We use a logarithmic transform to rewrite the PDE in terms of relative derivatives of $A$. Let
\begin{align*}
\ell(t,m) = \log A(t,m) \implies A(t,m) = e^{\ell(t,m)}.
\end{align*}
The direct substitution gives 
\begin{align}\label{eq:PDE_phi}
\ell_t + \frac{1}{2} R_t^2 \ell_{mm} + \left(\frac{1}{2} R_t^2 - \Gamma(t)\right) \ell_m^2 + (\alpha_A(t)m + \beta_A(t)) \ell_m + \kappa_A(t)m^2 + \gamma_A(t)m + \delta_A(t) = 0, 
\end{align}
with terminal condition $\ell(T,m)=0$. 

We strategically choose the temperature parameter as $\lambda_t=cR_t$, $c>0$, to ensure that the exploration cost scales proportionally with the current level of parameter uncertainty, thereby maintaining a consistent exploration-exploitation balance across the entire investment horizon. In this case, Equation \eqref{eq:PDE_phi} becomes
\begin{align*}
\ell_t + \frac{1}{2} R_t^2 \ell_{mm} + \frac{c-2}{2(c+2)}R_t^2\ell_m^2 + (\alpha_A(t)m + \beta_A(t)) \ell_m + \kappa_A(t)m^2 + \gamma_A(t)m + \delta_A(t) = 0, \quad \ell(T,m)=0.
\end{align*} 
Because all remaining dependence on $m$ is quadratic, we postulate the ansatz
\begin{align*}
\ell(t,m) = \frac{1}{2} P(t) m^2 + Q(t) m + K(t).
\end{align*}
When \(c\neq 2\), the coefficient of the quadratic term in the Riccati system is nonzero. Solving the resulting system yields
\begin{align*}
&\frac{1}{2} \dot{P} + \frac{c-2}{2(c+2)}R_t^2 P^2 + \alpha_A(t) P + \kappa_A(t) = 0, \; P(T)=0;\\
&\dot{Q} + \frac{c-2}{c+2} R_t^2 P Q + \alpha_A(t) Q + \beta_A(t) P + \gamma_A(t) = 0, \; Q(T)=0;\\
&\dot{K} + \frac{1}{2} R_t^2 P + \frac{c-2}{2(c+2)}R_t^2 Q^2 + \beta_A(t) Q + \delta_A(t) = 0, \; K(T)=0.
\end{align*}
The above ODE system admits the explicit solution:
\begin{align*}
P(t)&=\frac{1}{R_t}\frac{r_+r_-\left(1-\eta_t\right)}{r_+-r_-\eta_t},
\\
Q(t)&=\frac{\rho^*}{R_t}\left[\frac{r_+-r_-}{r_+-r_-\eta_t}\left(\frac{R_t}{R_T}\right)^{\frac{c+2-\Delta_c}{2(c+2)}}-1\right],
\\
K(t)&=\int_t^T\bigg[\frac{1}{2}R_s^2P(s)+\frac{c-2}{2(c+2)}R_s^2Q^2(s)+\rho^*\frac{c-2}{c+2}R_sQ(s)+\frac{cR_s}{2}\log\frac{c+2}{c}-\frac{2(\rho^*)^2}{c+2}\bigg]ds,
\end{align*}
where $\Delta_c = \sqrt{17c^2+36c+4}$, $r_{\pm} = {(5c+2\pm\Delta_c)}/{(2c-4)}$ and $\eta_t = \left({R_t}/{R_T}\right)^{-{\Delta_c}/{(c+2)}}$.
The resulting optimal value function is
\begin{equation}\label{optimal value}
J(t,x,m)
=
\exp\left(
\frac{1}{2}P(t)m^2+Q(t)m+K(t)
\right)(x-w)^2.
\end{equation}
Accordingly, the optimal pair $(M^{\ast},V^{\ast})$ is given by
\begin{equation}\label{optimal policy}
   V_t^{\ast} =\frac{c}{c+2}R_t,\quad
   M_t^{\ast}=\frac{c+2R_tP(t)}{c+2}m+\frac{2\rho^*+2R_tQ(t)}{c+2}.
\end{equation}
The case \(c=2\) is degenerate as the coefficient of the nonlinear term vanishes. As a result, the preceding closed-form representation containing factors of \((c-2)^{-1}\) cannot be evaluated at \(c=2\) directly. Setting \(c=2\) in the ODE system gives
\begin{align*}
&\frac{1}{2} \dot{P} - 2R_t P + \frac{1}{2} = 0, \quad P(T)=0;\\
&\dot{Q} -2R_t Q  - \rho^* = 0, \quad Q(T)=0;\\
&\dot{K} + \frac{1}{2} R_t^2 P +R_t\log2 - \frac{(\rho^*)^2}{2} = 0, \quad K(T)=0.
\end{align*}
Solving this degenerate linear system gives
\begin{align*} 
& P(t)=\frac{1}{3R_t}\left[1-\left(\frac{R_T}{R_t}\right)^3\right],\quad Q(t)=-\frac{\rho^*}{R_t}\left(1-\frac{R_T}{R_t}\right),\\
& K(t)=\left(\frac{1}{6}+\log 2\right)\log\frac{R_t}{R_T}+\frac{1}{18}\left(\frac{R_T^3}{R_t^3}-1\right)+\frac{(\rho^*)^2}{2}\left(t-T\right).
\end{align*}
These expressions coincide with the limit of the solution when \(c\neq2\) as \(c\to2\), so \(c=2\) can be regarded as the degenerate limiting case. Accordingly, the optimal pair $(M^{\ast}, V^{\ast})$ simplifies to
\begin{equation}\label{simplify policy}
   V_t^{\ast} =\frac{R_t}{2},\quad
   M_t^{\ast} = \Big(\frac{2}{3}-\frac{1}{6}\frac{R_T^3}{R_t^3}\Big)m+\frac{1}{2}\frac{R_T}{R_t}\rho^*.
\end{equation}

Finally, it is useful to separate the dependence on the true Sharpe ratio \(\rho^\ast\). For every \(c>0\), including the degenerate case \(c=2\), the value function can be written as
\begin{equation}\label{eq:unified_value}
J(t,x,m) = \exp\left\{\frac12 P(t)m^2 +\bar{Q}(t)\rho^\ast m + \frac12 K_{1}(t)(\rho^\ast)^2 + K_{2}(t) \right\}(x-w)^2 .
\end{equation}
Here \({Q}(t)=\rho^\ast \bar{Q}(t)\) and \( K(t)= K_{1}(t)(\rho^\ast)^2/2+K_{2}(t).\)
The explicit formulas obtained above for \(c\neq2\) and for \(c=2\) are precisely the closed-form solutions corresponding to this decomposition.

We next verify the candidate is
admissible. For fixed $(t,x,m)$ with $t<T$, write
\[
a_s:=\frac{c+2R_sP(s)}{c+2},
\quad
b_s:=\frac{2+2R_s\overline Q(s)}{c+2},
\]
Thus $M_s^*=a_sm_s+b_s\rho^*$, consistently with the actor parametrization used
later.

\begin{proposition}
\label{prop:closed-form-minimizer}
For every $c>0$, every finite horizon $T$, and every initial state
$(t,x,m)$ with $0\le t<T$, the Gaussian policy
$q_s^*=N(M_s^*,V_s^*)$ in \eqref{optimal policy}, including the
specialization \eqref{simplify policy}, belongs to
$\mathcal Q(t,x,m)$. Consequently, $q^*$ is optimal and the function in
\eqref{optimal value} is the optimal value function.
\end{proposition}

The complete proof is given in
Appendix~\ref{app:closed-form-admissibility}, where all six admissibility
conditions are verified directly on every finite horizon. Consequently, Proposition~\ref{prop:closed-form-minimizer} gives
$q^*\in\mathcal Q(t,x,m)$, and
Theorem~\ref{thm:recursive-verification} then shows that $q^*$ is optimal
and that the closed-form candidate in \eqref{optimal value} is the optimal value
function.


\section{Actor-Critic Learning Algorithm and Convergence Analysis}
\subsection{Actor-Critic Algorithm}
The filter tracks within-episode uncertainty, but the closed-form actor and value function still contain the fixed unknown $\rho^*$. We parameterize these occurrences and learn them from episodic trajectories: the critic evaluates the current policy, the actor updates the policy center, and the multiplier enforces the terminal mean constraint. 

By Bellman consistency applied to
\eqref{eq:policy-value-q}, for a fixed sampling policy $q$, with $X$ denoting the macroscopic wealth under \eqref{eq:macroSDE}, we
have
\[
J^{q}(t,x,m) = \mathbb E \left[J^{q}(s,X_s,m_s)+\int_t^s\lambda_v J^{q}(v,X_v,m_v)D_{\mathrm{KL}}\left(q_v\Vert N(m_v,R_v)\right)dv\bigg|X_t=x,\ m_t=m\right].
\]
For the policy evaluation, we adopt the  method based on martingale condition. More precisely, if $J^q$ is the true value function associated with $q$, then the cost-adjusted process
\[
\mathcal M_t^q = J^q(t,X_t,m_t)+\int_0^t\lambda_sJ^q(s,X_s,m_s)D_{\mathrm{KL}}\left(q_s\Vert N(m_s,R_s)\right)ds,
\]
is a martingale. Hence, for any predictable test function $\xi_t$ satisfying $\mathbb E\int_0^T\xi_t^2\,d\langle\mathcal M^q\rangle_t<\infty$ \citep{jia2022policy},
\[
\mathbb E
\left[
\int_0^T \xi_t\, d\mathcal M_t^q
\right]
=0.
\]
For policy improvement, \citet[Theorem~5]{jia2022policy} give the following policy-gradient representation, specialized here to an undiscounted additive KL criterion:
\[
\partial_\phi J(0,x,m;q^\phi)
=
\E\left[
\int_0^T
\partial_\phi\log q_t^\phi(\rho_t\mid m_t,R_t)
\,d\mathcal M_t^\phi
+
\int_0^T
\lambda_t\partial_\phi
D_{\mathrm{KL}}\left(q_t^\phi\Vert N(m_t,R_t)\right)dt
\right].
\]
Here $\mathcal M_t^\phi$ is the policy value plus accumulated running cost along the sampled trajectory in the additive model. We next adapt these principles to normalized actor--critic updates.

For the sampled-action updates, fix one episode and write $X_t$ for the microscopic wealth generated by the realized draw $\rho_t$. Denote by $J^\theta$ and $q^\phi$, respectively, the parameterized value function and sampling distribution. Along this trajectory, define the semimartingale residual
\begin{equation}\label{eq:main-parametric-residual}
d\mathcal R_t^{\theta,\phi}
:=
dJ^\theta(t,X_t,m_t)
+\lambda_tJ^\theta(t,X_t,m_t)
D_{\mathrm{KL}}\left(q_t^\phi\Vert N(m_t,R_t)\right)dt.
\end{equation}
This differential is the primitive object used in both the critic and actor updates. 

To approximate the residual in an implementable algorithm, we discretize $[0,T]$ into equal-length intervals \(0=t_0<t_1<\cdots<t_K=T\), with $t_i=i\Delta t$ and $\Delta t=T/K$. Along one episode, we collect samples \(D=\{(t_i,x_i,m_i,R_i)\}_{i=0}^K\). At time $t_i$, a Sharpe-ratio sample $\rho_i$ is drawn from $q^\phi_{t_i}(\cdot\mid m_i,R_i)$, and the corresponding risky allocation is \(u_i=-{\rho_i}(x_i-w)/{\widehat\sigma}\). The market simulator then generates the next wealth $x_{i+1}$, while the posterior state is updated using the observed return increment. Define the one-step residual increment
\begin{equation*}
\Delta \mathcal R_i^{\theta,\phi} = J^\theta(t_{i+1},x_{i+1},m_{i+1})-J^\theta(t_i,x_i,m_i)+\lambda_{t_i}J^\theta(t_i,x_i,m_i)D_{\mathrm{KL}}\left(q^\phi_{t_i}\Vert N(m_i,R_i)\right)\Delta t .
\end{equation*}
This is the finite-grid approximation of the increment of \eqref{eq:main-parametric-residual} over $[t_i,t_{i+1}]$.

We now specify the structured function approximators. In order to facilitate faster learning and convergence, we follow the structural form suggested by the closed-form solution in the previous section. According to Equation \eqref{optimal policy}, we focus on a Gaussian sampling distribution whose mean is affine in the posterior mean $m_t$ and in the learnable actor parameter $\phi$. Specifically, define
\[
q_t^\phi(\rho\mid m_t,R_t)
=
N\left(
\rho
\ \bigg|\
M^\phi(t,m_t,R_t),
V_t
\right),
\]
where
\begin{equation}\label{parameterized policy}
M^\phi(t,m_t,R_t)
=
\frac{c+2R_tP(t)}{c+2}m_t
+
\frac{2+2R_t\overline Q(t)}{c+2}\phi,
\quad
V_t
=
\frac{cR_t}{c+2}.
\end{equation}
Here $\phi$ is a learnable proxy for the unknown oracle Sharpe ratio $\rho^\ast$ in the actor. For later use, write $M_t^\phi:=M^\phi(t,m_t,R_t)$. The parameter $\phi$ is not assumed to equal $\rho^\ast$ during learning.

Similarly, motivated by Equation \eqref{eq:unified_value}, we parameterize the critic by
\[
J^\theta(t,x,m)
=
\exp
\left\{
\frac12P(t)m^2
+
\overline Q(t)\theta m
+
\frac12K_1(t)\theta^2
+
K_2(t)
\right\}
(x-w)^2.
\]
The critic parameter $\theta$ is another learnable proxy for the unknown Sharpe ratio in the value-function approximation. The parameter $\theta$ is also not required to equal $\rho^\ast$ during learning.

We next construct the actor--critic estimating equations, placing both updates on the same normalized-residual scale. Although the martingale-orthogonality framework of \citet{jia2022policy} permits various predictable test functions and a common gradient-type choice for the critic is  $\partial_\theta J^\theta$, we instead adopt a model-specific normalized test function
\begin{equation*}
\chi^\theta(t,x,m):=\frac{1}{J^\theta(t,x,m)^2}\frac{\partial J^\theta}{\partial\theta}(t,x,m)=\frac{1}{J^\theta(t,x,m)}\left(\overline Q(t)m+K_1(t)\theta\right).
\end{equation*}
Indeed, as
$\partial_\theta J^\theta=J^\theta\partial_\theta\log J^\theta$, its estimating increment satisfies
\[
\chi^\theta(t,X_t,m_t)d\mathcal R_t^{\theta,\phi}
=
\partial_\theta\log J^\theta(t,X_t,m_t)
\frac{d\mathcal R_t^{\theta,\phi}}{J^\theta(t,X_t,m_t)}.
\]
To see its effect, according to \eqref{eq:parametric-residual-decomposition} in Appendix~\ref{app:convergence-proofs}, write the normalized residual decomposition schematically as
\[
d\mathcal R_t^{\theta,\phi}
=J^\theta(t,X_t,m_t)
\left(G_t^{\rho_t^{(n)},\theta,\phi}dt+\Gamma_t^{\rho_t^{(n)},\theta,\phi}dW_t\right).
\]
As $\partial_\theta J^\theta=J^\theta\partial_\theta\log J^\theta$, the unnormalized choice $\partial_\theta J^\theta$ would produce
\[
\partial_\theta J^\theta\,d\mathcal R_t^{\theta,\phi}
=
(J^\theta)^2\partial_\theta\log J^\theta
\left(G_t^{\rho_t^{(n)},\theta,\phi}dt+\Gamma_t^{\rho_t^{(n)},\theta,\phi}dW_t\right).
\]
Because $(J^\theta)^2=(A^\theta)^2(X_t-w)^4$, this moment heavily weights trajectories with large wealth deviations and nearly suppresses critic learning when $X_t$ is close to $w$. By contrast, the proposed test function gives
\[
\frac{\partial_\theta J^\theta}{(J^\theta)^2}
d\mathcal R_t^{\theta,\phi}
=
\partial_\theta\log J^\theta
\left(G_t^{\rho_t^{(n)},\theta,\phi}dt+\Gamma_t^{\rho_t^{(n)},\theta,\phi}dW_t\right),
\]
so the homogeneous wealth factor cancels. The critic direction is consequently invariant under a common multiplicative rescaling of the value function and is less sensitive to the absolute wealth scale. The actor normalization by $(J^\theta)^{-1}$ performs the analogous cancellation, placing both updates on the common relative-residual scale $d\mathcal R_t^{\theta,\phi}/J^\theta$. These normalizations are therefore designed to improve scale stability in the present homogeneous mean--variance model.

In continuous time, the corresponding normalized critic estimating direction is
\begin{equation}\label{eq:main-normalized-critic-direction}
\widehat H(\theta,\phi)
:=
\int_0^T
\chi^\theta(t,X_t,m_t)\,d\mathcal R_t^{\theta,\phi}
=
\int_0^T
\frac{
\partial_\theta J^\theta(t,X_t,m_t)
}{
\bigl[J^\theta(t,X_t,m_t)\bigr]^2
}
d\mathcal R_t^{\theta,\phi}.
\end{equation}
The finite-grid approximation of \eqref{eq:main-normalized-critic-direction} gives the policy-evaluation update
\begin{equation}
\theta\leftarrow\theta+\alpha_\theta(n)\sum_i
\chi^\theta(t_i,x_i,m_i)\Delta \mathcal R_i^{\theta,\phi},
\end{equation}
where $\alpha_\theta(n)$ is the learning rate for the parameter $\theta$ at the $n$-th iteration.

For the actor, the offline policy-gradient representation motivates weighting the same normalized residual by the policy score. Since $q_t^\phi$ is Gaussian with mean $M_t^\phi$ and variance $V_t$, its score is
\begin{equation*}
\frac{\partial \log q_t^\phi}{\partial\phi}(\rho\mid m_t,R_t)
=\frac{\rho-M_t^\phi}{V_t}\cdot
\frac{2+2R_t\overline Q(t)}{c+2}.
\end{equation*}
The derivative of the KL-divergence term is
\begin{equation*}
\frac{\partial}{\partial\phi}D_{\mathrm{KL}}\left(q_t^\phi\Vert N(m_t,R_t)\right)
=\frac{M_t^\phi-m_t}{R_t}\cdot\frac{2+2R_t\overline Q(t)}{c+2}.
\end{equation*}
In continuous time, the corresponding normalized actor estimating direction is
\begin{equation*}
\widehat F(\phi,\theta)
=
\int_0^T
\frac{\partial\log q_t^\phi}{\partial\phi}(\rho_t\mid m_t,R_t)
\frac{d\mathcal R_t^{\theta,\phi}}{J^\theta(t,X_t,m_t)}
+
\int_0^T
\lambda_t\frac{\partial}{\partial\phi}
D_{\mathrm{KL}}\left(q_t^\phi\Vert N(m_t,R_t)\right)dt.
\end{equation*}

The two normalizations are well-defined under the continuous-time recursion used in the theory. For $x_0\ne w$, the sampled-action wealth deviation satisfies
\[
X_t-w
=
(x_0-w)
\exp\left\{
-\int_0^t\left(\rho^\ast\rho_s+\frac12\rho_s^2\right)ds
-\int_0^t\rho_s\,dW_s
\right\},
\]
and therefore never reaches zero. Since $J^\theta=A^\theta(t,m)(X_t-w)^2$ with $A^\theta>0$, we have $J^\theta(t,X_t,m_t)>0$ along every theoretical trajectory. The multiplier projections in Appendix~\ref{app:convergence-proofs} impose $w>x_0$, while the moment assumptions stated above and the polynomial bounds proved in the appendix ensure the required integrability after this cancellation.

We also center the normalized residual with respect to the current action before discretization. This centering step is necessary for exploration-scale stability. Write $\ell^\theta(t,m):=\log A^\theta(t,m)$ and let $\beta_{\mathrm{cv}}$ be a predictable scalar. At a sampling time, the action-independent baseline
\[
\left(G_t^{M_t^\phi,\theta,\phi}+\beta_{\mathrm{cv}}V_t\right)dt
+\Gamma_t^{M_t^\phi,\theta,\phi}dW_t
\]
has zero score-weighted conditional drift, because the conditional mean of the score is zero. The centered normalized residual becomes
\begin{equation}\label{eq:main-centered-normalized-residual}
d\mathcal C_t^{\theta,\phi}
=
\left[
(\rho_t-M_t^\phi)
\left(\rho_t+M_t^\phi-2m_t-2R_t\partial_m\ell^\theta(t,m_t)\right)
-\beta_{\mathrm{cv}}V_t
\right]dt
-\frac{2(\rho_t-M_t^\phi)}{R_t}\,dm_t.
\end{equation}
The corresponding continuous-time centered expression is
\begin{equation}\label{eq:main-centered-actor-direction}
\widehat F^{\mathrm c}(\phi,\theta)
=
\int_0^T
\frac{\partial\log q_t^\phi}{\partial\phi}(\rho_t\mid m_t,R_t)
d\mathcal C_t^{\theta,\phi}
+
\int_0^T
\lambda_t\frac{\partial}{\partial\phi}
D_{\mathrm{KL}}\left(q_t^\phi\Vert N(m_t,R_t)\right)dt.
\end{equation}
We take $\beta_{\mathrm{cv}}=1$ for literal conditional centering; the Gaussian score-weighted minimum-variance choice is $\beta_{\mathrm{cv}}=3$, which leaves the population direction unchanged. Importantly, centering is performed before discretization. Centering a finite-difference TD residual afterward leaves quadratic-variation remainders, such as $(\Delta W_i)^2-\Delta t$, which are action independent but are amplified by the score as $V_i\downarrow0$.
On the grid, define
\begin{equation}\label{eq:main-discrete-centered-increment}
\begin{aligned}
\Delta\mathcal C_i^{\theta,\phi}
:={}&
\Big[
(\rho_i-M_i^\phi)
\big(\rho_i+M_i^\phi-2m_i-2R_i\ell_{m,i}^\theta\big)
-\beta_{\mathrm{cv}}V_i
\Big]\Delta t-\frac{2(\rho_i-M_i^\phi)}{R_i}(m_{i+1}-m_i).
\end{aligned}
\end{equation}
The resulting finite-grid actor update is
\begin{equation}
\phi
\leftarrow
\phi
-
\alpha_\phi(n)
\sum_i
\left[
\frac{\partial \log q^\phi_{t_i}}{\partial\phi}(\rho_i\mid m_i,R_i)
\Delta\mathcal C_i^{\theta,\phi}
+
\lambda_{t_i}
\frac{\partial}{\partial\phi}
D_{\mathrm{KL}}
\left(
q^\phi_{t_i}\Vert N(m_i,R_i)
\right)
\Delta t
\right].
\end{equation}
Here $\alpha_\phi(n)$ is the actor learning rate. The action-independent KL drift cancels in the centered residual, while the explicit KL derivative is retained.

Our construction thus adopts the martingale-orthogonality and score-weighted policy-gradient principles of \citet{jia2022policy}, but does not use their estimating directions verbatim. The factor $(J^\theta)^{-2}$ in the critic, together with the factor $(J^\theta)^{-1}$ and the centering procedure in the actor, alters both the conditional means of the critic and actor estimating directions and the associated noise terms. Consequently, the convergence analysis cannot be obtained by directly transplanting an existing actor--critic proof. Instead, Appendix~\ref{app:convergence-proofs} rederives the conditional mean actor direction $F_n(\phi,\theta)$ defined in~\eqref{eq:two-part-decomposition} and establishes conditional second-moment bounds for the sampled actor estimating direction, while explicitly retaining the cross-variation between the wealth process and the filter.

Finally, we provide a scheme for learning the Lagrange multiplier $w$. The terminal mean constraint \(\mathbb E[X_T^*]=z\) suggests the stochastic approximation update
\[
w_{r+1}
=
w_r-\alpha_w(r)(\overline X_{T,r+1}-z),
\]
where $r$ indexes multiplier-update batches, $\alpha_w(r)>0$ is the corresponding learning rate, and $\overline X_{T,r+1}$ is the average of the $N_w$ terminal wealth values in batch $r+1$. This update is statistically self-correcting. If the average terminal wealth is above the target $z$, the update decreases $w$, which reduces the subsequent aggressive allocation. If the average terminal wealth is below the target, the update increases $w$. 

Algorithm~\ref{alg:offline_learning} presents the complete discrete-time algorithm. The episodic learning scheme carries the filtering mean forward across rounds and sets the initial filtering variance by a prescribed annealing schedule:
\begin{equation}
m_0^{(n+1)}=m_T^{(n)},
\quad
R_0^{(n)}
=
\frac{R_0}{(n+1)^\alpha},
\quad
R_0>0,\quad \alpha>0,\quad n\ge0.
\label{eq:episode-decay}
\end{equation}
Thus, within the $n$th episode, $R_t^{(n)}$ evolves according to the filtering Riccati equation. When entering the next round, the initial variance $R_0^{(n)}$ is reset according to \eqref{eq:episode-decay}. This reset separates the update of the belief mean from the schedule for the exploration intensity: the posterior mean is passed across episodes as the learned estimate, while the initial variance is specified by a deterministic annealing rule to control the exploration scale and the stochastic-approximation error.

\begin{algorithm}[t]
\caption{Offline-Episodic Actor--Critic Mean-Variance Algorithm}
\label{alg:offline_learning}
\textbf{Inputs}: initial wealth $x_0$, initial posterior pair $(m_0,R_0)$, horizon $T$, time step $\Delta t$, number of episodes $N$, number of time steps $K=\lfloor T/\Delta t\rfloor$, target payoff $z$, risk-free rate $r$, multiplier batch size $N_w$, initial learning rates $\alpha_\sigma,\alpha_\theta,\alpha_\phi,\alpha_w$ and temperature parameter $\lambda_t$. 

\textbf{Required program}: a market simulator $x'=\mathrm{Market}_{\Delta t}(t,x,m,a)$ that takes the current time-state pair $(t,x,m)$ and action $a$ as inputs and generates the next wealth state $x'$ at time $t+\Delta t$. The posterior state is updated from the observed return using the filtering update. 

\textbf{Learning procedure}: 
\begin{algorithmic}
    \STATE Initialize $\widehat\sigma$, $\theta$, $\phi$, and $w$.
    
    \FOR{episode $n=1$ \TO $N$}
        \STATE Initialize $k = 0$. Observe the initial state $x_{0} \leftarrow x_0$, $m_{0} \leftarrow m_0$, $R_{0} \leftarrow R_0$ and store $(t_0, x_{0}, m_{0})$.

        \WHILE{$k<K$}
            \STATE Compute and store the test function $\chi^\theta(t_k,x_k,m_k) = \frac{1}{J^\theta(t_k,x_k,m_k)^2}\frac{\partial J^\theta}{\partial\theta}(t_k,x_k,m_k)$.

            \STATE Generate a Sharpe ratio sample $\rho_{k} \sim q_{t_k}^\phi(\cdot | m_{k}, R_{k})$.

            \STATE Apply $u_k=-\frac{\rho_{k}}{\hat{\sigma}}(x_k-w)$ to the market simulator and observe the new wealth state $x_{k+1}$.

            \STATE Compute the plug-in whitened excess-return increment $\Delta\widetilde S_k = \frac{1}{\widehat\sigma}\left(\frac{S_{t_{k+1}}-S_{t_k}}{S_{t_k}}-r\Delta t\right)$.

            \STATE Update the posterior state $m_{k+1} = m_k+R_k(\Delta\widetilde S_k-m_k\Delta t)$, $R_{k+1} = {R_k}/{(1+R_k\Delta t)}$.

            \STATE Store transition $(t_{k+1},x_{k+1},m_{k+1},R_{k+1},\rho_{k})$.

            \STATE Update $k\leftarrow k+1$.
        \ENDWHILE

        \STATE Compute accumulated gradients:
        \begin{align*}
        &\Delta\theta =\sum_{i=0}^{K-1}\chi^\theta(t_i,x_i,m_i)\Delta\mathcal R_i^{\theta,\phi},\\
        &\Delta\phi =\sum_{i=0}^{K-1}\Bigg[\frac{\partial\log q_{t_i}^\phi}{\partial\phi}(\rho_{i}\mid m_i,R_i)\Delta\mathcal C_i^{\theta,\phi}+ \lambda_{t_i}\frac{\partial}{\partial\phi}D_{\mathrm{KL}}\left(q_{t_i}^\phi\Vert N(m_i,R_i)\right)\Delta t \Bigg].
        \end{align*}

        \STATE Update $\theta$ by $\theta \leftarrow \theta+ \alpha_\theta(n-1)\Delta\theta$.

        \STATE Update $\phi$ by $\phi\leftarrow\phi-\alpha_\phi(n-1)\Delta\phi$.

        \STATE Store the terminal wealth $X_T^{(n)}\leftarrow x_K$.

        \STATE \textbf{Update $w$ and $\widehat\sigma$ every $N_w$ episodes}:
        \IF{$n\equiv 0 \pmod {N_w}$}
            \STATE Update the Lagrange multiplier $w\leftarrow w - \alpha_w\big({n}/{N_w}-1\big) \left(\frac1{N_w} \sum_{\ell=n-N_w+1}^{n} X_T^{(\ell)} - z \right)$.

            \STATE Update the plug-in volatility estimate: $
            \widehat\sigma\leftarrow(1-\alpha_\sigma)\widehat\sigma+ \alpha_\sigma\widehat\sigma_{\rm MLE}$, where $\widehat\sigma_{\rm MLE}$ is estimated from the most recent observed returns.
        \ENDIF
    \ENDFOR
\end{algorithmic}
\end{algorithm}

\subsection{Convergence Analysis and Sharpe Ratio}

We next establish their convergence and portfolio implications. 
The results proceed from filtering consistency to actor and policy-center convergence, multiplier convergence, and the resulting Sharpe-ratio bound. There is an important distinction between the theoretical filter and its implementation. Equations \eqref{eq:SDE_mt}--\eqref{eq:SDE_Rt} are exact when $\sigma^\ast$ is known. When volatility is unknown, the standardized return uses the current estimate $\widehat\sigma_n$, so $(m_t^{(n)},R_t^{(n)})$ is a plug-in filtering state. In the ideal continuous-observation model, quadratic variation identifies volatility exactly through $T^{-1}\langle\log S\rangle_T=(\sigma^\ast)^2$. With discrete observations and no market microstructure noise, realized quadratic variation converges at rate $O_p(\Delta^{1/2})$ as the grid width $\Delta\to0$ \citep[Theorem~1]{barndorff2002realized}. Thus the continuous-time analysis below isolates the learning mechanism from a discretization error that vanishes under infill sampling. The feasibility of this approach, along with all other proofs, is provided in the Appendix~\ref{app:convergence-proofs}.

We first make explicit the theoretical sampled-action dynamics corresponding to Algorithm~\ref{alg:offline_learning}.  In the continuous-observation analysis, the sampled-action wealth $\widetilde X_t^\pi$ satisfies
\begin{equation}\label{eq:main-sampled-wealth-dynamics}
d\widetilde X_t^\pi
=-\rho^\ast\rho_t(\widetilde X_t^\pi-w)dt
-\rho_t(\widetilde X_t^\pi-w)dW_t.
\end{equation}
Since the filtering mean in \eqref{eq:SDE_mt} has diffusion coefficient $R_t$ with respect to the same Brownian motion $W$, the sampled wealth and the filter have the explicit cross-variation
\begin{equation*}
d\langle \widetilde X_t^\pi,m\rangle_t
=
-\rho_tR_t
\bigl(\widetilde X_t^\pi-w\bigr)dt.
\end{equation*}
Two state dynamics play distinct roles in the analysis. The policy value, the HJB equation, and the portfolio moments are defined under the macroscopic relaxed-control dynamics \eqref{eq:macroSDE}, whose coefficients are averaged over the sampling distribution. By contrast, the score-function actor statistic must retain the realized sampled action and is therefore constructed from the sampled-action microscopic dynamics.

Let $\mathcal F_n$ contain the information available before that episode. Let $\widehat H_{n+1}$ and $\widehat F_{n+1}^c$ denote the episode-indexed versions of the continuous-time estimating directions \eqref{eq:main-normalized-critic-direction} and \eqref{eq:main-centered-actor-direction} (with $\beta_{\mathrm{cv}}=1$), respectively, evaluated along $(X_t^{(n)},m_t^{(n)},R_t^{(n)},\rho_t^{(n)})_{0\le t\le T}$ with parameters $(\theta_n,\phi_n)$. The projected recursions are
\begin{equation}\label{eq:main-projected-actor-recursions}
\begin{aligned}
\theta_{n+1}
&=\Pi_\Theta\left[
\theta_n+\alpha_\theta(n)\widehat H_{n+1}(\theta_n,\phi_n)
\right],\\
\phi_{n+1}
&=\Pi_{\Phi_{n+1}}\left[
\phi_n-\alpha_\phi(n)\widehat F_{n+1}^c(\phi_n,\theta_n)
\right],\\
w_{r+1}&=\Pi_{\mathcal W_{r+1}}
\left[w_r-\alpha_w(r)(\overline X_{T,r+1}-z)\right].
\end{aligned}
\end{equation}
Here $\Theta=[-c_\theta,c_\theta]$ is a fixed compact interval with $c_\theta>|\rho^\ast|$. To make the expanding projections explicit, set
\begin{equation}\label{eq:main-expanding-actor-projection}
c_n:=1\vee\left[\log\log(n+\mathrm e^{\mathrm e})\right]^{1/4},
\quad
\Phi_n:=[-c_n,c_n],
\end{equation}
and choose $\underline w\in(x_0,w^\ast)$. For the batch counter $r$, set
\begin{equation}\label{eq:main-expanding-w-projection}
c_{w,r}:=\underline w+
\left(1\vee\left[\log\log(r+\mathrm e^{\mathrm e})\right]^{1/8}\right),
\quad
\mathcal W_r:=[\underline w,c_{w,r}].
\end{equation}
Thus, $n$ always indexes episodes, whereas $r$ indexes multiplier updates. Finally, define $\tau_t^{(n)}:=\log({R_t^{(n)}}/{R_T^{(n)}})$.
When the coefficient functions from Section~3.2 are evaluated along
$R^{(n)}$, denote them by $P_n,\overline Q_n,K_{1,n},K_{2,n}$.
Write $J_n^\theta$ for the corresponding critic. For simplicity, we define and subsequently use $L_n:=1+\log\log(n+\mathrm e^{\mathrm e})$.

\begin{assumption}\label{ass:main-learning-rates}
Take constants $\gamma_1$, $\gamma_2$, $\gamma_3$ satisfying $1/2<\gamma_1,\gamma_2,\gamma_3\le1$,
and take the filtering variance decay exponent $\alpha$ satisfying $0<\alpha\leq 1$, with $R_0T>1/2$ when $\alpha=1$.
Take constants $\alpha_\theta,\alpha_\phi,\alpha_w,\beta>0$, and let
\begin{equation*}
\alpha_\theta(n)=\frac{\alpha_\theta}{(n+\beta)^{\gamma_1}},
\quad
\alpha_\phi(n)=\frac{\alpha_\phi}{(n+\beta)^{\gamma_2}},
\quad
\alpha_w(r)=\frac{\alpha_w}{\bigl((r+1)N_w+\beta\bigr)^{\gamma_3}}.
\end{equation*}
By the above choices of $\gamma_1$, $\gamma_2$, $\gamma_3$ and $\alpha$,
we have
\begin{equation}
\begin{aligned}
&\sum_n\alpha_\theta(n)=\sum_n\alpha_\phi(n)=\infty,\quad\sum_r\alpha_w(r)=\infty,\\
&\sum_n\Big((\alpha_\theta(n))^2+(\alpha_\phi(n))^2\Big)<\infty,\quad \sum_r(\alpha_w(r))^2<\infty,\\
&\sum_n(\alpha_\phi(n))^2 c_n^4<\infty,
\quad
\sum_r(\alpha_w(r))^2(1+c_{w,r}^2)<\infty.
\end{aligned}
\label{eq:single-time-scale-summability}
\end{equation}
\end{assumption}
Assumption~\ref{ass:main-learning-rates} coordinates three distinct stochastic-approximation mechanisms. The nonsummability conditions allow the critic, actor, and multiplier to keep learning, whereas square summability controls their accumulated martingale noise. The remaining restriction $0<\alpha\leq1$ guarantees consistency of the episodic filter under the variance-reset schedule.

Next, Lemma~\ref{lem:main-filter} shows that passing the posterior mean across episodes while annealing the initial variance does not create an error floor.

\begin{lemma}\label{lem:main-filter}
Let $e_t^{(n)}:=m_t^{(n)}-\rho^\ast$ and $p_n:=1/R_0^{(n)}$. Then
\begin{equation}\label{eq:main-filter-recursion}
e_0^{(n+1)}
=\frac{p_n}{p_n+T}e_0^{(n)}
+\frac{1}{p_n+T}W_T^{(n)}.
\end{equation}
Each $e_0^{(n)}$ is Gaussian, $\sup_n\mathbb E|e_0^{(n)}|^q<\infty$ for every finite $q\ge1$, and constants $C,c_\alpha>0$ exist such that
\[
|\mathbb E[e_0^{(n)}]|
\le C(n+1)^{-\alpha/2},
\quad
\mathbb E[(e_0^{(n)})^2]\le C(n+1)^{-\alpha},
\quad
e_0^{(n)}\to0\quad\text{a.s.}
\]
\end{lemma}

Lemma~\ref{lem:main-filter} shows that both the filtering bias and mean-square error vanish under the variance-reset schedule. The unified bounds cover the inverse-linear endpoint $\alpha=1$ as well as $0<\alpha<1$, and the mean-square rate is $O((n+1)^{-\alpha})$. Thus the posterior component of the policy center converges to $\rho^\ast$ with uniformly controlled moments.

The conditional mean actor direction and noise calculations are given in Appendix~\ref{app:convergence-proofs}. They establish a uniformly positive stabilizing actor-drift coefficient, denote its lower bound by $\kappa_{c,R_0}>0$.  In particular, after the homogeneous wealth factor is cancelled by self-normalization, the conditional mean actor direction is independent of the absolute wealth scale and of the current predictable multiplier. Based on the stability analysis of the sampled actions, we obtain the following convergence result.

\begin{theorem}\label{thm:main-actor}
The actor recursion $\phi_n$ satisfies
\[
\phi_n\longrightarrow\rho^\ast\quad\text{a.s.}
\]
If $\beta=1$ and
$\kappa_{c,R_0}\alpha_\phi>\max\{\alpha,1/2\}$, then
\begin{equation}\label{eq:main-actor-rate}
\mathbb E|\phi_n-\rho^\ast|^2
\le CL_n\left[(n+1)^{-2\alpha}+(n+1)^{-\gamma_2}\right].
\end{equation}
\end{theorem}

Theorem~\ref{thm:main-actor} gives both consistency and a quantitative decomposition of the actor error. The term $(n+1)^{-2\alpha}$ is the squared contribution of the vanishing critic-induced perturbation, whereas $(n+1)^{-\gamma_2}$ is the stochastic-approximation error generated by the actor recursion. Hence, up to the slowly varying factor $L_n$, the mean-square convergence is governed by the slower of the two exponents $2\alpha$ and $\gamma_2$. The result also shows that convergence of the actor does not require the critic parameter itself to converge to $\rho^\ast$: the fixed projection keeps the critic bounded, while its influence on the actor field vanishes at rate $(n+1)^{-\alpha}$.

However, actor convergence alone does not yet identify the policy actually sampled, because its mean combines the posterior center and the actor parameter.  Let $\mathcal S(m,\phi;R)$ denote the theoretical Sharpe ratio generated by prior mean $m$, actor parameter $\phi$, and filtering covariance $R$. Define
\begin{equation*}
SR_{n}
:=
\mathcal S
\left(
m_0^{(n)},\phi_n;R^{(n)}
\right).
\end{equation*}
Denote the two weights in the zero filtering-variance limit by $a_c:=\frac{c}{c+2}$, $b_c:=\frac2{c+2}$, and
$a_c+b_c=1$. In a finite episode, the sum of the two weights in the actual policy mean
$M_t^{\phi_n}=a_t^{(n)}m_t^{(n)}+b_t^{(n)}\phi_n$ in \eqref{eq:ab-general} is generally not equal to one, but under the decaying variance rule these two weights converge respectively to $a_c=c/(c+2)$ and $b_c=2/(c+2)$, and their sum converges to one. The next result closes this gap.
\begin{theorem}\label{thm:main-policy-center}
Under the conditions of Theorem~\ref{thm:main-actor}, uniformly over $t\in[0,T]$,
\begin{equation}\label{eq:main-policy-coefficient-rate}
|a_t^{(n)}-a_c|+|b_t^{(n)}-b_c|=O((n+1)^{-\alpha}),
\end{equation}
and
\begin{equation}\label{eq:Mt-limit-center-variance}
\mathbb E|M_t^{\phi_n}-\rho^\ast|^2
\le CL_n\left[(n+1)^{-\alpha}+(n+1)^{-\gamma_2}\right].
\end{equation}
Hence, $\mathbb E[M_t^{\phi_n}]\to\rho^\ast$ and
$\operatorname{Var}(M_t^{\phi_n})\to0$.
\end{theorem}

The result implies concentration of the entire Gaussian sampling policy around the oracle Sharpe ratio.
Once the sampling policy concentrates at the oracle target, it remains to enforce the terminal mean constraint. The multiplier recursion operates on the slower batch counter $r$ and uses the terminal-wealth average to solve the corresponding scalar moment equation. The next theorem verifies that this additional constraint-learning step is compatible with the actor--filter limit.

\begin{theorem}\label{thm:main-w}
The multiplier recursion $w_r$ satisfies
\[
w_r\longrightarrow w^\ast
=\frac{ze^{(\rho^\ast)^2T}-x_0}{e^{(\rho^\ast)^2T}-1}
\quad\text{a.s.}
\]
\end{theorem}

Theorem~\ref{thm:main-w} shows that the batch recursion converges almost surely to the unique complete-information multiplier $w^\ast$. Consequently, the limiting policy not only has the oracle risky-asset direction but also satisfies the prescribed terminal mean target $\mathbb E[X_T]=z$.

The preceding results establish convergence of every component needed to implement the portfolio, but parameter convergence alone does not quantify economic performance. We therefore pass from the learned filter, actor, and policy center to the induced terminal-wealth moments and Sharpe ratio. 

For this calculation, suppress the episode superscript and consider a fixed actor parameter $\phi$ and a deterministic covariance path $R=\{R_t:0\le t\le T\}$. The mean and variance of the sampling distribution are
\begin{equation}\label{eq:ab-general}
M_t^\phi:=a_t m_t+b_t\phi,
\quad
V_t=\frac{cR_t}{c+2}.
\end{equation}
The actor parameter \(\phi\) is kept fixed throughout this subsection. Let $Y_t^\pi=X_t^\pi-w$. Then \eqref{eq:macroSDE} becomes
\[
dY_t^\pi
= -\rho^\ast M_t^\phi Y_t^\pi\,dt
-M_t^\phi Y_t^\pi\,dW_t
-\sqrt{V_t}\,Y_t^\pi\,dB_t .
\]
Assume $x_0\neq w$. Define $G(t,m):=\mathbb{E}\left[{Y_T^\pi}/{Y_t^\pi}|m_t=m\right]$. Then 
\[
\mathbb{E}\left[Y_T^\pi\mid Y_0^\pi=x_0-w,m_0\right]
=(x_0-w)G(0,m_0).
\]
By the Feynman--Kac formula, $G(t,m)$ satisfies
\[
G_t+\frac12R_t^2G_{mm}
+R_t(\rho^\ast-m-M_t^\phi)G_m
-\rho^\ast M_t^\phi G=0,
\quad
G(T,m)=1.
\]
As $M_t^\phi=a_tm+b_t\phi$ when the state variable is $m$, this PDE can be written as
\[
G_t+\frac12R_t^2G_{mm}
+R_t\bigl[\rho^\ast-(1+a_t)m-b_t\phi\bigr]G_m
-\rho^\ast(a_tm+b_t\phi)G=0,
\quad
G(T,m)=1.
\]
We use the exponential-affine ansatz
\[
G(t,m)=\exp\{p_G(t)m+q_G(t)\}.
\]
Substituting this ansatz into the PDE and matching the coefficients of $m$ and the constant terms gives
\begin{align*}
&p_G'(t)-R_t(1+a_t)p_G(t)-\rho^\ast a_t=0,
\quad p_G(T)=0,\\
&q_G'(t)+\frac12R_t^2p_G^2(t)
+R_t(\rho^\ast-b_t\phi)p_G(t)
-\rho^\ast b_t\phi=0,
\quad q_G(T)=0.
\end{align*}
Solving the above linear equations backward yields
\begin{align*}
p_G(t)
&=-\int_t^T \exp \left\{-\int_t^s R_u(1+a_u)du\right\}\rho^*a_sds,\\
q_G(t)
&=
\int_t^T
\left[
\frac12R_s^2p_G^2(s)
+R_s(\rho^\ast-b_s\phi)p_G(s)
-\rho^\ast b_s\phi
\right]ds.
\end{align*}
Therefore, $G(0,m_0)=\exp\{p_G(0)m_0+q_G(0)\}$, and hence
\begin{equation*}
\mathbb{E}[X_T^\pi]
=
w+(x_0-w)\exp\{p_G(0)m_0+q_G(0)\}.
\end{equation*}

We next compute the second moment of $Y_T^\pi$. Define $H(t,m):=
\mathbb{E}\left[{(Y_T^\pi)^2}/{(Y_t^\pi)^2}|m_t=m\right]$.
Then
\[
\mathbb{E}[(Y_T^\pi)^2\mid Y_0^\pi=x_0-w,m_0]=
(x_0-w)^2H(0,m_0).
\]
By the Feynman--Kac formula, $H(t,m)$ satisfies
\[
H_t+\frac12R_t^2H_{mm}
+R_t\bigl[\rho^\ast-(1+2a_t)m-2b_t\phi\bigr]H_m
+\left[
(a_tm+b_t\phi)^2+V_t
-2\rho^\ast(a_tm+b_t\phi)
\right]H=0,
\,
H(T,m)=1.
\]
We use the exponential-quadratic ansatz
\[
H(t,m)=
\exp\left\{
\frac12A_H(t)m^2+B_H(t)m+C_H(t)
\right\}.
\]
Substituting this ansatz into the PDE and matching the coefficients of $m^2$, $m$, and the constant terms gives
\begin{align*}
&A_H'(t)+R_t^2A_H^2(t)
-2R_t(1+2a_t)A_H(t)+2a_t^2=0,
\quad A_H(T)=0,\\
&B_H'(t)+R_t^2A_H(t)B_H(t)
+R_t(\rho^\ast-2b_t\phi)A_H(t)
-R_t(1+2a_t)B_H(t)
+2a_tb_t\phi-2\rho^\ast a_t=0,
\quad B_H(T)=0,\\
&C_H'(t)+\frac12R_t^2\bigl(A_H(t)+B_H^2(t)\bigr)
+R_t(\rho^\ast-2b_t\phi)B_H(t)
+b_t^2\phi^2+V_t-2\rho^\ast b_t\phi=0,
\quad C_H(T)=0.
\end{align*}
The first equation is a Riccati equation. Equivalently, by introducing an auxiliary function \(U(t)\) such that
\begin{equation}\label{A_H(t)}
A_H(t)=\frac{U'(t)}{R_t^2U(t)},
\end{equation}
we obtain the second-order linear terminal value problem
\[
U''(t)-4a_tR_tU'(t)+2a_t^2R_t^2U(t)=0,
\quad
U(T)=1,\quad U'(T)=0.
\]
Thus, \(A_H(t)\) is characterized by \eqref{A_H(t)}. Given \(A_H(t)\), the remaining two equations are linear and admit the following backward integral representations:
\begin{align*}
B_H(t)
&=
\int_t^T
\exp\left\{
\int_t^s
\left[
R_u^2A_H(u)-R_u(1+2a_u)
\right]du
\right\}\left(R_s(\rho^\ast-2b_s\phi)A_H(s)+2a_sb_s\phi-2\rho^\ast a_s\right)ds,\\
C_H(t)
&=
\int_t^T
\left[
\frac12R_s^2\bigl(A_H(s)+B_H^2(s)\bigr)
+R_s(\rho^\ast-2b_s\phi)B_H(s)
+b_s^2\phi^2
+V_s
-2\rho^\ast b_s\phi
\right]ds.
\end{align*}
It follows that the second moment around the target level $w$ is
\begin{equation*}
\mathbb{E}\left[(X_T^\pi-w)^2\right]
=
(x_0-w)^2
\exp\left\{
\frac12A_H(0)m_0^2+B_H(0)m_0+C_H(0)
\right\}.
\end{equation*}
Noting $\operatorname{Var}(X_T^\pi)=\operatorname{Var}(Y_T^\pi)=\mathbb{E}[(Y_T^\pi)^2]-\mathbb{E}[Y_T^\pi]^2$, we obtain
\[
\operatorname{Var}(X_T^\pi)
=
(x_0-w)^2
\left[
\exp\left\{
\frac12A_H(0)m_0^2+B_H(0)m_0+C_H(0)
\right\}
-
\exp\{2p_G(0)m_0+2q_G(0)\}
\right].
\]
By definition, $SR=\frac{\mathbb{E}[X_T^\pi]-x_0}{\sqrt{\operatorname{Var}(X_T^\pi)}}$, which can be expressed as
\begin{align*}
SR
&=
\frac{
(w-x_0)+(x_0-w)\exp\{p_G(0)m_0+q_G(0)\}
}{
\sqrt{
(x_0-w)^2
\left[
\exp\left\{
\frac12A_H(0)m_0^2+B_H(0)m_0+C_H(0)
\right\}
-
\exp\{2p_G(0)m_0+2q_G(0)\}
\right]
}
} \\
&=
\operatorname{sgn}(x_0-w)
\frac{
\exp\{p_G(0)m_0+q_G(0)\}-1
}{
\sqrt{
\exp\left\{
\frac12A_H(0)m_0^2+B_H(0)m_0+C_H(0)
\right\}
-
\exp\{2p_G(0)m_0+2q_G(0)\}
}
}.
\end{align*}
The preceding expression depends on the entire filtering path through the coefficient ODEs. The oracle Sharpe ratio is $SR_{\mathrm{or}}:=\sqrt{e^{(\rho^*)^2T}-1}$. 
When filtering variance is zero, the effective policy mean becomes the constant $\overline M_c(m,\phi)={(cm+2\phi)}/{(c+2)}$, 
and the corresponding Sharpe ratio is
\[
SR_\infty(M)
=\frac{e^{\rho^*MT}-1}{\sqrt{e^{M^2T}-1}}.
\]
This limiting expression separates parameter-learning error from the finite-$R$ perturbation of the coefficient ODEs.

\begin{lemma}\label{lem:main-sr-local}
Assume $\rho^*\ne0$ and $w>x_0$. It holds that
\[
SR_\infty(\rho^*)=SR_{\mathrm{or}},
\quad
SR_\infty'(\rho^*)=0,
\]
and
\[
SR_\infty''(\rho^*)
=-T\frac{e^{(\rho^*)^2T}
\left(e^{(\rho^*)^2T}-1-(\rho^*)^2T\right)}
{(e^{(\rho^*)^2T}-1)^{3/2}}<0.
\]
Consequently, for $|M-\rho^*|\le|\rho^*|/4$,
\[
|SR_\infty(M)-SR_{\mathrm{or}}|
\le C|M-\rho^*|^2.
\]
Moreover, there exists $\delta_{\mathrm{ODE}}>0$ such that, on
$|m-\rho^*|\vee|\phi-\rho^*|\le\delta_{\mathrm{ODE}}$, the limiting
normalized variance factor
$H^{(\infty)}(m,\phi)-G^{(\infty)}(m,\phi)^2$ has a positive lower
bound. For this $\delta_{\mathrm{ODE}}$, there exist $C<\infty$ and an
integer $n_{\mathrm{SR}}$ such that, for $n\ge n_{\mathrm{SR}}$,
\begin{equation}\label{eq:main-sr-ode-stability}
\left|
\mathcal S(m,\phi;R^{(n)})
-SR_\infty\left(\overline M_c(m,\phi)\right)
\right|
\le \frac{C}{(n+1)^\alpha}.
\end{equation}
\end{lemma}

Lemma~\ref{lem:main-sr-local} justifies the performance-robustness mechanism. Because $SR_\infty'(\rho^\ast)=0$, the first-order policy-center error cancels and the zero-covariance Sharpe loss is locally quadratic. Hence filtering and actor estimation errors affect the limiting Sharpe ratio only through their squared magnitude. The remaining $O(n^{-\alpha})$ term is the first-order perturbation caused by the nonzero filtering covariance and exploration variance. Combining these two effects with Lemma~\ref{lem:main-filter} and Theorem~\ref{thm:main-actor} yields the portfolio-level guarantee below.

\begin{theorem}\label{thm:main-sr}
Under conditions of Theorem~\ref{thm:main-policy-center} and Lemma~\ref{lem:main-sr-local}, it holds that
\begin{equation}\label{eq:main-sr-rate}
|SR_n-SR_{\mathrm{or}}|
=O_p\left(
\frac{L_n}{(n+1)^\alpha}
+\frac{L_n}{(n+1)^{\gamma_2}}
\right).
\end{equation}
\end{theorem}

Theorem~\ref{thm:main-sr} underscores the parameter-convergence in the portfolio-performance guarantee. The term $L_n(n+1)^{-\alpha}$ collects the effects of filtering uncertainty and the finite-covariance perturbation of the moment ODEs, whereas $L_n(n+1)^{-\gamma_2}$ is inherited from the actor mean-square error. Consequently, the slower of $\alpha$ and $\gamma_2$ determines the asymptotic Sharpe-ratio rate. Lemma~\ref{lem:main-sr-local} explains why these componentwise errors transfer in this form: the oracle is a stationary point of the limiting Sharpe ratio, so policy-center errors enter locally through their squared magnitude, while the remaining finite-covariance perturbation contributes at first order.

In a stationary market, observations from successive episodes concern the
same fixed market parameter, and the natural variance-reset schedule is
usually
\[
\alpha=1,
\quad
R_0^{(n)}=\frac{R_0}{n+1}.
\]
Take also the harmonic actor step size
$\alpha_\phi(n)=\alpha_\phi/(n+1)$, so that $\gamma_2=1$. Under the $\alpha=1$ case in Lemma~\ref{lem:main-filter} (which assumes $R_0T>1/2$ and gives the required $O((n+1)^{-1})$ filtering
mean-square rate), and provided
$\kappa_{c,R_0}\alpha_\phi>1$, the rate in
\eqref{eq:main-sr-rate} becomes
\[
|SR_n-SR_{\mathrm{or}}|
=O_p\left(\frac{\log\log n}{n}\right).
\]
The critic and multiplier exponents $\gamma_1,\gamma_3\in(1/2,1]$ ensure their respective stochastic-approximation conditions but do not enter the displayed Sharpe-rate bound directly.

\section{Simulation Experiments}
In this section, we conduct numerical simulations to compare the performance of the proposed posterior-guided actor--critic method with several benchmark methods for the continuous-time mean--variance problem. We compare three methods: Rolling MLE, Entropy-Regularized Bayesian Filtering and the proposed Actor--Critic method. Rolling MLE represents the classical plug-in/adaptive-control benchmark. Entropy-Regularized Bayesian Filtering refers to the entropy-regularized Bayesian drift-uncertainty method, in which the posterior state $(m_t,R_t)$ is updated by filtering and the portfolio allocation is sampled directly from the entropy-regularized Gaussian policy. Finally, the Actor--Critic method learns the sampling distribution through the martingale actor--critic updates developed in Section 4. Throughout the simulations, we set the weight parameter in the KL-regularized formulation to $c=2$. 

\subsection*{Maximum likelihood estimation}

Maximum-likelihood estimation (MLE) is a standard method for estimating the parameters $\mu$ and $\sigma$ in the geometric Brownian motion model (\ref{eq:St}). We refer the interested reader to \cite[Subsection~9.3.2]{campbell1997econometrics} for a detailed description of this method. At each decision-making time $t_i$, the MLE estimators for $\mu$ and $\sigma$ are calculated based on the most recent 100 data points of the price. One can then substitute the MLE estimators into the optimal allocation (\ref{eq:u*}) and the expression of the Lagrange multiplier $w=\frac{z e^{\rho^{2}T}-x_0}{e^{\rho^{2}T}-1}$ to compute the allocation $u_i\in\mathbb{R}$ in the risky asset. Such a two-phase procedure is commonly used in adaptive control, where the first phase is identification and the second one is optimization (see, e.g., \cite{bielecki2019adaptive,kumar2015stochastic}). The real-time estimation procedure also allows the MLE to be applicable in nonstationary markets with time-varying $\mu$ and $\sigma$.

\subsection*{Bayesian filtering}
The Bayesian Filtering benchmark in this section refers to the entropy-regularized Bayesian drift-uncertainty method \citep{au2026entropy}. It combines Kalman--Bucy filtering for the unknown Sharpe ratio with entropy-regularized relaxed control. This is different from simple posterior sampling: the benchmark samples the risky allocation directly from the entropy-regularized optimal Gaussian policy.

As the posterior learning mechanism is the same as that adopted in our framework, we use the posterior mean $m_i$ and posterior variance $R_i$ generated by the filtering equations. At each decision-making time $t_i$, given the current wealth $x_i$ and posterior state $(m_i,R_i)$, the entropy-regularized Bayesian policy is 
\begin{equation}
\pi_{t_i}^{\mathrm{Ent}}
=
N\left(
\bar{u}_{t_i}^{\mathrm{Ent}},
\left(\varsigma_{t_i}^{\mathrm{Ent}}\right)^2
\right),
\end{equation}
where
\begin{align*}
&\bar{u}_{t_i}^{\mathrm{Ent}}
=
-\frac{m_{i}\left[1+2R_{i}\alpha^{\mathrm{Ent}}(t_i)\right]}{\sigma}
\left(x_{i}-w\right),\;\alpha^{\mathrm{Ent}}(t)=-\frac{T-t}{1+2R_t(T-t)};\\
&\left(\varsigma_{t_i}^{\mathrm{Ent}}\right)^2
=
\frac{\varepsilon_{\mathrm{Ent}}}{2\sigma^2 A^{\mathrm{Ent}}(t_i,m_{i})},\; A^{\mathrm{Ent}}(t,m)=\exp\left\{\frac{(t-T)m^2}{1+2R_t(T-t)}+\frac{1}{2} \log\frac{(1+R_0t)\{1+R_0(2T-t)\}}{(1+R_0T)^2}\right\}.
\end{align*}
Here $\varepsilon_{\mathrm{Ent}}>0$ is the entropy-regularization parameter. In the theoretical entropy-regularized Bayesian model, $\sigma$ is assumed known. In our implementation, when volatility is estimated, we use the current plug-in estimate $\widehat\sigma_n$ in the allocation formula, held fixed within episode $n$.

At each time $t_i$, we sample $ u_i^{\rm Ent}$ and apply this sampled dollar allocation to the wealth dynamics. The posterior state is then updated from the observed return. The entropy-regularized Bayesian formulation identifies $w$ as the Lagrange multiplier associated with the terminal mean constraint $\mathbb E[X_T^*]=z$, but it does not provide a separate numerical updating scheme for $w$. To ensure a fair comparison, we therefore update $w$ using the same batch-based stochastic approximation procedure employed in our proposed method. This benchmark isolates the effect of entropy regularization under Bayesian drift uncertainty. In particular, its mean control coincides with the deterministic Bayesian Markowitz feedback, while the entropy term affects the policy variance. Therefore, it should not be described as simple Bayesian posterior sampling.


\subsection*{The stationary market case}
We first perform numerical simulations in a stationary market environment, where the price process is simulated according to the geometric Brownian motion (\ref{eq:St}) with constant $\mu$ and $\sigma$. We take $T = 1$ and $\Delta t = {1}/{252}$, indicating that the MV problem is considered over a 1-year period, with daily rebalancing. For the simulations, we select plausible values of the annualized return and volatility from the sets $\mu \in \{-0.5, -0.3, -0.1, 0.0, 0.1, 0.3, 0.5\}$ and $\sigma \in \{0.1, 0.2, 0.3, 0.4\}$, respectively. These values are commonly used to model a typical stock in simulation studies \citep{hutchinson1994nonparametric}. The annualized interest rate is taken to be $r = 0.02$.

\begin{figure}[htbp]
    \centering
    \includegraphics[width=0.9\textwidth]{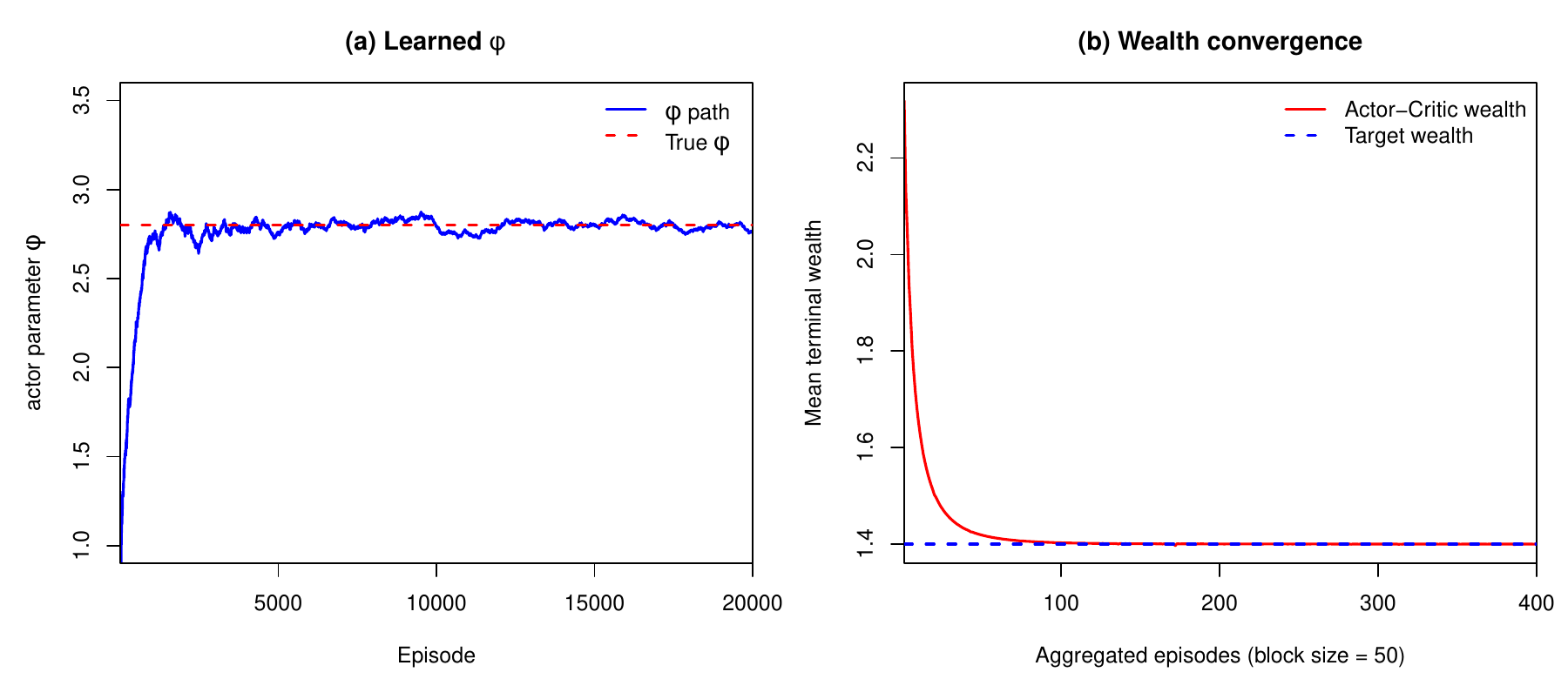}
    \caption{The learned parameters and average reward along one sample trajectory ($\mu=0.3$, $\sigma=0.1$).}
    \label{fig:random}
\end{figure}

\begin{figure}[htbp]
    \centering
    \includegraphics[width=0.9\textwidth]{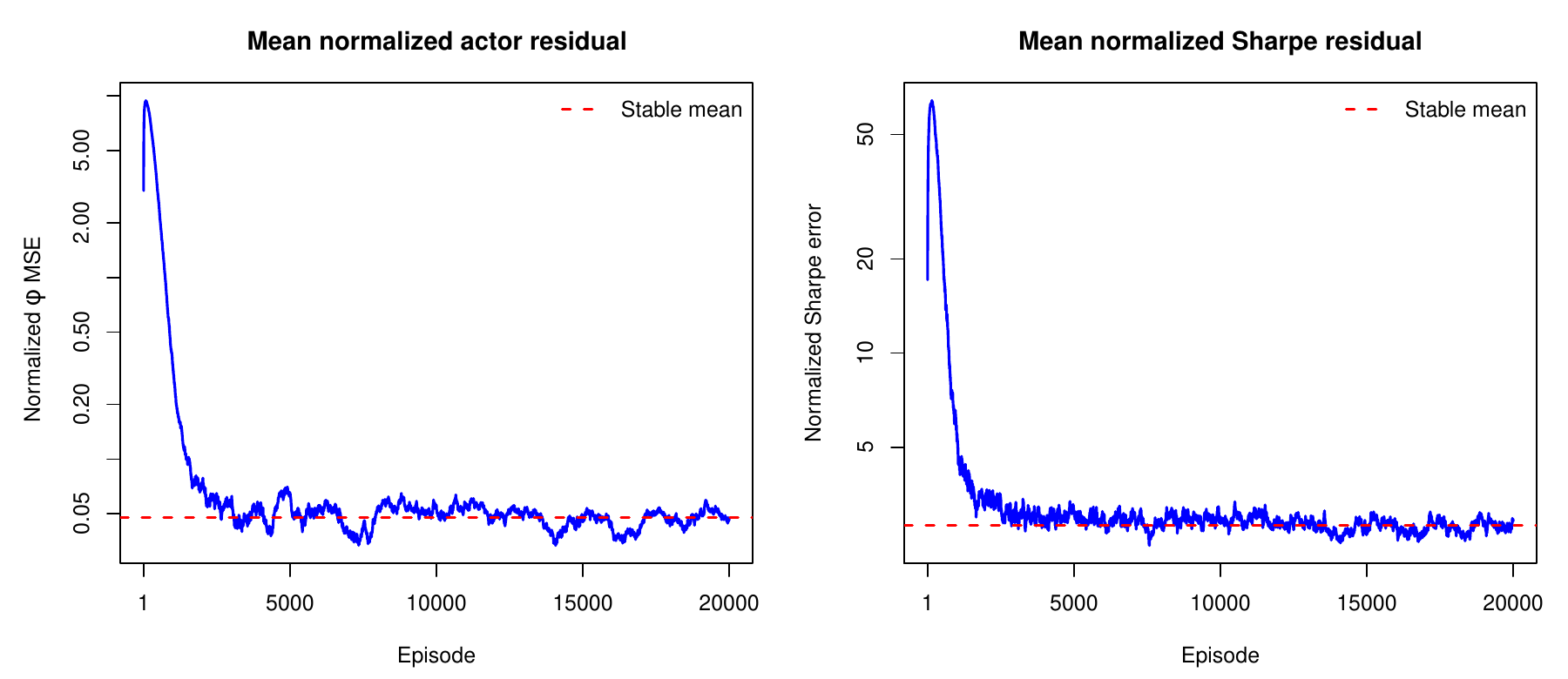}
    \caption{Rate-normalized simulation errors based on $100$ independent
    repetitions ($\mu=0.3$, $\sigma=0.1$).  The red dashed lines indicate the tail averages over the stable
    region $n\geq 10000$.}
    \label{fig:convergence}
\end{figure}

We consider the problem with a $40\%$ annualized target return on the terminal wealth starting from a normalized initial wealth $x_0 = 1$, hence, $z = 1.4$. These model parameters will be fixed for all the simulations considered in this section. For the algorithms, we take the total training episodes $N = 20000$, and take the sample size $ N_w = 10$ for learning the Lagrange multiplier $w$. Across all the simulations in this section,  the learning rates are fixed at $\alpha_\theta=\alpha_\phi=0.1$, \(\alpha_{\sigma}=0.1\) and \(\alpha_w=0.5\), with learning-rate schedules \(\alpha=\gamma_1=\gamma_2=\gamma_3=0.51\). The comparative results under the setting are reported in Table \ref{table:random} and Figure \ref{fig:random}. To assess numerically whether the observed errors are consistent with the
rates established in Theorems~\ref{thm:main-actor} and~\ref{thm:main-sr},
we adopt the rate-normalization plotting approach used in
\citet{banholzer2022rates}. Figure~\ref{fig:convergence} is constructed
from $100$ independent repetitions of our algorithm. In the left panel,
the squared actor error $|\phi_n-\rho^*|^2$ is first averaged across the
repetitions and then divided by the theoretical actor rate
$L_n[(n+1)^{-2\alpha}+(n+1)^{-\gamma_2}]$. In the right panel, the
absolute Sharpe-ratio error $|SR_n-SR_{\mathrm{or}}|$ is averaged in the
same way and divided by the theoretical Sharpe-ratio rate
$L_n[(n+1)^{-\alpha}+(n+1)^{-\gamma_2}]$. 

As shown in Figure~\ref{fig:random}, the learned actor parameter $\phi$ and the mean terminal wealth rapidly approach the true value and the target level, respectively. This indicates that actor learning and multiplier adjustment jointly achieve the desired terminal-wealth constraint. Table~\ref{table:random} shows that Actor--Critic records the highest or tied-highest Sharpe ratio in most parameter settings. Bayesian Filtering typically exhibits higher terminal-wealth variance and a lower Sharpe ratio, while Rolling MLE is markedly less stable, producing very large variances and near-zero Sharpe ratios in many cases. In Figure~\ref{fig:convergence}, the relatively large values at the beginning reflect initialization and the finite-sample transient before the algorithm enters its asymptotic regime. After this initial period, both normalized errors fluctuate around finite levels and show no systematic upward trend. This tail behavior is consistent with the theoretical convergence bounds.

\section{Empirical Analysis}
In this section, we conduct some empirical analysis using real market data. The raw panel was downloaded from Yahoo Finance and contains daily adjusted
closing prices for 499 firms in the contemporaneous S\&P~500 constituent
list, from January~2, 1990 to December~31, 2024. We retain only stocks with no missing adjusted closing price
over the entire 1990--2024 experimental span, leaving 235 eligible stocks.

The split is chronological. The period January~1, 1990--December~31, 2009 is
used exclusively for initialization and training. Initial volatility and
Sharpe-ratio estimates are computed from all daily basket excess returns in
this period. A common training pool is then formed by drawing 20,000 random
starting positions with replacement and constructing one-year windows of 252
return increments (253 price observations). Each learned method resamples
20,000 episodes with replacement from this pool. Consequently, training episodes may overlap, but none contains a
post-2009 observation. The holdout period January~1, 2010--December~31, 2024
is divided chronologically into 15 non-overlapping complete 252-increment
episodes. No holdout return is used
to initialize or bootstrap-train a method.

For each of six experiments, ten stocks are sampled without replacement
within the basket; the same stock may appear in different baskets. Dividing
each stock price by its first observation removes differences in price level,
after which the one-dimensional basket required by the model is obtained by
taking the equally weighted average of the ten normalized price series. The
basket path is rebased to one at the beginning of every episode. Daily simple
returns are then computed from consecutive basket values, and the daily
portion of the 2\% annual risk-free rate, based on 252 trading days per year,
is subtracted to obtain excess returns. These excess returns are the market
observations supplied to all learning rules.

The empirical target is
$z=1.15$, corresponding to 15\% one-year wealth growth. Transaction costs are
omitted. Since the 20,000 training episodes are not completely mutually independent, the filtering learning rate is set to $\alpha=0.2$. All other parameter settings and stochastic methods follow the implementation conventions in Section~5. For each method, we record compounded terminal wealth,
the mean and standard deviation of one-year episode returns, their ratio,
and whether any episode has non-finite or non-positive wealth.
Table~\ref{tab:empirical-aggregate} reports cross-experiment means.  It shows that Actor--Critic achieves the highest Sharpe-like ratio, $1.376$, together with the lowest episode-return standard deviation, $0.122$, among all reported methods. Although buy-and-hold produces the largest compounded terminal wealth, the Actor--Critic result indicates a more favorable risk-adjusted trade-off and no bankruptcy events. Rolling MLE attains the largest mean episode return but combines it with a standard deviation of $4.250$, zero terminal wealth, and a $100\%$ bankruptcy rate, revealing substantial tail instability relative to the Bayesian and Actor--Critic methods.

\begin{table}[htbp]
\centering
\begin{threeparttable}
\caption{Aggregate out-of-sample performance across six ten-stock S\&P~500 baskets.}
\label{tab:empirical-aggregate}
\small
\setlength{\tabcolsep}{4.5pt}
\renewcommand{\arraystretch}{1.15}
\begin{tabularx}{\textwidth}{>{\raggedright\arraybackslash}p{3.2cm} *{5}{>{\centering\arraybackslash}X}}
\toprule
Method
& \makecell{Terminal\\wealth}
& \makecell{Mean episode\\return}
& \makecell{Episode-return\\s.d.}
& \makecell{Sharpe-like\\ratio}
& \makecell{Bankruptcy\\rate} \\
\midrule
Buy-and-hold Basket    & {\bf 19.111} & 0.183 & 0.172 & 1.092 & 0.00 \\
Rolling MLE      & 0.000 & {\bf 0.696} & 4.250 & 0.188 & 1.00 \\
Bayesian Filtering      &  9.492 & 0.161 & 0.165 & 0.996 & 0.00 \\
Actor--Critic       & 7.863 & 0.148 & {\bf 0.122} & {\bf 1.376} & 0.00 \\
\bottomrule
\end{tabularx}
\begin{tablenotes}[flushleft]
\footnotesize
\item Notes: Entries are means across six experiments; the workbook does not
report standard errors. Bold marks the best entry under each criterion. 
\end{tablenotes}
\end{threeparttable}
\end{table}

\section{Conclusions}
This paper develops a novel posterior-guided reinforcement-learning framework for continuous-time mean--variance portfolio optimization. Instead of randomizing the portflio policies with exogenous noise under standard entropy regularization, the method samples a Sharpe ratio from a Gaussian distribution anchored to the Bayesian filter by a KL penalty and maps that draw into the feedback portfolio control. The resulting formulation preserves the common return innovation shared by wealth and beliefs, leading to a  Gaussian stochastic policy with declining degree of exploration guided by the accumulated filtering information of the unknown drift. 

For implementation, normalized martingale estimating equations update the critic and actor, while a separate batch stochastic-approximation recursion enforces the terminal mean constraint. The convergence analysis establishes consistency of the episodic filter, actor, policy center, and Lagrange multiplier. It further yields a last-iterate Sharpe-ratio error of $O_p\left(L_n(n+1)^{-\alpha}+L_n(n+1)^{-\gamma_2}\right)$; in the stationary-market specialization with harmonic filtering and actor schedules, this becomes $O_p(\log\log n/n)$. These results explain how the filtering and actor learning jointly determine the data-driven portfolio decision with dynamically adjusted exploration.

\begin{sidewaystable}[ht!]
\centering
\caption{Comparison of mean terminal wealth, terminal-wealth variance, and Sharpe ratio under Actor--Critic, Bayesian filtering, and Rolling MLE over 100 independent repetitions. Each method is evaluated using 1,000 independent test paths. Entries report the across-repetition mean followed by the standard deviation in parentheses.}
\label{table:random}
\renewcommand{\arraystretch}{1.2}
\normalsize
\setlength{\tabcolsep}{3.45pt}
\begin{tabular}{c *{9}{c}}
    \toprule
    \multirow{2}{*}{$(\mu,\sigma)$}
    & \multicolumn{3}{c}{\textbf{Actor--Critic}}
    & \multicolumn{3}{c}{\textbf{Bayesian Filtering}}
    & \multicolumn{3}{c}{\textbf{Rolling MLE}} \\
    \cmidrule(lr){2-4}\cmidrule(lr){5-7}\cmidrule(lr){8-10}
    & \textbf{Mean} & \textbf{Variance} & \textbf{Sharpe ratio}
    & \textbf{Mean} & \textbf{Variance} & \textbf{Sharpe ratio}
    & \textbf{Mean} & \textbf{Variance} & \textbf{Sharpe ratio} \\
    \midrule
    $(-0.5, 0.1)$ & 1.40(0.00) & 0.00(0.00) & \(\bm{3.20(7.48)\times 10^{14}}\) & 1.40(0.01) & 0.03(0.17) & 3.92(1.05) & 1.40(0.01) & 0.08(0.27) & 5.57(5.57) \\
    $(-0.3, 0.1)$ & 1.40(0.00) & 0.00(0.00) & \(\bm{7160.00(6840.00)}\) & 1.40(0.01) & 0.03(0.04) & 2.65(0.69) & 1.41(0.06) & 4.78(3.24) & 0.21(0.06) \\
    $(-0.1, 0.1)$ & 1.40(0.01) & 0.05(0.03) & \(\bm{1.90(0.40)}\) & 1.40(0.01) & 0.10(0.03) & 1.31(0.18) & 1.19(0.34) & 98.40(38.80) & 0.02(0.03) \\
    $(0, 0.1)$ & 1.21(0.05) & 1.19(0.39) & \(\bm{0.20(0.03)}\) & 1.21(0.06) & 1.20(0.56) & 0.20(0.03) & 1.03(0.38) & 170.00(38.60) & 0.00(0.03) \\
    $(0.1, 0.1)$ & 1.40(0.02) & 0.18(0.03) & \(\bm{0.97(0.09)}\) & 1.40(0.03) & 0.22(0.04) & 0.85(0.08) & 1.09(0.34) & 132.00(31.80) & 0.01(0.03) \\
    $(0.3, 0.1)$ & 1.40(0.00) & 0.00(0.00) & \(\bm{476.00(326.00)}\) & 1.40(0.01) & 0.04(0.04) & 2.43(0.56) & 1.39(0.10) & 11.10(7.41) & 0.13(0.05) \\
    $(0.5, 0.1)$ & 1.40(0.00) & 0.00(0.00) & \(\bm{2.35(5.98)\times 10^{11}}\) & 1.40(0.01) & 0.02(0.03) & 3.85(0.96) & 1.40(0.02) & 0.31(0.72) & 1.79(1.70) \\
    \addlinespace
    $(-0.5, 0.2)$ & 1.40(0.00) & 0.00(0.00) & \(\bm{170.00(120.00)}\) & 1.40(0.01) & 0.05(0.07) & 2.23(0.57) & 1.43(0.13) & 15.20(7.61) & 0.12(0.04) \\
    $(-0.3, 0.2)$ & 1.40(0.00) & 0.02(0.02) & \(\bm{4.32(1.46)}\) & 1.40(0.01) & 0.06(0.03) & 1.69(0.27) & 1.33(0.31) & 65.50(28.60) & 0.04(0.04) \\
    $(-0.1, 0.2)$ & 1.41(0.03) & 0.39(0.06) & \(\bm{0.66(0.06)}\) & 1.41(0.03) & 0.44(0.08) & 0.62(0.06) & 1.02(0.36) & 153.00(33.10) & 0.00(0.03) \\
    $(0, 0.2)$ & 1.07(0.03) & 0.64(0.34) & \(\bm{0.09(0.03)}\) & 1.07(0.04) & 0.64(0.52) & 0.09(0.03) & 1.06(0.39) & 170.00(37.60) & 0.01(0.03) \\
    $(0.1, 0.2)$ & 1.38(0.05) & 0.84(0.16) & \(\bm{0.42(0.04)}\) & 1.37(0.06) & 0.85(0.20) & 0.41(0.04) & 1.10(0.35) & 167.00(40.70) & 0.01(0.03) \\
    $(0.3, 0.2)$ & 1.40(0.01) & 0.02(0.02) & \(\bm{2.83(0.57)}\) & 1.40(0.01) & 0.07(0.03) & 1.53(0.21) & 1.24(0.30) & 79.80(21.20) & 0.03(0.03) \\
    $(0.5, 0.2)$ & 1.40(0.00) & 0.00(0.00) & \(\bm{62.00(32.00)}\) & 1.40(0.01) & 0.05(0.05) & 2.17(0.48) & 1.43(0.15) & 19.10(6.82) & 0.10(0.04) \\
    \addlinespace
    $(-0.5, 0.3)$ & 1.40(0.00) & 0.01(0.02) & \(\bm{6.05(2.32)}\) & 1.40(0.01) & 0.06(0.03) & 1.78(0.30) & 1.33(0.23) & 54.70(15.30) & 0.05(0.03) \\
    $(-0.3, 0.3)$ & 1.40(0.01) & 0.08(0.03) & \(\bm{1.51(0.26)}\) & 1.40(0.01) & 0.13(0.03) & 1.16(0.15) & 1.12(0.32) & 110.00(32.20) & 0.01(0.03) \\
    $(-0.1, 0.3)$ & 1.38(0.04) & 0.86(0.16) & \(\bm{0.42(0.04)}\) & 1.39(0.05) & 0.91(0.22) & 0.41(0.04) & 1.06(0.38) & 156.00(30.10) & 0.00(0.03) \\
    $(0, 0.3)$ & 1.03(0.02) & 0.45(0.27) & 0.05(0.03) & 1.04(0.03) & 0.42(0.43) & \(\bm{0.06(0.03)}\) & 0.95(0.40) & 170.00(30.20) & 0.00(0.03) \\
    $(0.1, 0.3)$ & 1.28(0.05) & 1.11(0.30) & \(\bm{0.27(0.04)}\) & 1.27(0.07) & 1.09(0.39) & 0.27(0.04) & 1.05(0.40) & 178.00(89.40) & 0.00(0.03) \\
    $(0.3, 0.3)$ & 1.40(0.01) & 0.11(0.03) & \(\bm{1.22(0.14)}\) & 1.40(0.02) & 0.16(0.03) & 1.01(0.10) & 1.03(0.33) & 119.00(32.50) & 0.00(0.03) \\
    $(0.5, 0.3)$ & 1.40(0.00) & 0.01(0.01) & \(\bm{4.37(1.11)}\) & 1.40(0.01) & 0.06(0.04) & 1.70(0.27) & 1.31(0.24) & 64.40(24.50) & 0.04(0.03) \\
    \addlinespace
    $(-0.5, 0.4)$ & 1.40(0.01) & 0.04(0.03) & \(\bm{2.29(0.56)}\) & 1.40(0.01) & 0.09(0.03) & 1.42(0.21) & 1.23(0.33) & 89.50(29.60) & 0.03(0.03) \\
    $(-0.3, 0.4)$ & 1.40(0.02) & 0.18(0.04) & {\bf 0.96(0.11)} & 1.40(0.02) & 0.23(0.04) & 0.85(0.09) & 1.11(0.41) & 129.00(28.00) & 0.01(0.04) \\
    $(-0.1, 0.4)$ & 1.32(0.05) & 1.14(0.27) & \(\bm{0.30(0.03)}\) & 1.32(0.06) & 1.17(0.37) & 0.30(0.03) & 0.98(0.44) & 173.00(45.70) & 0.00(0.03) \\
    $(0, 0.4)$ & 1.02(0.02) & 0.37(0.23) & 0.03(0.03) & 1.02(0.02) & 0.33(0.37) & \(\bm{0.04(0.03)}\) & 0.95(0.38) & 171.00(36.40) & 0.00(0.03) \\
    $(0.1, 0.4)$ & 1.20(0.05) & 1.03(0.38) & \(\bm{0.20(0.03)}\) & 1.19(0.06) & 0.97(0.45) & 0.20(0.04) & 1.06(0.45) & 174.00(39.70) & 0.00(0.03) \\
    $(0.3, 0.4)$ & 1.40(0.02) & 0.25(0.04) & \(\bm{0.81(0.07)}\) & 1.40(0.03) & 0.30(0.05) & 0.74(0.06) & 1.13(0.40) & 140.00(41.20) & 0.01(0.03) \\
    $(0.5, 0.4)$ & 1.40(0.01) & 0.05(0.02) & \(\bm{1.93(0.30)}\) & 1.40(0.01) & 0.10(0.03) & 1.32(0.16) & 1.23(0.33) & 102.00(29.20) & 0.02(0.03) \\
    \bottomrule
\end{tabular}
\end{sidewaystable}

\ \\
\ \\
\textbf{Acknowledgements}: Xiang Yu is supported by the Hong Kong RGC General Research Fund (GRF) under grants no. 15211524 and no. 15214125. \\

\clearpage
\bibliographystyle{abbrvnat}
\bibliography{ref}

\clearpage
\appendix

\section{Verification and Admissibility Proofs}

\subsection*{Proof of Theorem~\ref{thm:recursive-verification}}
\label{app:verification-theorem-proof}

\begin{proof}
Fix $(t,x,m)$ with $x\ne w$ and an arbitrary
$q\in\mathcal Q(t,x,m)$. For brevity, all the derivatives of $v$ below
are evaluated at $(s,X_s^\pi,m_s)$. Under \eqref{eq:macroSDE}, the
controlled generator applied to $v$ is
\begin{align}
\mathcal L^qv={}&-\rho^*M_s(X_s^\pi-w)v_x
+\frac12(X_s^\pi-w)^2(M_s^2+V_s)v_{xx}
+R_s(\rho^*-m_s)v_m+\frac12R_s^2v_{mm}
\notag\\
&
-M_sR_s(X_s^\pi-w)v_{xm}.
\label{eq:appendix-controlled-generator}
\end{align}
Since the infimum of the Hamiltonian is bounded above by its value at the specific density $q_s$, the HJB equation implies
\begin{equation}
\Bigl(v_t+\mathcal L^qv
+\lambda_sD_{\mathrm{KL}}
\bigl(q_s\Vert N(m_s,R_s)\bigr)v\Bigr)
(s,X_s^\pi,m_s)\ge0.
\label{eq:appendix-hjb-inequality}
\end{equation}
Recall the integrating factor $\Lambda_s^q$ defined in Definition~\ref{def:admissible-distributional-control}. Then
\[
d\Lambda_s^q
=\lambda_sD_{\mathrm{KL}}
\bigl(q_s\Vert N(m_s,R_s)\bigr)\Lambda_s^qds,
\quad \Lambda_t^q=1.
\]
It\^o's formula yields
\begin{align}
dv(s,X_s^\pi,m_s)
={}&(v_t+\mathcal L^qv)(s,X_s^\pi,m_s)ds
+\Bigl[R_sv_m-M_s(X_s^\pi-w)v_x\Bigr]
(s,X_s^\pi,m_s)dW_s
\notag\\
&-\sqrt{V_s}(X_s^\pi-w)v_x(s,X_s^\pi,m_s)dB_s.
\label{eq:appendix-ito-v}
\end{align}
Because this integrating factor $d\Lambda_s^q$ has finite variation, its quadratic
covariation with $v(s,X_s^\pi,m_s)$ is zero. The product rule
gives
\begin{align}
d\left[
\Lambda_s^qv(s,X_s^\pi,m_s)\right]
={}&\Lambda_s^q
\Bigl(v_t+\mathcal L^qv
+\lambda_sD_{\mathrm{KL}}
\bigl(q_s\Vert N(m_s,R_s)\bigr)v\Bigr)
(s,X_s^\pi,m_s)ds
\notag\\
+&\Lambda_s^q
\Bigl[R_sv_m-M_s(X_s^\pi-w)v_x\Bigr]
(s,X_s^\pi,m_s)dW_s
\notag\\
-&\Lambda_s^q
\sqrt{V_s}(X_s^\pi-w)v_x
(s,X_s^\pi,m_s)dB_s.
\label{eq:appendix-integrating-factor-ito}
\end{align}
Choose the explicit localization sequence
\begin{align}
\tau_n:={}&\left(T-\frac{T-t}{n}\right)
\wedge\inf\bigg\{s\ge t:\int_t^s
(\Lambda_v^q)^2\bigg(
\Bigl[R_vv_m(v,X_v^\pi,m_v)
-M_v(X_v^\pi-w)v_x(v,X_v^\pi,m_v)\Bigr]^2
\notag\\
&\hspace{31mm}
+V_v(X_v^\pi-w)^2v_x^2(v,X_v^\pi,m_v)
\bigg)dv\ge n\bigg\},
\label{eq:appendix-verification-stopping-times}
\end{align}
with the convention $\inf\varnothing=\infty$. The stochastic
integrals in \eqref{eq:appendix-integrating-factor-ito}, stopped at
$\tau_n$, are square-integrable martingales, because their joint quadratic
variation is bounded by $n$. Moreover, $\tau_n\uparrow T$ almost surely.
Indeed, the stochastic-integral coefficients are locally square-integrable
on every compact subinterval of $[t,T)$. Also, by
\eqref{eq:wealth-explicit-admissibility}, $X_s^\pi-w$ does not reach zero
when $x\ne w$. Integrating \eqref{eq:appendix-integrating-factor-ito} from $t$ to
$\tau_n$ and taking $\E$, we get
\begin{align}
&\E\left[
\Lambda_{\tau_n}^qv(\tau_n,X_{\tau_n}^\pi,m_{\tau_n})
\right]-v(t,x,m)
\notag\\
&=
\E\int_t^{\tau_n}\Lambda_s^q
\Bigl(v_t+\mathcal L^qv
+\lambda_sD_{\mathrm{KL}}
\bigl(q_s\Vert N(m_s,R_s)\bigr)v\Bigr)
(s,X_s^\pi,m_s)ds
\ge0,
\label{eq:appendix-stopped-verification-inequality}
\end{align}
where the last inequality follows from
\eqref{eq:appendix-hjb-inequality}. Hence, for every $n$,
\begin{equation}
v(t,x,m)\le
\E\left[
\Lambda_{\tau_n}^qv(\tau_n,X_{\tau_n}^\pi,m_{\tau_n})
\right].
\label{eq:appendix-stopped-upper-bound}
\end{equation}
As $\tau_n\uparrow T$, continuity of $v$ and the terminal condition imply
\begin{align*}
\Lambda_{\tau_n}^qv(\tau_n,X_{\tau_n}^\pi,m_{\tau_n})
\longrightarrow
\Lambda_T^qv(T,X_T^\pi,m_T)=\Lambda_T^q(X_T^\pi-w)^2
\quad\text{almost surely}.
\end{align*}
Condition~(vi) of
Definition~\ref{def:admissible-distributional-control} asserts precisely that
the random variables on the left, as $\tau_n$ ranges over stopping times,
form a uniformly integrable family. Therefore, by Vitali's convergence
theorem, the almost-sure convergence above and condition~(vi) imply
$L^1$ convergence. Hence,
\begin{align}
&\lim_{n\to\infty}\E\left[
\Lambda_{\tau_n}^qv(\tau_n,X_{\tau_n}^\pi,m_{\tau_n})
\right]
=\E\left[\Lambda_T^q(X_T^\pi-w)^2\right].
\label{eq:appendix-ui-limit}
\end{align}
It follows from \eqref{eq:appendix-stopped-upper-bound} and
\eqref{eq:appendix-ui-limit} that
\begin{equation}
v(t,x,m)\le
\E\left[\Lambda_T^q(X_T^\pi-w)^2\right].
\label{eq:appendix-verification-terminal-bound}
\end{equation}

It remains to identify the expectation in
\eqref{eq:appendix-verification-terminal-bound} with the recursively
defined policy value. Condition~(v) gives directly
\[
\E\left[\Lambda_T^q(X_T^\pi-w)^2\right]<\infty.
\]
The $L^1$ representation \eqref{eq:l1-policy-value-representation} then
shows that $\Lambda^qJ^q$ is the uniformly integrable martingale with this
terminal value. Consequently,
\begin{align}
\Lambda_s^qJ_s^q
=\E\left[
\Lambda_T^q(X_T^\pi-w)^2\mid\mathcal F_s\right],
\quad
J_t^q
=\E\left[\Lambda_T^q(X_T^\pi-w)^2\right].
\label{eq:appendix-bsde-terminal-representation}
\end{align}
It holds by \eqref{eq:appendix-verification-terminal-bound} and
\eqref{eq:appendix-bsde-terminal-representation} that
\[
v(t,x,m)\le J^q(t,x,m)
\quad\text{for every }q\in\mathcal Q(t,x,m).
\]
For $q=q^*$, the drift in
\eqref{eq:appendix-integrating-factor-ito} is identically zero. Therefore
\eqref{eq:appendix-stopped-verification-inequality} is an equality for
every $n$, and the same uniformly integrable limit gives
$v(t,x,m)=J^{q^*}(t,x,m)$. 
It follows that
\[
v(t,x,m)\le\inf_{q\in\mathcal Q(t,x,m)}J^q(t,x,m)
\le J^{q^*}(t,x,m)=v(t,x,m),
\]
which proves $v=J$ and the optimality of $q^*$.

Finally, if $x=w$, then \eqref{eq:macroSDE} implies
$X_s^\pi-w=0$ for all $s\in[t,T]$. The terminal variable in
\eqref{eq:l1-policy-value-representation} is therefore zero, and that
representation gives $J^q(t,w,m)=0$ for every admissible $q$. By the continuous
extension stipulated in the theorem, $v(t,w,m)=0$ as well. Hence the same
identity holds.
\end{proof}

\subsection*{Proof of Proposition~\ref{prop:closed-form-minimizer}}
\label{app:closed-form-admissibility}

\begin{proof}
Fix $(t,x,m)$ with $t<T$. We verify the six items in
Definition~\ref{def:admissible-distributional-control} directly on the
given finite horizon.

\paragraph{Step 1: the elementary admissibility conditions.}
Because $V_s^*={cR_s}/{(c+2)}>0$ on $[t,T]$, $q_s^*=N(M_s^*,V_s^*)$ is a probability density. To make the
measurability and moment claims explicit, divide the filtering equation by
$R_s$ and use $d(R_s^{-1})=ds$:
\[
d\left(\frac{m_s-\rho^*}{R_s}\right)=dW_s.
\]
Hence
\begin{equation}
m_s=\rho^*+\frac{R_s}{R_t}(m-\rho^*)+R_s(W_s-W_t),
\label{eq:appendix-filter-explicit}
\end{equation}
which is adapted and Gaussian. In view that
$M_s^*=a_sm_s+b_s\rho^*$ and $a_s,b_s,V_s^*$ are bounded deterministic
functions on $[t,T]$, it holds that $q^*$ is progressively measurable and
$\E\int_t^T\bigl((M_s^*)^2+V_s^*\bigr)ds<\infty$. This verifies conditions~(i)--(iii). The same integrability ensures that the each term
in \eqref{eq:wealth-explicit-admissibility} is well defined, and the
existence and pathwise uniqueness of the linear wealth SDE also follow. 

\paragraph{Step 2: a global auxiliary estimate for the wealth moments.}
We use the time scale already employed in the learning analysis
$\tau_s:=\log(R_s/R_T)$ and write
$p_c(\tau_s):=R_sP(s)$. Since $d\tau_s/ds=-R_s$, the Riccati equation for
$P$ gives
\begin{equation}
p_c'
=\frac{(c-2)p_c^2-(5c+2)p_c+2c}{c+2},
\quad p_c(0)=0.
\label{eq:appendix-p-phase-line}
\end{equation}
For $c\ne2$, the right-hand side of
\eqref{eq:appendix-p-phase-line} is positive at zero and equals $-2$ at
one. Its first positive root is $r_-$, so phase-line analysis gives
$0\le p_c(\tau)<r_-<1$ and $p_c'(\tau)\ge0$. For $c=2$, the same facts
follow from $p_2(\tau)=(1-e^{-3\tau})/3$. Consequently,
\begin{equation}
\frac{c}{c+2}\le a_s=\frac{c+2p_c(\tau_s)}{c+2}<1,
\quad \frac{d a_s}{d\tau_s}\ge0.
\label{eq:appendix-a-range}
\end{equation}

For later use, let $G_X:=G^{(T)}$ be the maximal backward solution of
\begin{equation}
\frac12\dot G_X
-R_s(1+4a_s)G_X
+\frac12R_s^2G_X^2
+6a_s^2=0,
\quad G_X(T)=0.
\label{eq:riccati-wealth-integrability}
\end{equation}
Writing $\widehat G_X(\tau_s):=R_sG_X(s)$ transforms this equation into
\begin{equation}
\widehat G_X'
=\widehat G_X^2-(1+8a_s)\widehat G_X+12a_s^2,
\quad \widehat G_X(0)=0.
\label{eq:appendix-forward-riccati}
\end{equation}
The nonnegative half-line is invariant because the right-hand side of
\eqref{eq:appendix-forward-riccati} equals $12a_s^2\ge0$ at zero. Viewed
as a quadratic polynomial in $\widehat G_X$, its smaller positive zero is $h_s:=\frac{1+8a_s-\sqrt{1+16a_s+16a_s^2}}{2}$. This expression is increasing in $a_s$ as
$\sqrt{1+16a_s+16a_s^2}>1+2a_s$, and hence is nondecreasing in
$\tau_s$ by \eqref{eq:appendix-a-range}. Its defining identity gives
\[
\frac{dh_s}{d\tau_s}
\ge h_s^2-(1+8a_s)h_s+12a_s^2=0.
\]
Thus $h_s$ is a supersolution of \eqref{eq:appendix-forward-riccati}.
As $\widehat G_X(0)=0<h_T$, scalar comparison yields
\begin{equation*}
0\le\widehat G_X(\tau_s)
\le h_s
=\frac{1+8a_s-\sqrt{1+16a_s+16a_s^2}}{2}
<\infty
\end{equation*}
on every finite interval. As a result, $G_X$ is finite on $[t,T]$ for every
$c>0$ and every finite horizon, without any restriction on
$\log(R_t/R_T)$.
\paragraph{Step 3: the wealth moments.}
Set $Y_s^\ast=X_s^{\pi^*}-w$. Under the candidate policy,
\begin{align}
dY_s^\ast={}&-\rho^*M_s^*Y_s^\ast ds
-M_s^*Y_s^\ast\,dW_s-\sqrt{V_s^*}Y_s^\ast\,dB_s,
\notag\\
dm_s={}&R_s(\rho^*-m_s)ds+R_s\,dW_s.
\label{eq:appendix-joint-state}
\end{align}
The shared Brownian motion gives
$d\langle Y^\ast,m\rangle_s=-M_s^*R_sY_s^\ast ds$. For each fixed
$u\in[t,T]$, consider the conditional fourth-moment function
\[
(s,y,z)\longmapsto
\E\bigl[|Y_u|^4\mid Y_s^\ast=y,\ m_s=z\bigr],
\quad t\le s\le u.
\]
Its exponential-quadratic candidate is
\[
U^{(u)}(s,y,z)=y^4
\exp\left\{\frac12G^{(u)}(s)z^2+H^{(u)}(s)z+I^{(u)}(s)\right\},
\quad U^{(u)}(u,y,z)=y^4.
\]
The multiplicative wealth dynamics make the fourth moment homogeneous of
degree four in $y$, while the affine dependence of $M_s^*$ on the
Gaussian filter gives the exponential-quadratic dependence on $z$. The
superscript records the terminal date. Before the required integrability is established,
$U^{(u)}$ is used only as a positive auxiliary function for the stopped
argument below rather than being assumed a priori to equal the displayed
conditional expectation. The drift in It\^o's formula is obtained from
\begin{align*}
U^{(u)}_s
-\rho^*M_s^*yU^{(u)}_y+R_s(\rho^*-z)U^{(u)}_z
+\frac12y^2\bigl((M_s^*)^2+V_s^*\bigr)U^{(u)}_{yy}
+\frac12R_s^2U^{(u)}_{zz}-M_s^*R_syU^{(u)}_{yz}.
\end{align*}
Using $M_s^*=a_sz+b_s\rho^*$, the coefficient of $z^2U^{(u)}$ vanishes when
\[
\frac12\dot G^{(u)}-R_s(1+4a_s)G^{(u)}
+\frac12R_s^2\bigl(G^{(u)}\bigr)^2+6a_s^2=0,
\quad G^{(u)}(u)=0.
\]
For $u=T$, this is \eqref{eq:riccati-wealth-integrability}, so
$G^{(T)}=G_X$. For an earlier terminal date, scalar comparison after the
time change in Step~2 shows that
\[
0\le G^{(u)}(s)\le G_X(s),\quad t\le s\le u.
\]
Thus $G^{(u)}$ is finite, uniformly over $u\in[t,T]$. Equating the
coefficients of $zU^{(u)}$ and $U^{(u)}$ to zero gives, respectively,
\begin{align*}
\dot H^{(u)}
+\bigl[R_s^2G^{(u)}-R_s(1+4a_s)\bigr]H^{(u)}
+\rho^*\bigl[-4a_s+12a_sb_s+R_sG^{(u)}(1-4b_s)\bigr]
={}&0,
\quad H^{(u)}(u)=0,\\
\dot I^{(u)}
+6V_s^*+\frac12R_s^2\bigl(G^{(u)}+(H^{(u)})^2\bigr)
+R_s\rho^*(1-4b_s)H^{(u)}
+(\rho^*)^2(6b_s^2-4b_s)
={}&0,
\quad I^{(u)}(u)=0.
\end{align*}
Once $G^{(u)}$ is fixed, the first equation is linear in $H^{(u)}$, and
the second determines $I^{(u)}$ after $H^{(u)}$ is known. The bound by
$G_X$ and the boundedness
of the deterministic coefficients therefore imply, by the standard
linear-ODE estimates, that $H^{(u)}$ and $I^{(u)}$ are also uniformly
bounded over $u\in[t,T]$. Consequently, $U^{(u)}$ solves the backward
Kolmogorov equation for the joint process $(Y^\ast,m)$ with terminal value
$y^4$.

For completeness, stop the joint state at
$\sigma_n^{(u)}:=\inf\{s\ge t:|Y_s^\ast|+|m_s|\ge n\}\wedge u$. With the three coefficient equations just described, the stopped process
$U^{(u)}(s\wedge\sigma_n^{(u)},Y_{s\wedge\sigma_n^{(u)}}^\ast,
m_{s\wedge\sigma_n^{(u)}})$ is a true
martingale because all its stochastic-integral coefficients are bounded.
Consequently,
\[
U^{(u)}(t,x-w,m)
\ge
\E\left[\mathbf 1_{\{\sigma_n^{(u)}=u\}}|Y_u^\ast|^4\right].
\]
Continuity of $(Y^\ast,m)$ gives $\sigma_n^{(u)}\uparrow u$ almost surely, and Fatou's
lemma yields $\E|Y_u^\ast|^4<\infty$. The uniform bounds on
$G^{(u)},H^{(u)},I^{(u)}$ further give
$\sup_{t\le s\le T}\E|Y_s^\ast|^4<\infty$. Moreover, \eqref{eq:appendix-filter-explicit} implies
$\E [\sup_{t\le s\le T}|m_s|^2]<\infty$ and
$\sup_s\E|M_s^*|^4<\infty$. Applying the
BDG inequality to
\eqref{eq:appendix-joint-state}, followed by Cauchy--Schwarz, gives
\begin{align*}
\E\sup_{t\le s\le T}|Y_s^\ast|^2
\le C\bigg(&|x-w|^2
+\int_t^T
\bigl(\E|M_s^*|^4\bigr)^{1/2}
\bigl(\E|Y_s^\ast|^4\bigr)^{1/2}ds+\int_t^TV_s^*\E|Y_s^\ast|^2ds\bigg)<\infty.
\end{align*}
In view of $X_s^{\pi^*}=Y_s^\ast+w$, these estimates verify the condition~(iv).

\paragraph{Step 4: the weighted terminal loss and uniform integrability.}
The Gaussian KL formula, the affine representation
$M_s^*=a_sm_s+b_s\rho^*$, and continuity of $m$ show first that
\[
\int_t^T\lambda_sD_{\mathrm{KL}}
\bigl(q_s^*\Vert N(m_s,R_s)\bigr)ds<\infty
\quad\text{almost surely}.
\]
Let $v$ be the closed-form candidate in
\eqref{eq:unified_value}. For $x\ne w$, normalize the process used in the
verification theorem as
\[
Z_s:=
\frac{\exp\left\{\int_t^s\lambda_vD_{\mathrm{KL}}
(q_v^*\Vert N(m_v,R_v))dv\right\}
v(s,X_s^{\pi^*},m_s)}
{v(t,x,m)}.
\]
Indeed, the closed-form expression for $v$
gives, along the state process,
\[
\frac{v_x}{v}=\frac{2}{Y_s^\ast},
\quad
\frac{v_m}{v}=P(s)m_s+\rho^*\overline Q(s).
\]
Because $q^*$ attains the HJB minimum, the HJB equation implies that the
drift of $v(s,X_s^{\pi^*},m_s)$ is
$-\lambda_sD_{\mathrm{KL}}(q_s^*\Vert N(m_s,R_s))v\,ds$. Its $dW_s$ and
$dB_s$ coefficients are, respectively,
$-M_s^*Y_s^\ast v_x+R_sv_m$ and $-\sqrt{V_s^*}Y_s^\ast v_x$. Hence the product rule
shows that the exponential integrating factor in $Z$ cancels the drift.
Using $M_s^*=a_sm_s+b_s\rho^*$ then gives
\begin{equation}
\frac{dZ_s}{Z_s}
=\left[
\bigl(R_sP(s)-2a_s\bigr)m_s
+\rho^*\bigl(R_s\overline Q(s)-2b_s\bigr)
\right]dW_s
-2\sqrt{V_s^*}\,dB_s.
\label{eq:appendix-candidate-density}
\end{equation}
The coefficients are locally square integrable, so this equation makes
$Z$ a continuous local martingale. Moreover, the multiplicative equation
for $Y^\ast$ preserves $Y_s^\ast\ne0$ when $x\ne w$. Hence, $v>0$ and $Z$ is
strictly positive. All deterministic coefficients in \eqref{eq:appendix-candidate-density}
are bounded on $[t,T]$. Hence its squared diffusion coefficient is bounded
by $C(1+m_s^2)$.

We now prove that $Z$ is a true martingale. Set
$\sigma_n:=\inf\{s\ge t:|m_s|\ge n\}\wedge T$. The stopped coefficients
in \eqref{eq:appendix-candidate-density} are bounded, so
$Z_{s\wedge\sigma_n}$ is a true martingale and
$\E[Z_{\sigma_n}]=1$. Hence
$\left.d \mathbb Q_n/d \mathbb P\right|_{\mathcal F_T}=Z_{\sigma_n}$ defines a
probability measure.
The two-dimensional Girsanov theorem states that subtracting the time
integral of each stopped diffusion coefficient in
\eqref{eq:appendix-candidate-density} gives Brownian motions under
$\mathbb Q_n$. In particular,
\[
W_s-\int_t^{s\wedge\sigma_n}
\left[
\bigl(R_vP(v)-2a_v\bigr)m_v
+\rho^*\bigl(R_v\overline Q(v)-2b_v\bigr)
\right]dv
\]
is a Brownian motion under $\mathbb Q_n$. Substituting the
corresponding decomposition of $dW_s$ into the filtering equation shows
that, under $\mathbb Q_n$, the drift of $m_s$ is
\[
R_s\left[
\rho^*-m_s
+\mathbf 1_{\{s\le\sigma_n\}}
\left(
\bigl(R_sP(s)-2a_s\bigr)m_s
+\rho^*\bigl(R_s\overline Q(s)-2b_s\bigr)
\right)
\right].
\]
This change of measure is used below through
$\E[Z_{\sigma_n}\log Z_{\sigma_n}]
=\E^{\mathbb Q_n}[\log Z_{\sigma_n}]$.
Denote the displayed $\mathbb Q_n$-Brownian motion by
$W^{\mathbb Q_n}$. The boundedness of $R,P,\overline Q,a,b$ gives, with $C$ independent of $n$,
\[
\left|R_s\left[\rho^*-m_s+
\mathbf 1_{\{s\le\sigma_n\}}
\left(\bigl(R_sP(s)-2a_s\bigr)m_s
+\rho^*\bigl(R_s\overline Q(s)-2b_s\bigr)\right)\right]\right|
\le C(1+|m_s|),
\quad |R_s|\le C.
\]
For $s\in[t,T]$, It\^o's formula applied to $|m_s|^2$ gives
\begin{align*}
|m_s|^2
={}&|m|^2+2\int_t^sm_vR_v\left[\rho^*-m_v+
\mathbf 1_{\{v\le\sigma_n\}}
\left(\bigl(R_vP(v)-2a_v\bigr)m_v
+\rho^*\bigl(R_v\overline Q(v)-2b_v\bigr)\right)\right]dv
\notag\\
&+\int_t^sR_v^2dv+2\int_t^sR_vm_v\,dW_v^{\mathbb Q_n}.
\end{align*}
Taking the supremum and expectation and using
$2|m_v|(1+|m_v|)\le C(1+|m_v|^2)$ to control the drift gives
\begin{align*}
\E^{\mathbb Q_n}\sup_{t\le u\le s}|m_u|^2
\le{}&|m|^2+C\int_t^s\left(
1+\E^{\mathbb Q_n}\sup_{t\le u\le v}|m_u|^2\right)dv
+2\E^{\mathbb Q_n}\sup_{t\le u\le s}
\left|\int_t^uR_vm_v\,dW_v^{\mathbb Q_n}\right|.
\end{align*}
For the last term, the BDG inequality \citep[Theorem~4.1]{revuz1999continuous} first replaces
the martingale maximum by its quadratic variation; boundedness of $R$ and
Cauchy--Schwarz give the second inequality below, and Young's inequality
implies the last:
\begin{align*}
\E^{\mathbb Q_n}\sup_{t\le u\le s}
\left|\int_t^uR_vm_v\,dW_v^{\mathbb Q_n}\right|
&\le C\E^{\mathbb Q_n}
\left(\int_t^sR_v^2|m_v|^2dv\right)^{1/2}\\
&\le C\left(
\E^{\mathbb Q_n}\sup_{t\le u\le s}|m_u|^2\right)^{1/2}\le \frac12\E^{\mathbb Q_n}\sup_{t\le u\le s}|m_u|^2+C.
\end{align*}
After absorbing the last expectation into the left-hand side, we obtain
\[
\E^{\mathbb Q_n}\sup_{t\le u\le s}|m_u|^2
\le C(1+|m|^2)+C\int_t^s
\E^{\mathbb Q_n}\sup_{t\le u\le v}|m_u|^2dv.
\]
Gronwall's inequality then yields, with a constant independent of
$n$,
\begin{equation*}
\sup_n\E^{\mathbb Q_n}
\left[\sup_{t\le s\le T}|m_s|^2\right]
\le C(1+|m|^2)e^{C(T-t)}<\infty.
\end{equation*}
Using \eqref{eq:appendix-candidate-density} and Girsanov's theorem once
more, we get that
\begin{align}
\E\left[Z_{\sigma_n}\log Z_{\sigma_n}\right]
={}&\frac12\E^{\mathbb Q_n}\int_t^{\sigma_n}
\bigg(
\left[
\bigl(R_vP(v)-2a_v\bigr)m_v
+\rho^*\bigl(R_v\overline Q(v)-2b_v\bigr)
\right]^2
+4V_v^*
\bigg)dv
\notag\\
\le{}&C\left(
1+\E^{\mathbb Q_n}\int_t^T|m_v|^2dv
\right)\le C.
\label{eq:appendix-benes-entropy-bound}
\end{align}
As $x\log x\ge-e^{-1}$ for $x\ge0$,
it follows from \eqref{eq:appendix-benes-entropy-bound} that, with
$\log^+x:=\max\{\log x,0\}$,
\[
\sup_n\E\left[Z_{\sigma_n}\log^+Z_{\sigma_n}\right]
\le \sup_n\E\left[Z_{\sigma_n}\log Z_{\sigma_n}\right]+e^{-1}
<\infty.
\]
As $x\log^+x$ is increasing, convex, and grows faster than linearly, the
de la Vall\'ee--Poussin criterion \citep{de1915integrale} makes
$\{Z_{\sigma_n}\}_{n\ge1}$ uniformly integrable. Continuity ensures 
$\sigma_n\uparrow T$ and $Z_{\sigma_n}\to Z_T$ almost surely; uniform
integrability then gives convergence in $L^1$ and $\E[Z_T]=1$. Hence the
positive local martingale $Z$ is a true martingale on $[t,T]$. Fatou's
lemma and conditional Jensen's inequality further give, for every
$[t,T]$-valued stopping time $\tau$,
\[
\E\left[Z_\tau\log^+Z_\tau\right]
\le \E\left[Z_T\log^+Z_T\right]
\le \liminf_n\E\left[Z_{\sigma_n}\log^+Z_{\sigma_n}\right]<\infty.
\]
The same criterion then shows that
$\{Z_\tau:\tau\le T\}$ is uniformly integrable.

By virtue of $v(T,X_T^{\pi^*},m_T)=(X_T^{\pi^*}-w)^2$, we also obtain
\[
\E\left[
\exp\left\{\int_t^T\lambda_sD_{\mathrm{KL}}
(q_s^*\Vert N(m_s,R_s))ds\right\}
|X_T^{\pi^*}-w|^2
\right]
=v(t,x,m)<\infty.
\]
This proves condition~(v). Moreover, by the definition of $Z$, for every
$[t,T]$-valued stopping time $\tau$,
\[
\exp\left\{\int_t^\tau\lambda_sD_{\mathrm{KL}}
(q_s^*\Vert N(m_s,R_s))ds\right\}
v(\tau,X_\tau^{\pi^*},m_\tau)
=v(t,x,m)Z_\tau.
\]
As $v(t,x,m)<\infty$ is deterministic, multiplication by this constant
preserves uniform integrability. Thus the uniform integrability of
$\{Z_\tau:\tau\le T\}$ proves condition~(vi). If $x=w$, the wealth
deviation and the weighted candidate process are identically zero, and
the same conclusion follows directly.

We have  verified conditions~(i)--(vi) on every finite horizon, and
hence $q^*\in\mathcal Q(t,x,m)$ without an additional restriction.
Theorem~\ref{thm:recursive-verification} then shows that $q^*$ is optimal
and that the closed-form candidate is the value function.
\end{proof}

\section{Convergence Proofs and Auxiliary Estimates}\label{app:convergence-proofs}

\subsection*{Basic Model Setup}

In the conditional derivation of the conditional mean actor direction, the Lagrange multiplier $w$ is first treated as a given constant. Its convergence is proved after the actor analysis.
The policy value, filtering, and Sharpe ratio use the macroscopic relaxed-control dynamics in Equation \eqref{eq:macroSDE}.
To remain consistent with the discrete actor update \eqref{eq:main-sampled-wealth-dynamics}, the actor--critic convergence part uses the sampled-action microscopic dynamics before averaging over actions.

For the $n$th episode, let $\F_n$ denote the information available before the start of that episode.
The original market excess-return process is
\[
dS_t^{(n)}/S_t^{(n)}
=
\sigstar\left(\rhostar dt+dW_t^{(n)}\right).
\]
When $\sigma^\ast$ is unknown, we first use the historical data available before the episode to construct an $\F_n$-measurable volatility estimator $\widehat\sigma_n$.
For a nonempty closed convex set $C\subset\mathbb R$, denote the Euclidean projection by
\[
\Pi_C(z):=\arg\min_{y\in C}|z-y|.
\]
Correspondingly, define the distance from a point $z$ to $C$ by
\[
\operatorname{dist}(z,C)
:=\inf_{y\in C}|z-y|
=|z-\Pi_C(z)|.
\]
Take constants $0<\underline\sigma<\sigstar<\overline\sigma<\infty$, and let the volatility projection set be
$\Sigma:=[\underline\sigma,\overline\sigma]$.
After the $j$th episode ends, define the nonnegative volatility statistic using the time-averaged continuous quadratic variation
\begin{equation*}
U_j
:=
\frac1T\int_0^T d\langle \ln S^{(j)}\rangle_t.
\end{equation*}
The corresponding projected stochastic-approximation recursion is written as
\begin{equation}
\widehat\sigma_{n+1}
=
\Pi_{\Sigma}
\left[
\left(
\frac1{n+1}\sum_{j=0}^{n}U_j
\right)^{1/2}
\right].
\label{eq:sigma-update}
\end{equation}
For the initial episode, set $\chi_0=1$ by convention. For $n\ge1$, let $\chi_n=\sigma^*/\widehat\sigma_n$ denote the ratio between the true and estimated standard deviations. Therefore, the standardized market excess-return observation satisfies
\begin{align}
d\widetilde S_t^{(n)}
&=\chi_n\left(\rhostar\,dt+dW_t^{(n)}\right),
\notag\\
dm_t^{(n)}
&=R_t^{(n)}
\left[
d\widetilde S_t^{(n)}-m_t^{(n)}dt
\right],\label{eq:filter-m-n}\\
dR_t^{(n)}
&=-\bigl(R_t^{(n)}\bigr)^2dt,
\quad 0\le t\le T.\notag
\end{align}
The Brownian increments generated in different episodes are mutually independent conditioned on past information.
Across episodes, the posterior mean and filtering variance are initialized according to \eqref{eq:episode-decay}.

For the continuous-time sampled-action model, writing $Y_t^{(n)}:=\widetilde X_t^{\pi,n}-w$, the wealth dynamics is

\begin{equation}
dY_t^{(n)}
=-
\chi_n\rhostar\rho_t^{(n)}Y_t^{(n)}dt
-
\chi_n\rho_t^{(n)}Y_t^{(n)}dW_t^{(n)}.
\label{micro wealth}
\end{equation}
Thus, the volatility estimator in fact enters filtering and wealth separately through two channels: the standardized observation and the physical allocation.
In the continuous-observation model, the quadratic-variation identity holds directly, namely $U_j=(\sigma^*)^2$.
Therefore, in this setting the volatility estimator is equivalent to the oracle volatility, and $\chi_n$ is identically equal to $1$, so that the filtering, wealth, actor--critic, Sharpe-ratio, and multiplier dynamics used below all reduce exactly to their forms under oracle volatility.

In addition, in an infill regime without market microstructure noise, discrete realized quadratic variation converges to continuous quadratic variation at the rate $O_p(\Delta^{1/2})$ as the grid width $\Delta\to0$ \citep[Theorem~1]{barndorff2002realized}. Let $\widehat\sigma_{n,\Delta}$ denote the estimator obtained from discrete realized quadratic variation, and let
$\chi_{n,\Delta}:=\sigstar/\widehat\sigma_{n,\Delta}$.
By the consistency above, finite averaging, and continuity of the projection, for each fixed episode $n\ge1$,
\[
\widehat\sigma_{n,\Delta}\xrightarrow{p}\sigstar,
\quad
\chi_{n,\Delta}\xrightarrow{p}1.
\]
Let $\widetilde m^{(n,\Delta)},\widetilde Y^{(n,\Delta)}$ denote the system using $\chi_{n,\Delta}$, coupled with the same Brownian motion, sampled action, and initial values.
By the explicit solution of the scalar linear filtering equation,
\[
\widetilde m_t^{(n,\Delta)}-m_t^{(n)}
=R_t^{(n)}(\chi_{n,\Delta}-1)
(\rhostar t+W_t^{(n)}),
\]
so this filter perturbation converges to zero uniformly in probability on $[0,T]$.
For $Y_0^{(n)}\ne0$, the explicit wealth solution gives
\[
\begin{aligned}
\log\left|\frac{\widetilde Y_t^{(n,\Delta)}}{Y_t^{(n)}}\right|
=-(\chi_{n,\Delta}-1)\rhostar\int_0^t\rho_s^{(n)}ds
-\frac{\chi_{n,\Delta}^2-1}{2}\int_0^t(\rho_s^{(n)})^2ds-(\chi_{n,\Delta}-1)\int_0^t\rho_s^{(n)}dW_s^{(n)}.
\end{aligned}
\]
For fixed $n$, the coefficients $a_s^{(n)},b_s^{(n)}$ and $V_s^{(n)}$ are bounded on $[0,T]$. Moreover, $|\phi_n|\le c_n<\infty$, and \eqref{eq:m-explicit} gives $\sup_{s\le T}\E|m_s^{(n)}|^2<\infty$. Thus the conditional Gaussian law \eqref{eq:actor-draw-centered} and the tower property yield
\[
\begin{aligned}
\E|\rho_s^{(n)}|^2
=\E\left[|M_s^{\phi_n}|^2+V_s^{(n)}\right]\le 2|a_s^{(n)}|^2\E|m_s^{(n)}|^2
+2|b_s^{(n)}|^2c_n^2+V_s^{(n)}
\le C_n, \quad 0\le s\le T,
\end{aligned}
\]
where $C_n$ is a deterministic constant for this fixed $n$. By joint measurability and Tonelli's theorem,
\[
\E\left[\int_0^T|\rho_s^{(n)}|^2ds\right]
=\int_0^T\E|\rho_s^{(n)}|^2ds
\le TC_n<\infty.
\]
Since a nonnegative random variable with finite expectation is finite almost surely, $\int_0^T|\rho_s^{(n)}|^2ds<\infty$ almost surely. Hence these integrals and $\sup_{t\le T}|Y_t^{(n)}|$ are finite almost surely. Since $\chi_{n,\Delta}\to1$ in probability,
\[
\sup_{0\le t\le T}|\widetilde Y_t^{(n,\Delta)}-Y_t^{(n)}|\xrightarrow{p}0.
\]
For $Y_0^{(n)}=0$, both wealth deviations vanish identically.

For the approximation of the continuous-time representation by the discrete sampled dynamics and the corresponding policy-gradient estimator, see \citet{jia2026accuracy}. Therefore, the subsequent continuous-time convergence proofs do not track the discrete execution error again.

To write the policy coefficients for general $c>0$, let $\tau_t^{(n)} := \log({R_t^{(n)}}/{R_T^{(n)}})$.
Define $p_c,h_c$ as the solutions on $\tau\ge0$ of the following initial-value problems:
\begin{align}
\frac{dp_c}{d\tau}
&=
\frac{(c-2)p_c^2-(5c+2)p_c+2c}{c+2},
\quad p_c(0)=0,
\label{eq:pc-ode}\\
\frac{dh_c}{d\tau}
&=
-\frac{2c-(c-2)p_c}{c+2}h_c,
\quad h_c(0)=1.
\label{eq:hc-ode}
\end{align}
They come from the HJB coefficient equations for general $c$. Indeed, if
$p_c=R_tP(t)$ and $h_c=1+R_t\overline Q(t)$, then using
$\dot R_t=-R_t^2$ and
$d\tau_t/dt=-R_t$, the equations for $P,\overline Q$
\begin{align*}
0
&=
\dot P-2R_tP+R_t^2P^2+c
-\frac{(c+2R_tP)^2}{c+2}, \quad P(T)=0,\\
0
&=
\dot{\overline Q}+R_t(P-\overline Q)
+R_t^2P\overline Q
-\frac{2(1+R_t\overline Q)(c+2R_tP)}{c+2}, \quad \overline{Q}(T)=0,
\end{align*}
reduce exactly to \eqref{eq:pc-ode}--\eqref{eq:hc-ode}. Therefore, the actor sampling distribution in the $n$th episode is
\begin{equation*}
q_t^\phi
=
N\left(
M_t^\phi,V_t^{(n)}
\right),
\quad
M_t^\phi
=
a_t^{(n)}m_t^{(n)}+b_t^{(n)}\phi,
\quad
V_t^{(n)}=\frac{c}{c+2}R_t^{(n)},
\end{equation*}
where
\begin{equation}
a_t^{(n)}
=
\frac{c+2p_c(\tau_t^{(n)})}{c+2},
\quad
b_t^{(n)}
=
\frac{2h_c(\tau_t^{(n)})}{c+2}.
\label{eq:ab-general-c}
\end{equation}

\subsection*{Proof of Lemma~\ref{lem:main-filter}}
Because the update of the filtering state depends only on the market observation process, after giving the basic setup above, we first analyze the evolution of the filtering state.

\begin{proof}
In addition to notations in Lemma~\ref{lem:main-filter}, set $g_n:=\frac1{p_n+T}$ and $
a_n:=\frac{p_n}{p_n+T}$. Then
\begin{equation}
R_t^{(n)}
=\frac1{p_n+t}
=\frac{R_0}
{(n+1)^\alpha+R_0t},
\quad 0\le t\le T,
\label{eq:decay-R-path}
\end{equation}
Combining this with \eqref{eq:filter-m-n} gives
\begin{equation*}
de_t^{(n)}
=
-\frac{e_t^{(n)}}{p_n+t}\,dt
+
\frac1{p_n+t}\,dW_t^{(n)}.
\end{equation*}
As $p_n+t$ is a deterministic $C^1$ finite-variation process and $e_t^{(n)}$ is a continuous It\^o process, both are continuous semimartingales. Applying the It\^o product rule to $(p_n+t)e_t^{(n)}$ gives
\begin{align*}
d\left[(p_n+t)e_t^{(n)}\right]
=
(p_n+t)\,de_t^{(n)}
+e_t^{(n)}\,d(p_n+t)
+d\langle p_n+t,e^{(n)}\rangle_t
=dW_t^{(n)}.
\end{align*}
Integrating the above identity and using the initial error $e_0^{(n)}$ gives
\begin{equation}
(p_n+t)e_t^{(n)}
=p_n e_0^{(n)}
+W_t^{(n)},
\quad\text{i.e. }\;
m_t^{(n)}-\rhostar
=
\frac{p_n e_0^{(n)}
+W_t^{(n)}}
{p_n+t}.
\label{eq:m-explicit}
\end{equation}
Setting $t=T$ and using $m_0^{(n+1)}=m_T^{(n)}$ gives \eqref{eq:main-filter-recursion}. It follows from iteration that
\begin{equation*}
e_0^{(n)}
=
A_{0,n}e_0^{(0)}
+\sum_{j=0}^{n-1}A_{j+1,n}g_jW_T^{(j)},
\quad
A_{k,n}:=\prod_{\ell=k}^{n-1}a_\ell,
\quad A_{n,n}:=1.
\end{equation*}
Consequently, $e_0^{(n)}$ is Gaussian. Write $\mu_n:=\E[e_0^{(n)}]$ and
$v_n:=\Var(e_0^{(n)})$.
Then
\begin{align}
\mu_{n+1}=a_n\mu_n,
\quad
v_{n+1}=a_n^2v_n+Tg_n^2.
\label{eq:e-mean-variance-recursion}
\end{align}
As $0<g_n\le1/T$,
\[
v_{n+1}-v_n
=-Tg_n(2-Tg_n)v_n+Tg_n^2.
\]
If $v_n\ge T^{-1}$, then
\[
v_{n+1}-v_n
\le
-Tg_n(2-Tg_n)T^{-1}+Tg_n^2
=-2g_n+2Tg_n^2
\le0.
\]
If $v_n<T^{-1}$, then by \eqref{eq:e-mean-variance-recursion} and
$a_n=1-Tg_n$,
\[
v_{n+1}
=(1-Tg_n)^2v_n+Tg_n^2
\le
(1-Tg_n)^2T^{-1}+Tg_n^2
=T^{-1}-2g_n+2Tg_n^2
\le T^{-1}.
\]
Therefore, induction gives $\sup_n v_n\le\max\{v_0,T^{-1}\}$.
Hence the Gaussian moment formula gives the uniform moment bound in Lemma~\ref{lem:main-filter}.

We next prove the rates for the mean and the second moment. Since
$g_n\asymp(n+1)^{-\alpha}$ and $0<\alpha<1$,
\[
\sum_{j=0}^{n-1}g_j\ge c_\alpha n^{1-\alpha}.
\]
Therefore
\[
|\mu_n|
=|e_0^{(0)}|\prod_{j=0}^{n-1}(1-Tg_j)
\le
|e_0^{(0)}|
\exp\left\{-T\sum_{j=0}^{n-1}g_j\right\}
\le
Ce^{-c_\alpha n^{1-\alpha}}.
\]
Also, since
\[
g_n-g_{n+1}
=O((n+1)^{-\alpha-1})
=o(g_n^2),
\]
we can take $N$ sufficiently large so that
$
g_n-g_{n+1}\le Tg_n^2$ and $T g_n\le 1/2$ hold for
$ n\ge N$.
Since $v_N<\infty$ and $g_N>0$, choose a finite constant
$C\ge\max\{2,v_N/g_N\}$. Then $v_N\le Cg_N$.
If $v_n\le Cg_n$ for some $n\ge N$, then
\eqref{eq:e-mean-variance-recursion} gives
\begin{align*}
v_{n+1}
=(1-Tg_n)^2v_n+Tg_n^2
\le
Cg_n-(2TC-T)g_n^2+CT^2g_n^3.
\end{align*}
Because $0<g_n\le1/T$ and $C\ge2$,
\[
v_{n+1}
\le
Cg_n-CTg_n^2
\le
Cg_n-C(g_n-g_{n+1})
=Cg_{n+1}.
\]
Induction yields $v_n\le Cg_n$ for all $n\ge N$. Hence
\[
v_n=O(g_n)=O((n+1)^{-\alpha}).
\]
Combining this with the exponential bound for the mean gives
\[
\E[(e_0^{(n)})^2]
=v_n+\mu_n^2
\le C(n+1)^{-\alpha},
\]
This proves the bounds on the mean and the second moment when
$0<\alpha<1$.
Now suppose $\alpha=1$, so that $R_0T>1/2$. Then
\[
g_n=\frac{R_0}{n+1+R_0T},
\quad
a_n=\frac{n+1}{n+1+R_0T}.
\]
The same product estimate for $\mu_n$, together with
$T\sum_{j=0}^{n-1}g_j\ge R_0T\log(n+1)-C$, yields
\[
|\mu_n|\le C(n+1)^{-R_0T}\le C(n+1)^{-1/2}.
\]
Also, choose $C\ge v_0$ sufficiently large. If
$v_n\le C/(n+1)$, then
\[
v_{n+1}
\le
\frac{C(n+1)+TR_0^2}{(n+1+R_0T)^2}
\le \frac{C}{n+2},
\]
where the last inequality holds after choosing $C$ sufficiently large,
because $2R_0T-1>0$. Induction therefore gives
$v_n=O((n+1)^{-1})$. In both cases,
$\E[(e_0^{(n)})^2]=v_n+\mu_n^2\le C(n+1)^{-\alpha}$,
which proves the two rate estimates in Lemma~\ref{lem:main-filter}.

Finally, we prove a.s. convergence. For any $\varepsilon>0$, when $n$ is sufficiently large,
$|\mu_n|\le\varepsilon/2$. By the Gaussian tail bound,
\[
\mathbb P(|e_0^{(n)}|>\varepsilon)
\le
2\exp\{-c_\varepsilon(n+1)^\alpha\}.
\]
The right-hand side is summable, so the Borel--Cantelli lemma gives
$e_0^{(n)}\to0$ a.s.
\end{proof}

Furthermore, \eqref{eq:m-explicit} directly implies the conditional moment bound needed below. For $n\ge0$,
\begin{align*}
\E\left[
|m_t^{(n)}|^2
\mid\F_n
\right]
&\le
C_{R_0,T}\left(
1+|e_0^{(n)}|^2
\right),
\quad
0\le t\le T,
\end{align*}
where $C_{R_0,T}<\infty$ is independent of $n,e_0^{(n)},\phi_n$.

\subsection*{Conditional Second-Moment Bound for Actor Noise}
Having given the convergence analysis for filtering, we next give the convergence analysis for the parameters in the actor--critic iteration.
First, we derive the actor-noise bound. Let $\G_t^{(n)}$ denote the information before the current action is sampled. Then
$\phi_n$, $m_t^{(n)}$, $R_t^{(n)}$, $M_t^{\phi_n}$, and $V_t^{(n)}$ are all $\G_t^{(n)}$-measurable, and write
\begin{equation}
\rho_t^{(n)}-M_t^{\phi_n}
\mid\G_t^{(n)}
\sim
N\left(0,V_t^{(n)}\right).
\label{eq:actor-draw-centered}
\end{equation}
The Gaussian score is
\begin{equation}
s_t^{(n)}
:=
\partial_\phi
\log q_t^{\phi_n}
\left(
\rho_t^{(n)}
\right)
=
b_t^{(n)}
\frac{
\rho_t^{(n)}-M_t^{\phi_n}
}{
V_t^{(n)}
}.
\label{eq:score}
\end{equation}
By the conditional Gaussian moment formula,
\begin{equation*}
\E\left[|s_t^{(n)}|^2\mid\G_t^{(n)}\right]
=
\frac{|b_t^{(n)}|^2}{V_t^{(n)}}.
\end{equation*}
For all $n\ge0$,
\[
\tau_t^{(n)}
=\log\frac{p_n+T}{p_n+t}
\in[0,\overline\tau_0],
\quad
\overline\tau_0:=\log(1+R_0T).
\]
Let
\[
f_c(p):=\frac{(c-2)p^2-(5c+2)p+2c}{c+2}.
\]
When $c\ne2$, it holds that
\[
f_c(0)>0,\quad
f_c(p)>0\ (0\le p<r_-),\quad f_c(r_-)=0.
\]
Thanks to the uniqueness of the solution and phase-line analysis for the one-dimensional autonomous ODE
\citep[Lemma~1.1]{teschl2012ode}, we get
\[
0\le p_c(\tau)<r_-,\quad p_c(\tau)\uparrow r_-.
\]
When $c=2$,
\[
p_2(\tau)=\frac13(1-e^{-3\tau})\uparrow\frac13.
\]
Moreover,
\[
2c-(c-2)p_c(\tau)>0,
\]
because when $0<c\le2$ its left-hand side is no smaller than $2c$. When $c>2$, it follows that
\[
f_c\left(\frac{2c}{c-2}\right)
=-\frac{4c}{c-2}<0
\;\Longrightarrow\;
p_c(\tau)<r_-<\frac{2c}{c-2}.
\]
Consequently, 
\[
h_c(\tau)
=\exp\left\{-\int_0^\tau
\frac{2c-(c-2)p_c(s)}{c+2}\,ds\right\}
\in(0,1].
\]
Then $0\le p_c(\tau)\le p_c(\overline\tau_0),
\quad
h_c(\overline\tau_0)\le h_c(\tau)\le1$. By \eqref{eq:ab-general-c}, we obtain
\begin{equation}
\frac{c}{c+2}
\le
a_t^{(n)}
\le
\frac{c+2 p_c(\overline\tau_0)}{c+2},
\quad
\frac{2h_c(\overline\tau_0)}{c+2}
\le
b_t^{(n)}
\le
\frac2{c+2}.
\label{eq:ab-uniform}
\end{equation}
Moreover,
\begin{equation}
\frac{cR_0}
{(c+2)\bigl((n+1)^\alpha+R_0T\bigr)}
\le
V_t^{(n)}
\le
\frac{cR_0}{(c+2)(n+1)^\alpha}.
\label{eq:V-bounds}
\end{equation}
Therefore the conditional second moment of the score is controlled as follows:
\begin{equation*}
\E\left[
|s_t^{(n)}|^2
\mid\G_t^{(n)}
\right]
\le
C_{c,R_0}(n+1)^\alpha,
\end{equation*}
where $C_{c,R_0}<\infty$ is a constant. Moreover, uniformly for $t\in[0,T]$,
\[
\tau_t^{(n)}
=\int_t^T\frac{R_0^{(n)}}{1+sR_0^{(n)}}\,ds
=(T-t)R_0^{(n)}+O\left((R_0^{(n)})^2\right),
\quad
0\le\tau_t^{(n)}\le TR_0^{(n)}.
\]
Since $R_0^{(n)}=R_0(n+1)^{-\alpha}$ and the ODE solution is smooth near $\tau=0$,
\begin{equation}
p_c(\tau_t^{(n)})=O((n+1)^{-\alpha}),
\quad
1-h_c(\tau_t^{(n)})=O((n+1)^{-\alpha}).
\label{eq:pc-hc-small}
\end{equation}
For any fixed $c>0$, define
\begin{equation*}
J_n^\theta(t,x,m)
=
A_n^\theta(t,m)(x-w)^2,
\quad
A_n^\theta(t,m)>0,
\end{equation*}
where
\begin{align}
\log A_n^\theta(t,m)
&=
\frac12P_n(t)m^2+\overline Q_n(t)\theta m
+\frac12K_{1,n}(t)\theta^2+K_{2,n}(t),
\label{eq:parametric-critic-log}\\
P_n(t)
&=
\frac{p_c(\tau_t^{(n)})}{R_t^{(n)}},
\quad
\overline Q_n(t)
=
\frac{h_c(\tau_t^{(n)})-1}{R_t^{(n)}},
\label{eq:parametric-critic-PQ}\\
\dot K_{1,n}(t)
&=
1-\frac{c-2}{c+2}h_c(\tau_t^{(n)})^2,
\quad K_{1,n}(T)=0,
\label{eq:parametric-critic-K1}\\
\dot K_{2,n}(t)
&=
-\frac12R_t^{(n)}p_c(\tau_t^{(n)})
-\frac{cR_t^{(n)}}2\log\frac{c+2}{c},
\quad K_{2,n}(T)=0.
\label{eq:parametric-critic-K2}
\end{align}
\eqref{eq:pc-hc-small} shows that
$p_c(\tau_t^{(n)})/R_t^{(n)}$ and
$(h_c(\tau_t^{(n)})-1)/R_t^{(n)}$ are uniformly bounded. By
\eqref{eq:parametric-critic-PQ}--\eqref{eq:parametric-critic-K2}
and the smoothness of the ODE, all coefficients and their time derivatives are bounded on $[0,T]$.
Therefore there exists $C_{c,R_0,T}<\infty$ such that
\begin{equation}
\sup_{n\ge0}\sup_{0\le t\le T}
\left(
\begin{aligned}
&|P_n(t)|+|\overline Q_n(t)|+|K_{1,n}(t)|+|K_{2,n}(t)|\\
&\quad
+|\dot P_n(t)|+|\dot{\overline Q}_n(t)|
+|\dot K_{1,n}(t)|+|\dot K_{2,n}(t)|
\end{aligned}
\right)
\le C_{c,R_0,T}.
\label{eq:critic-coefficient-bound}
\end{equation}

\subsection*{Conditional Mean Actor Direction under a Fixed Critic Projection}

To remain consistent with the discrete update, after the sampled action $\rho_t^{(n)}$ is given, the wealth process used in the proof of actor--critic convergence is based on the microscopic dynamics \eqref{micro wealth}.
Therefore, the explicit cross-variation is $d\langle Y^{(n)},m^{(n)}\rangle_t
=
-\rho_t^{(n)}R_t^{(n)}Y_t^{(n)}dt$.
Let $\ell_n^\theta:=\log A_n^\theta$. Applying It\^o formula to $J_n^\theta=A_n^\theta (Y^{(n)})^2$, we may write
\begin{equation}
dJ_n^\theta
=
J_n^\theta\left(\mathcal D_t^{\rho_t^{(n)},\theta,n}dt
+\Gamma_t^{\rho_t^{(n)},\theta,n}dW_t^{(n)}\right).
\label{eq:micro-Ito-J}
\end{equation}
where
\begin{align*}
\mathcal D_t^{\rho,\theta,n}
=&
\partial_t\ell_n^\theta
+R_t^{(n)}(\rhostar-m_t^{(n)})\partial_m\ell_n^\theta
+\frac12(R_t^{(n)})^2
\left(
\partial_{mm}\ell_n^\theta
+(\partial_m\ell_n^\theta)^2
\right)
\\
&+\rho^2
-2\rho
\left(
\rhostar+R_t^{(n)}\partial_m\ell_n^\theta
\right),\\
\Gamma_t^{\rho,\theta,n}
=&
-2\rho+R_t^{(n)}\partial_m\ell_n^\theta.
\end{align*}
After adding the KL drift, define the sampled-action critic residual
\begin{equation*}
d\mathcal R_t^{\theta,n}(\phi)
:=
dJ_n^\theta
+\lambda_tJ_n^\theta
D_{\mathrm{KL}}
\left(q_t^\phi\Vert N(m_t^{(n)},R_t^{(n)})\right)dt.
\end{equation*}
The corresponding normalized residual drift is defined by
\begin{equation}
G_t^{\rho,\theta,n}(\phi)
=
\mathcal D_t^{\rho,\theta,n}
+
\lambda_tD_{\mathrm{KL}}
\left(q_t^\phi\Vert N(m_t^{(n)},R_t^{(n)})\right).
\label{eq:sampled-residual-drift}
\end{equation}
By \eqref{eq:micro-Ito-J} and \eqref{eq:sampled-residual-drift},
\begin{equation}
\frac{d\mathcal R_t^{\theta,n}(\phi)}{J_n^\theta}
=
G_t^{\rho_t^{(n)},\theta,n}(\phi)dt
+
\Gamma_t^{\rho_t^{(n)},\theta,n}dW_t^{(n)}.
\label{eq:parametric-residual-decomposition}
\end{equation}
For the actor, using
$dm_t^{(n)}=R_t^{(n)}(\rhostar-m_t^{(n)})dt+R_t^{(n)}dW_t^{(n)}$
in \eqref{eq:main-centered-normalized-residual}, with
$\beta_{\mathrm{cv}}=1$, gives the episode-wise centered decomposition
\begin{equation}
d\mathcal C_t^{\theta_n,\phi_n}
=
\overline G_t^{\theta_n,\phi_n,n}dt
+\overline\Gamma_t^{\phi_n,n}dW_t^{(n)},
\label{eq:centered-residual-decomposition}
\end{equation}
where
\begin{align}
\overline G_t^{\theta,\phi,n}
&:=
(\rho_t-M_t^\phi)\left(
\rho_t+M_t^\phi-2\rhostar
-2R_t^{(n)}\partial_m\ell_n^\theta(t,m_t^{(n)})
\right)-V_t^{(n)},
\label{eq:centered-residual-drift}\\
\overline\Gamma_t^{\phi,n}
&:=-2(\rho_t-M_t^\phi).
\label{eq:centered-residual-diffusion}
\end{align}
Thus the continuous-time actor increment corresponding to the adjusted update is
\begin{align}
\widehat F^{\mathrm c}_{n+1}(\phi_n,\theta_n)
=
\int_0^T
s_t^{(n)}
\overline G_t^{\theta_n,\phi_n,n}dt
+
\int_0^T
s_t^{(n)}
\overline\Gamma_t^{\phi_n,n}
dW_t^{(n)}
+
\int_0^T
\lambda_t\partial_\phi
D_{\mathrm{KL}}
\left(q_t^{\phi_n}\Vert N(m_t^{(n)},R_t^{(n)})\right)dt.
\label{eq:centered-actor-increment}
\end{align}
The second integral is the centered corresponding martingale-noise term.
For the third integral, the Gaussian KL formula gives
\begin{equation}
\lambda_t\partial_\phi
D_{\mathrm{KL}}
\left(q_t^\phi\Vert N(m_t^{(n)},R_t^{(n)})\right)
=
c\,b_t^{(n)}(M_t^\phi-m_t^{(n)}).
\label{eq:correct-KL-gradient}
\end{equation}
Its stochastic-approximation decomposition is
\begin{equation}
\widehat F_{n+1}^c(\phi_n,\theta_n)
=
F_n(\phi_n,\theta_n)+\xi_{n+1},
\quad
F_n(\phi,\theta)
:=
\E\left[
\widehat F_{n+1}^c(\phi,\theta)
\mid\F_n
\right],
\label{eq:two-part-decomposition}
\end{equation}
where
$\xi_{n+1}
:=
\widehat F_{n+1}^c(\phi_n,\theta_n)-F_n(\phi_n,\theta_n)$ and
$\E[\xi_{n+1}\mid\F_n]=0$.
Note that $F_n$ depends on the episode-initial state $(m_0^{(n)},R_0^{(n)})$. The sampled-action estimating statistic corresponding to critic update \eqref{eq:main-normalized-critic-direction} is defined by
\begin{equation*}
\widehat H_{n+1}(\theta_n,\phi_n)
:=
\int_0^T
\frac{\partial_\theta J_n^{\theta_n}}
{(J_n^{\theta_n})^2}
d\mathcal R_t^{\theta_n,n}
(\phi_n;\rho_t^{(n)}).
\end{equation*}
If the terminal wealth of the $n$th episode is denoted by $X_T^{(n)}$, then a multiplier update is triggered at each time $n\equiv0\pmod {N_w}$, and the single sample $X_T$ is replaced by the sample mean of the most recent $N_w$ terminal wealth values. That is, the $(r+1)$st multiplier update occurs at episode $(r+1)N_w$ and
\[
{\mathcal I}_{r+1}:=\{rN_w+1,\ldots,(r+1)N_w\},
\]
where $N_w \in \mathbb{Z}^+$. The batch terminal-wealth statistic used for the update at this time is
\begin{equation}
\overline X_{T,r+1}
:=
\frac1{|\mathcal I_{r+1}|}
\sum_{\ell\in\mathcal I_{r+1}}X_T^{(\ell)}
=
\frac1{N_w}\sum_{\ell=rN_w+1}^{(r+1)N_w}X_T^{(\ell)}.
\label{eq:terminal-wealth-batch-statistic}
\end{equation}

\begin{lemma}
\label{lem:actor-field}
For all $n\ge0$, $\phi\in\mathbb R$, and $\theta\in\Theta$,
\begin{equation}
F_n(\phi,\theta)
:=
\E[\widehat F_{n+1}^c(\phi,\theta)\mid\F_n]
=
K_n(\phi-\rhostar)+B_n(\theta-\rhostar)
\label{eq:critic-biased-field}
\end{equation}
where
\begin{equation}
K_n:=\frac4{c+2}\int_0^T h_c(\tau_t^{(n)})^2dt,
\quad
B_n:=\frac4{c+2}\int_0^T
h_c(\tau_t^{(n)})(1-h_c(\tau_t^{(n)}))dt.
\label{eq:K-B-def}
\end{equation}
and
\begin{equation}
K_n
\ge
\kappa_{c,R_0}
:=
\frac{4Th_c^2(\overline\tau_0)}{c+2}
>0,
\quad
0\le B_n\le C_{c,R_0,T}(n+1)^{-\alpha}.
\label{eq:explicit-curvature-lower-bound}
\end{equation}
\end{lemma}

\begin{proof}
Fix $t$ and condition on $\G_t^{(n)}$. Write
$M:=M_t^\phi$, $V:=V_t^{(n)}$,
$b:=b_t^{(n)}$, $Z:=\rho_t^{(n)}-M$. By \eqref{eq:actor-draw-centered} and \eqref{eq:score}, these abbreviations
satisfy $Z\mid\G_t^{(n)}\sim N(0,V)$ and $s_t^{(n)}=bZ/V$. By
\eqref{eq:centered-residual-drift}--\eqref{eq:centered-residual-diffusion},
\begin{align}
s_t^{(n)}\overline G_t^{\theta,\phi,n}
&=
b\left[
2\frac{Z^2}{V}
\left(
M-\rhostar
-R_t^{(n)}\partial_m\ell_n^\theta(t,m_t^{(n)})
\right)
+\frac{Z^3}{V}-Z
\right],
\label{eq:centered-score-drift-product}\\
s_t^{(n)}\overline\Gamma_t^{\phi,n}
&=-2b\frac{Z^2}{V}.
\label{eq:centered-score-diffusion-product}
\end{align}
The Gaussian identities
$\E[Z]=\E[Z^3]=0$, $\E[Z^2]=V$, $\E[Z^4]=3V^2$, $\E[Z^6]=15V^3$ imply
\begin{equation}
\E\left[
s_t^{(n)}\overline G_t^{\theta,\phi,n}
\middle|\G_t^{(n)}
\right]
=2b_t^{(n)}
\left[
M_t^\phi-\rhostar
-R_t^{(n)}\partial_m\ell_n^\theta(t,m_t^{(n)})
\right].
\label{eq:centered-score-mean}
\end{equation}
For the conditional second-moment estimates, use the policy formula
preceding \eqref{eq:ab-general-c} and differentiate
\eqref{eq:parametric-critic-log} with respect to $m$.
The uniform bounds on $a_t^{(n)}$, $b_t^{(n)}$,
$R_t^{(n)}P_n(t)$, and $R_t^{(n)}\overline Q_n(t)$, together with the
fixed projection $\theta\in\Theta$, give
\[
\left|
M-\rhostar
-R_t^{(n)}\partial_m\ell_n^\theta(t,m_t^{(n)})
\right|
\le C_{c,R_0,T,\Theta}
\left(1+|m_t^{(n)}|+|\phi|\right).
\]
Also, \eqref{eq:V-bounds} gives $V\le C_{c,R_0}$. Applying
$|x+y+z|^2\le3(|x|^2+|y|^2+|z|^2)$ to
\eqref{eq:centered-score-drift-product} and using the Gaussian fourth
and sixth moments yields
\begin{align}
\E\left[
\left|s_t^{(n)}\overline G_t^{\theta,\phi,n}\right|^2
\middle|\G_t^{(n)}
\right]
&\le
3b^2\left[
4\frac{\E[Z^4\mid\G_t^{(n)}]}{V^2}
\left|
M-\rhostar
-R_t^{(n)}\partial_m\ell_n^\theta(t,m_t^{(n)})
\right|^2
+\frac{\E[Z^6\mid\G_t^{(n)}]}{V^2}
+\E[Z^2\mid\G_t^{(n)}]
\right]\notag\\
&=3b^2\left[
12\left|
M-\rhostar
-R_t^{(n)}\partial_m\ell_n^\theta(t,m_t^{(n)})
\right|^2
+16V
\right]\notag\\
&\le C_{c,R_0,T,\Theta}
\left(1+|m_t^{(n)}|^2+|\phi|^2\right),
\label{eq:centered-score-drift-bound}
\end{align}
Similarly, \eqref{eq:centered-score-diffusion-product} and the Gaussian
fourth-moment identity above give
\begin{equation}
\E\left[
\left|s_t^{(n)}\overline\Gamma_t^{\phi,n}\right|^2
\middle|\G_t^{(n)}
\right]
=4b^2\frac{\E[Z^4\mid\G_t^{(n)}]}{V^2}
=12b^2
\le C_{c,R_0,T}.
\label{eq:centered-score-diffusion-bound}
\end{equation}
To verify the conditional integrability of the drift integral, use the
conditional second-moment bound following Lemma~\ref{lem:main-filter}.
That bound, the conditional Tonelli theorem and the tower property, then give
\[
\E\left[
\int_0^T
\left|s_t^{(n)}\overline G_t^{\theta,\phi,n}\right|^2dt
\middle|\F_n
\right]
=\int_0^T
\E\left[
\E\left[
\left|s_t^{(n)}\overline G_t^{\theta,\phi,n}\right|^2
\middle|\G_t^{(n)}
\right]
\middle|\F_n
\right]dt
\le C_{c,R_0,T,\Theta}
\left(1+|e_0^{(n)}|^2+|\phi|^2\right)<\infty.
\]
Thus the integral Cauchy--Schwarz inequality makes the drift integral in
\eqref{eq:centered-actor-increment} conditionally integrable.

For the stochastic integral, the conditional Tonelli theorem, the tower
property  and \eqref{eq:centered-score-diffusion-bound} give
\[
\E\left[
\int_0^T
\left|s_t^{(n)}\overline\Gamma_t^{\phi,n}\right|^2dt
\middle|\F_n
\right]=
\int_0^T
\E\left[
\E\left[
\left|s_t^{(n)}\overline\Gamma_t^{\phi,n}\right|^2
\middle|\G_t^{(n)}
\right]
\middle|\F_n
\right]dt
\le C_{c,R_0,T}T<\infty.
\]
Thus the stochastic integral is square integrable and has conditional
mean zero given $\F_n$.

For the KL-gradient integral, \eqref{eq:correct-KL-gradient}, the policy
formula preceding \eqref{eq:ab-general-c}, and the uniform bound on
$b_t^{(n)}$ imply that
\[
\left|
\lambda_t\partial_\phi D_{\mathrm{KL}}
\left(q_t^\phi\Vert N(m_t^{(n)},R_t^{(n)})\right)
\right|^2
\le
C_{c,R_0,T}
\left(|m_t^{(n)}|^2+|\phi|^2\right).
\]
Therefore, the conditional second-moment bound following
Lemma~\ref{lem:main-filter}, the conditional Tonelli theorem and the tower
property imply
\begin{align*}
&\E\left[
\int_0^T
\left|
\lambda_t\partial_\phi D_{\mathrm{KL}}
\left(q_t^\phi\Vert N(m_t^{(n)},R_t^{(n)})\right)
\right|^2dt
\middle|\F_n
\right]\le
C_{c,R_0,T}
\left(1+|e_0^{(n)}|^2+|\phi|^2\right)<\infty.
\end{align*}
The integral Cauchy--Schwarz inequality then gives
\[
\E\left[
\int_0^T
\left|
\lambda_t\partial_\phi D_{\mathrm{KL}}
\left(q_t^\phi\Vert N(m_t^{(n)},R_t^{(n)})\right)
\right|dt
\middle|\F_n
\right]<\infty.
\]

Adding the explicit KL derivative
$c\,b_t^{(n)}(M_t^\phi-m_t^{(n)})$ from
\eqref{eq:correct-KL-gradient} to \eqref{eq:centered-score-mean} gives the
conditional mean actor direction per unit time
\begin{equation}
b_t^{(n)}
\left[
(c+2)M_t^\phi-cm_t^{(n)}-2\rhostar
-2R_t^{(n)}\partial_m\log A_n^\theta(t,m_t^{(n)})
\right].
\label{eq:conditional-field-before-cancellation}
\end{equation}
By \eqref{eq:parametric-critic-log},
$\partial_m\log A_n^\theta=P_n(t)m+\overline Q_n(t)\theta$. Use the coefficient identities in \eqref{eq:parametric-critic-PQ} and \eqref{eq:ab-general-c}. Consequently, \eqref{eq:conditional-field-before-cancellation} equals
\[
\frac4{c+2}
\left[
h_c(\tau_t^{(n)})^2(\phi-\rhostar)
+h_c(\tau_t^{(n)})(1-h_c(\tau_t^{(n)}))
(\theta-\rhostar)
\right].
\]
Integration over $[0,T]$ gives
\eqref{eq:critic-biased-field}--\eqref{eq:K-B-def}. Finally,
$h_c(\overline\tau_0)\le h_c(\tau_t^{(n)})\le1$ and
\eqref{eq:pc-hc-small} yield
\eqref{eq:explicit-curvature-lower-bound}.
\end{proof}

\begin{remark}
The derivation of the conditional mean actor direction above first treats the multiplier $w$ as a given constant for the following reason:
the actor field uses the product of the normalized residual coefficient and the score
$s_t^{(n)}$. After normalization, the common quadratic factor $(X_t-w)^2$ has already canceled before taking the conditional expectation. Therefore
\eqref{eq:conditional-field-before-cancellation}
gives a conditional mean actor direction that does not depend on $w$. Similarly, the moment bound for the actor noise does not need to absorb $w$ into the constants.
Thus, Theorem~\ref{thm:main-actor} may be understood as holding for any episode-wise predictable multiplier sequence.
\end{remark}

\begin{lemma}
\label{lem:correct-actor-noise}
Under the fixed projection \eqref{eq:main-projected-actor-recursions}, let $\phi_n, \theta_n$ be $\F_n$-measurable. Then there exists $C<\infty$ such that, for all $n\ge0$,
\begin{align}
\E[|\widehat F_{n+1}^c(\phi_n,\theta_n)|^2\mid\F_n]
&\le
C(1+|e_0^{(n)}|^4+|\phi_n|^4),
\label{eq:correct-Fhat-bound}\\
\E[|\xi_{n+1}|^2\mid\F_n]
&\le
C(1+|e_0^{(n)}|^4+|\phi_n|^4).
\label{eq:correct-xi-bound}
\end{align}
\end{lemma}

\begin{proof}
Denote the three terms in \eqref{eq:centered-actor-increment} by
\[
I_{1,n}:=\int_0^T
s_t^{(n)}\overline G_t^{\theta_n,\phi_n,n}dt,
\quad
I_{2,n}:=\int_0^T
s_t^{(n)}\overline\Gamma_t^{\phi_n,n}dW_t^{(n)},
\]
and
\[
I_{3,n}:=\int_0^T K_t^{(n)}dt,
\quad
K_t^{(n)}:=
\lambda_t\partial_\phi D_{\mathrm{KL}}
\left(q_t^{\phi_n}\Vert N(m_t^{(n)},R_t^{(n)})\right).
\]
Thus
$\widehat F^{\mathrm c}_{n+1}=I_{1,n}+I_{2,n}+I_{3,n}$.
For $I_{1,n}$, the integral Cauchy--Schwarz inequality first gives the
pathwise bound
\[
|I_{1,n}|^2
\le
T\int_0^T
\left|s_t^{(n)}\overline G_t^{\theta_n,\phi_n,n}\right|^2dt.
\]
Taking conditional expectation, applying the conditional Tonelli theorem
and the tower property, and then using
\eqref{eq:centered-score-drift-bound} and the conditional second-moment
bound following Lemma~\ref{lem:main-filter}, we obtain
\begin{align*}
\E[|I_{1,n}|^2\mid\F_n]
&\le
T\int_0^T
\E\left[
\E\left[
\left|s_t^{(n)}\overline G_t^{\theta_n,\phi_n,n}\right|^2
\middle|\G_t^{(n)}
\right]
\middle|\F_n
\right]dt\\
&\le
C\int_0^T
\left(
1+\E[|m_t^{(n)}|^2\mid\F_n]+|\phi_n|^2
\right)dt\\
&\le
C\left(1+|e_0^{(n)}|^4+|\phi_n|^4\right).
\end{align*}

For $I_{2,n}$, by the predictable-version convention and the square-integrability
estimate established in Lemma~\ref{lem:actor-field}, the conditional
It\^o isometry gives
\[
\E[|I_{2,n}|^2\mid\F_n]
=
\E\left[
\int_0^T
|s_t^{(n)}\overline\Gamma_t^{\phi_n,n}|^2dt
\middle|\F_n
\right]
\le C.
\]

For $I_{3,n}$, the integral Cauchy--Schwarz inequality implies
$|I_{3,n}|^2
\le
T\int_0^T|K_t^{(n)}|^2dt$. The KL-gradient square-integrability estimate established in the proof of
Lemma~\ref{lem:actor-field}, applied with $\phi=\phi_n$, therefore yields
\[
\E[|I_{3,n}|^2\mid\F_n]
\le
CT\E\left[
\int_0^T|K_t^{(n)}|^2dt
\middle|\F_n
\right]
\le
C\left(1+|e_0^{(n)}|^2+|\phi_n|^2\right).
\]

Combining the three estimates with
$|x+y+z|^2\le3(|x|^2+|y|^2+|z|^2)$ gives
\begin{align*}
\E\left[
|\widehat F^{\mathrm c}_{n+1}(\phi_n,\theta_n)|^2
\middle|\F_n
\right]
\le
3\sum_{j=1}^3\E[|I_{j,n}|^2\mid\F_n]\le
C\left(1+|e_0^{(n)}|^4+|\phi_n|^4\right),
\end{align*}
which proves \eqref{eq:correct-Fhat-bound}. Finally, since
$\xi_{n+1}=\widehat F^{\mathrm c}_{n+1}-
\E[\widehat F^{\mathrm c}_{n+1}\mid\F_n]$, the conditional variance
identity gives
\begin{align*}
\E[|\xi_{n+1}|^2\mid\F_n]
&=
\E[|\widehat F^{\mathrm c}_{n+1}|^2\mid\F_n]
-\left|\E[\widehat F^{\mathrm c}_{n+1}\mid\F_n]\right|^2\\
&\le
\E[|\widehat F^{\mathrm c}_{n+1}|^2\mid\F_n],
\end{align*}
and \eqref{eq:correct-xi-bound} follows from
\eqref{eq:correct-Fhat-bound}.
\end{proof}

\subsection*{Proof of Theorem~\ref{thm:main-actor}}

\begin{proof}
Let $x_n:=\phi_n-\rhostar $ and $
d_n:=\theta_n-\rhostar$. We first prove almost sure convergence. As $c_n$ is nondecreasing and tends to infinity, $\rhostar\in\Phi_n\subseteq\Phi_{n+1}$ for all sufficiently large $n$. For a closed convex set $C\subset\mathbb R$, the Euclidean projection satisfies
$|\Pi_C(z)-y|\le |z-y|$,
 $y\in C$.
Therefore, we obtain that
\begin{align*}
|x_{n+1}|^2
&\le
\left|
x_n-\alpha_\phi(n)\widehat F^{\mathrm c}_{n+1}(\phi_n,\theta_n)
\right|^2\\
&=
|x_n|^2
-2\alpha_\phi(n) x_n\widehat F^{\mathrm c}_{n+1}(\phi_n,\theta_n)
+(\alpha_\phi(n))^2
\left|\widehat F^{\mathrm c}_{n+1}(\phi_n,\theta_n)\right|^2.
\end{align*}
As $x_n,\theta_n$ are $\F_n$-measurable, taking conditional expectation of the above inequality with respect to $\F_n$ and using \eqref{eq:two-part-decomposition} gives
\begin{align*}
\E[|x_{n+1}|^2\mid\F_n]
\le
|x_n|^2-2\alpha_\phi(n)x_nF_n(\phi_n,\theta_n)+
(\alpha_\phi(n))^2
\E\left[
|\widehat F_{n+1}^c(\phi_n,\theta_n)|^2
\mid\F_n
\right].
\end{align*}
For any fixed $\varepsilon>0$, let
\begin{equation*}
D_\varepsilon:=[-\varepsilon,\varepsilon],
\quad
V_n^{(\varepsilon)}:=\left(|x_n|-\varepsilon\right)_+^2.
\end{equation*}
If the closed intervals $D\subseteq C$, then
\begin{equation}
\operatorname{dist}(\Pi_C(z),D)
\le\operatorname{dist}(z,D).
\label{eq:nested-projection-distance}
\end{equation}
Moreover, since $\Pi_{D_\varepsilon}(x)\in D_\varepsilon$, the definition of distance gives
\begin{align}
\operatorname{dist}(x+h,D_\varepsilon)^2
&\le |x+h-\Pi_{D_\varepsilon}(x)|^2\notag\\
&=\operatorname{dist}(x,D_\varepsilon)^2
+2(x-\Pi_{D_\varepsilon}(x))h+h^2.
\label{eq:squared-distance-one-step}
\end{align}
For fixed $\varepsilon$, when $n$ is sufficiently large,
$D_\varepsilon\subseteq\Phi_{n+1}-\rhostar$. Applying
\eqref{eq:nested-projection-distance}--\eqref{eq:squared-distance-one-step}
to the actor update and taking conditional expectation given $\F_n$ yields
\begin{align}
\E[V_{n+1}^{(\varepsilon)}\mid\F_n]
\le&
V_n^{(\varepsilon)}
-2\alpha_\phi(n)(x_n-\Pi_{D_\varepsilon}(x_n))
\left(K_nx_n+B_nd_n\right)
\notag\\
&+
(\alpha_\phi(n))^2
\E\left[
|\widehat F^{\mathrm c}_{n+1}(\phi_n,\theta_n)|^2
\mid\F_n
\right].
\label{eq:actor-epsilon-RS-pre}
\end{align}
In view of $\theta_n\in\Theta$ and \eqref{eq:explicit-curvature-lower-bound}, $|B_nd_n|\longrightarrow0$.
Thus one may choose $N_\varepsilon$ such that, for all $n\ge N_\varepsilon$,
\[
|B_nd_n|\le
\frac{\kappa_{c,R_0}\varepsilon}{2}.
\]
When $|x_n|>\varepsilon$, $x_n-\Pi_{D_\varepsilon}(x_n)$ has the same sign as $x_n$, and
$|x_n|=|x_n-\Pi_{D_\varepsilon}(x_n)|+\varepsilon$. Write
$\Delta_n^{(\varepsilon)}:=x_n-\Pi_{D_\varepsilon}(x_n)$. It then holds that
\begin{align*}
\Delta_n^{(\varepsilon)}(K_nx_n+B_nd_n)
\ge
\kappa_{c,R_0}|\Delta_n^{(\varepsilon)}|
\left(|\Delta_n^{(\varepsilon)}|+\varepsilon\right)
-|\Delta_n^{(\varepsilon)}|\,|B_nd_n|\ge
\kappa_{c,R_0}|\Delta_n^{(\varepsilon)}|^2.
\end{align*}
When $|x_n|\le\varepsilon$, both sides are zero, so this lower bound holds for all $x_n$. On the other hand, by Lemma~\ref{lem:correct-actor-noise} and $|\phi_n|\le c_n$,
\begin{equation}
\E\left[
|\widehat F^{\mathrm c}_{n+1}(\phi_n,\theta_n)|^2
\mid\F_n
\right]
\le
C\left(1+|e_0^{(n)}|^4+c_n^4\right).
\label{quadratic limit of f}
\end{equation}
By the Gaussian fourth-moment formula, Lemma~\ref{lem:main-filter}, and
\eqref{eq:single-time-scale-summability},
\[
\E[|e_0^{(n)}|^4]
=
\mu_n^4+6\mu_n^2v_n+3v_n^2
\le C(n+1)^{-2\alpha},
\]
and hence
\begin{align*}
\E\left[
\sum_n(\alpha_\phi(n))^2|e_0^{(n)}|^4
\right]
&\le
C\sum_n(\alpha_\phi(n))^2(n+1)^{-2\alpha}<\infty.
\end{align*}
The deterministic part is controlled by
$\sum_n(\alpha_\phi(n))^2 c_n^4<\infty$. Therefore
\begin{equation*}
\sum_n C(\alpha_\phi(n))^2
\left(1+|e_0^{(n)}|^4+c_n^4\right)<\infty
\quad\text{a.s.}
\end{equation*}
Combining Lemma~\ref{lem:correct-actor-noise} with \eqref{eq:actor-epsilon-RS-pre} gives the final recursion
\begin{equation*}
\E[V_{n+1}^{(\varepsilon)}\mid\F_n]
\le
V_n^{(\varepsilon)}
-2\kappa_{c,R_0}\alpha_\phi(n)V_n^{(\varepsilon)}
+C(\alpha_\phi(n))^2
\left(1+|e_0^{(n)}|^4+c_n^4\right).
\end{equation*}
The Robbins--Siegmund theorem \citep{robbins1971convergence} shows that $V_n^{(\varepsilon)}$ converges a.s. and $\sum_n\alpha_\phi(n) V_n^{(\varepsilon)}<\infty$
a.s. As $\sum_n\alpha_\phi(n)=\infty$, it holds that
$V_n^{(\varepsilon)}\to0$ a.s. Taking the probability-one event simultaneously over the countable sequence $\varepsilon_m=1/m$ gives
\[
\limsup_{n\to\infty}|x_n|\le\frac1m
\quad\text{for all }m\ge1,
\]
so $x_n\to0$ a.s., that is, $\phi_n\to\rhostar$ a.s.

We next prove the mean-square rate. Fix $\delta\in(0,2)$.
When $\gamma_2=1$, choose $\delta>0$ sufficiently small that
$(2-\delta)\kappa_{c,R_0}\alpha_\phi>\max\{2\alpha,1\}$. Such a choice is possible because
$\kappa_{c,R_0}\alpha_\phi>\max\{\alpha,1/2\}$. When
$1/2<\gamma_2<1$, any fixed $\delta\in(0,2)$ may be used. Applying the
weighted Young inequality to the critic-bias term gives
\[
2B_n|x_n|\,|d_n|
\le
\delta K_n|x_n|^2
+\frac{B_n^2}{\delta K_n}|d_n|^2.
\]
Therefore,
\begin{align*}
-2x_nF_n(\phi_n,\theta_n)
&\le
-(2-\delta)K_n|x_n|^2
+\frac{B_n^2}{\delta K_n}|d_n|^2\le
-(2-\delta)\kappa_{c,R_0}|x_n|^2
+C(n+1)^{-2\alpha}.
\end{align*}
Substituting \eqref{quadratic limit of f} and the preceding estimate into the conditional recursion gives
\begin{align}
\E[|x_{n+1}|^2\mid\F_n]
\le &
(1-(2-\delta)\kappa_{c,R_0}\alpha_\phi(n))|x_n|^2
+C\alpha_\phi(n)(n+1)^{-2\alpha}\notag\\
&+C(\alpha_\phi(n))^2(1+|e_0^{(n)}|^4+c_n^4).
\label{eq:correct-RS-recursion}
\end{align}
Taking full expectation in \eqref{eq:correct-RS-recursion} gives
\begin{align*}
\E[|x_{n+1}|^2]
\le
\left(1-(2-\delta)\kappa_{c,R_0}\alpha_\phi(n)\right)\E[|x_n|^2]
+C\alpha_\phi(n)(n+1)^{-2\alpha}
+C(\alpha_\phi(n))^2L_n.
\end{align*}
Take $N_0$ sufficiently large so that
$0<(2-\delta)\kappa_{c,R_0}\alpha_\phi(n)<1$ for all $n\ge N_0$.
Iterating the above display from $N_0$ to $n-1$ gives
\begin{align*}
\E[|x_n|^2]
\le&
\prod_{j=N_0}^{n-1}
\left(1-(2-\delta)\kappa_{c,R_0}\alpha_\phi(j)\right)
\E[|x_{N_0}|^2]\\
&+
\sum_{k=N_0}^{n-1}
\prod_{j=k+1}^{n-1}
\left(1-(2-\delta)\kappa_{c,R_0}\alpha_\phi(j)\right)
\left[
C\alpha_\phi(k)(k+1)^{-2\alpha}
+C(\alpha_\phi(k))^2L_k
\right].
\end{align*}
First consider $1/2<\gamma_2<1$. Since $\beta=1$ and $\alpha_\phi(j)=\frac{\alpha_\phi}{(j+1)^{\gamma_2}}$, and since
\[
\prod_{j=k+1}^{n-1}
\left(1-(2-\delta)\kappa_{c,R_0}\alpha_\phi(j+1)^{-\gamma_2}\right)
\le
\exp\left\{-c\sum_{j=k+1}^{n-1}(j+1)^{-\gamma_2}\right\},
\]
since $\phi_{N_0}\in\Phi_{N_0}=[-c_{N_0},c_{N_0}]$,
$\E[|x_{N_0}|^2]\le(c_{N_0}+|\rhostar|)^2<\infty$.
The initial term is controlled by $C\exp\{-c n^{1-\gamma_2}\}$. For the convolution term, if
$N_0\le k\le\lfloor n/2\rfloor$, then likewise
$\sum_{j=k+1}^{n-1}(j+1)^{-\gamma_2}\ge c n^{1-\gamma_2}$,
so this part is also controlled by $C\exp\{-c n^{1-\gamma_2}\}$. If
$\lfloor n/2\rfloor<k\le n-1$, then $L_k\le L_n$, and the local discrete kernel satisfies
\begin{align*}
&\sum_{k=\lfloor n/2\rfloor}^{n-1}
\exp\left\{-c\sum_{j=k+1}^{n-1}(j+1)^{-\gamma_2}\right\}
(k+1)^{-\gamma_2-2\alpha}
\le
C(n+1)^{-2\alpha},\\
&\sum_{k=\lfloor n/2\rfloor}^{n-1}
\exp\left\{-c\sum_{j=k+1}^{n-1}(j+1)^{-\gamma_2}\right\}
L_k(k+1)^{-2\gamma_2}
\le
CL_n(n+1)^{-\gamma_2}.
\end{align*}
The two estimates here come from the same local scale: when $k$ is close to $n$, the effective window length of $\sum_{j=k+1}^{n-1}(j+1)^{-\gamma_2}$ is $(n+1)^{\gamma_2}$. Combining the terms above and absorbing the stretched exponential term into an arbitrarily high-order polynomial decay term, we obtain
\[
\E[|x_n|^2]
\le
CL_n\left[
(n+1)^{-2\alpha}
+(n+1)^{-\gamma_2}
\right].
\]
When $\gamma_2=1$, the same expansion degenerates to the harmonic case. By the choice of $\delta$ above,
$(2-\delta)\kappa_{c,R_0}\alpha_\phi>\max\{2\alpha,1\}$, and
\[
\prod_{j=k+1}^{n-1}
\left(1-(2-\delta)\kappa_{c,R_0}\alpha_\phi(j+1)^{-1}\right)
\le C\left(\frac{k+2}{n+1}\right)^{(2-\delta)\kappa_{c,R_0}\alpha_\phi},
\]
therefore the initial term is controlled by
$C(n+1)^{-(2-\delta)\kappa_{c,R_0}\alpha_\phi}$, and
\begin{align*}
&\sum_{k=N_0}^{n-1}
\left(\frac{k+2}{n+1}\right)^{(2-\delta)\kappa_{c,R_0}\alpha_\phi}
(k+1)^{-1-2\alpha}
\le C(n+1)^{-2\alpha},\\
&\sum_{k=N_0}^{n-1}
\left(\frac{k+2}{n+1}\right)^{(2-\delta)\kappa_{c,R_0}\alpha_\phi}
L_k(k+1)^{-2}
\le CL_n(n+1)^{-1}.
\end{align*}
Here the first term uses
$(2-\delta)\kappa_{c,R_0}\alpha_\phi>2\alpha$, and the second term uses
$(2-\delta)\kappa_{c,R_0}\alpha_\phi>1$. Hence, when $\gamma_2=1$, we likewise have
\[
\E[|x_n|^2]\le
CL_n\left[
(n+1)^{-2\alpha}
+(n+1)^{-\gamma_2}
\right].
\]
Furthermore, let $r_n:=
\sqrt{L_n}
\left[
(n+1)^{-\alpha}
+(n+1)^{-\gamma_2/2}
\right]$. 
For any $M>0$, Markov's inequality gives
\[
\mathbb P\left(
|\phi_n-\rhostar|>Mr_n
\right)
\le
\frac{\E[|\phi_n-\rhostar|^2]}{M^2r_n^2}
\le
\frac{C}{M^2}.
\]
Letting $M\to\infty$ gives
\begin{equation}
\phi_n-\rhostar
=
O_p\left(
\sqrt{L_n}
\left[
(n+1)^{-\alpha}
+(n+1)^{-\gamma_2/2}
\right]
\right).
\label{eq:phi-Op-rate}
\end{equation}
\end{proof}

\subsection*{Proof of Theorem~\ref{thm:main-policy-center}}
\begin{proof}
First prove the coefficient estimates. By \eqref{eq:pc-hc-small}, both required estimates hold uniformly for $0\le t\le T$. Hence
\[
a_t^{(n)}-a_c
=\frac{2p_c(\tau_t^{(n)})}{c+2}
=O((n+1)^{-\alpha}),
\]
and
\[
b_t^{(n)}-b_c
=\frac{2(h_c(\tau_t^{(n)})-1)}{c+2}
=O((n+1)^{-\alpha}).
\]
We next estimate $M_t^{\phi_n}-\rhostar$. By \eqref{eq:m-explicit},
\begin{equation*}
\begin{split}
M_t^{\phi_n}-\rhostar
=
\left(a_t^{(n)}+b_t^{(n)}-1\right)\rhostar
+\frac{a_t^{(n)}p_n}{p_n+t}e_0^{(n)}
+\frac{a_t^{(n)}}{p_n+t}W_t^{(n)}
+b_t^{(n)}(\phi_n-\rhostar).
\end{split}
\end{equation*}
By the two coefficient estimates above and $a_c+b_c=1$,
\[
|a_t^{(n)}+b_t^{(n)}-1|
\le |a_t^{(n)}-a_c|+|b_t^{(n)}-b_c|
\le C(n+1)^{-\alpha}.
\]
Then, by $(z_1+\cdots+z_4)^2\le4\sum_i z_i^2$ and the uniform boundedness of $a_t^{(n)}$, $b_t^{(n)}$, and ${p_n}/{(p_n+t)}$, we obtain
\begin{align*}
\E|M_t^{\phi_n}-\rhostar|^2
&\le
C(n+1)^{-2\alpha}
+C\E[(e_0^{(n)})^2]
+C\E\left[\frac{(W_t^{(n)})^2}{(p_n+t)^2}\right]
+C\E[|\phi_n-\rhostar|^2]\\
&\le
C(n+1)^{-2\alpha}
+C\E[(e_0^{(n)})^2]
+C\bigl(R_0^{(n)}\bigr)^2\E[(W_t^{(n)})^2]
+C\E[|\phi_n-\rhostar|^2].
\end{align*}
Here we used $(p_n+t)^{-1}\le p_n^{-1}=R_0^{(n)}$. Since
$\E[(W_t^{(n)})^2]=t\le T$, by Lemma~\ref{lem:main-filter} and
\eqref{eq:main-actor-rate},
\begin{align*}
\E|M_t^{\phi_n}-\rhostar|^2
&\le
C(n+1)^{-2\alpha}
+C(n+1)^{-\alpha}
+CL_n\left[
(n+1)^{-2\alpha}
+(n+1)^{-\gamma_2}
\right].
\end{align*}
Since $L_n\ge1$ and $(n+1)^{-2\alpha}\le(n+1)^{-\alpha}$, the right-hand side can be controlled by
\[
CL_n\left[
\frac1{(n+1)^\alpha}
+\frac1{(n+1)^{\gamma_2}}
\right].
\]
This proves \eqref{eq:Mt-limit-center-variance}. Finally, the right-hand side of
\eqref{eq:Mt-limit-center-variance} converges to zero because $\alpha,\gamma_2>0$ and
$L_n$ grows only at order $\log\log n$. Thus the Cauchy--Schwarz inequality gives
\[
|\E[M_t^{\phi_n}]-\rhostar|
=|\E[M_t^{\phi_n}-\rhostar]|
\le
\left(\E|M_t^{\phi_n}-\rhostar|^2\right)^{1/2}
\longrightarrow0.
\]
At the same time,
\[
\Var(M_t^{\phi_n})
=\E\left|M_t^{\phi_n}-\E[M_t^{\phi_n}]\right|^2
\le
\E|M_t^{\phi_n}-\rhostar|^2
\longrightarrow0.
\]
\end{proof}

\subsection*{Proof of Lemma~\ref{lem:main-sr-local}}

\begin{proof}
We first compute the derivative of $SR_\infty$ at $\rhostar$. At $M=\rhostar$, we have $SR_\infty(\rhostar)=\sqrt{e^{(\rhostar)^2T}-1}=SR_{\mathrm{or}}$. Let
$h(M):=SR_\infty'(M)/SR_\infty(M)$. Using the explicit formula of, we get $SR_\infty$,
\begin{equation}
h(M)
=
\frac{\rhostar T e^{\rhostar MT}}{e^{\rhostar MT}-1}
-
\frac{MT e^{M^2T}}{e^{M^2T}-1}.
\label{eq:SR-log-derivative}
\end{equation}
At $M=\rhostar$, $h(\rhostar)=0$, and hence
$SR_\infty'(\rhostar)=0$. We next compute the second derivative. Differentiating \eqref{eq:SR-log-derivative} leads to
\begin{align*}
h'(M)
={}&
-\frac{(\rhostar)^2T^2e^{\rhostar MT}}
{\left(e^{\rhostar MT}-1\right)^2}-
\frac{Te^{M^2T}}{e^{M^2T}-1}
+
\frac{2M^2T^2e^{M^2T}}
{\left(e^{M^2T}-1\right)^2}.
\end{align*}
Setting $M=\rhostar$, we get
\begin{align*}
h'(\rhostar)
&=
\frac{(\rhostar)^2T^2e^{(\rhostar)^2T}}
{(e^{(\rhostar)^2T}-1)^2}
-\frac{Te^{(\rhostar)^2T}}{e^{(\rhostar)^2T}-1}=
-T
\frac{
e^{(\rhostar)^2T}
\left(e^{(\rhostar)^2T}-1-(\rhostar)^2T\right)
}
{(e^{(\rhostar)^2T}-1)^2}.
\end{align*}
On the other hand, as $SR_\infty'=SR_\infty h$, it holds that $SR_\infty''=SR_\infty(h^2+h')$.
From $h(\rhostar)=0$ and
$SR_\infty(\rhostar)=\sqrt{e^{(\rhostar)^2T}-1}$, it follows that
\[
SR_\infty''(\rhostar)
=
-T
\frac{
e^{(\rhostar)^2T}
\left(e^{(\rhostar)^2T}-1-(\rhostar)^2T\right)
}
{(e^{(\rhostar)^2T}-1)^{3/2}}.
\]
Strict convexity ensures  $e^{(\rhostar)^2T}-1-(\rhostar)^2T>0$, and it thus holds that
$SR_\infty''(\rhostar)<0$.

We next prove the local Taylor bound. In view of $\rhostar\ne0$, when
$|M-\rhostar|\le|\rhostar|/4$ we have
$|M|\ge3|\rhostar|/4>0$. Hence $SR_\infty$ is twice continuously differentiable on this closed neighborhood, and $SR_\infty''$ is bounded. It holds by Taylor's theorem that
\[
|SR_\infty(M)-SR_{\mathrm{or}}|
\le
\frac12
\sup_{|z-\rhostar|\le|\rhostar|/4}|SR_\infty''(z)|
|M-\rhostar|^2.
\]

Finally, we prove the ODE stability estimate.
To make the continuous-dependence argument fully explicit, we write the coefficient ODEs from Section 4.2 below. For short, we write
\[
U^{(n)}
:=
\left(
p_G^{(n)},q_G^{(n)},
A_H^{(n)},B_H^{(n)},C_H^{(n)}
\right)
\]
and omit the time subscript. It satisfies the common terminal condition $U^{(n)}(T)=0$ and
\begin{align*}
\dot p_G
&=R(1+a)p_G+\rhostar a,\\
\dot q_G
&=-\frac12R^2p_G^2-R(\rhostar-b\phi)p_G+\rhostar b\phi,\\
\dot A_H
&=-R^2A_H^2+2R(1+2a)A_H-2a^2,\\
\dot B_H
&=-R^2A_HB_H-R(\rhostar-2b\phi)A_H
+R(1+2a)B_H-2ab\phi+2\rhostar a,\\
\dot C_H
&=-\frac12R^2(A_H+B_H^2)
-R(\rhostar-2b\phi)B_H-b^2\phi^2-V+2\rhostar b\phi.
\end{align*}
Here $R=R_t^{(n)}$, $V=V_t^{(n)}$, $a=a_t^{(n)}$, and $b=b_t^{(n)}$. Let
$F^{(n)}(t,U;\phi)$ denote the right-hand side of the ODE above. Substituting $R=V=0$, $a=a_c$, and $b=b_c$ gives the limiting vector field
\[
F^{(\infty)}(t,U;\phi)
=
\left(
\rhostar a_c,
\rhostar b_c\phi,
-2a_c^2,
2a_c(\rhostar-b_c\phi),
2\rhostar b_c\phi-b_c^2\phi^2
\right),
\]
where the superscript $(\infty)$ denotes the limiting system with $R=V=0$.
This vector field is independent of $U$, and the explicit solution under the same terminal condition is
\begin{align}
p_G^{(\infty)}(t)
&=
-\rhostar a_c(T-t),
\label{eq:p0-sol}\\
q_G^{(\infty)}(t)
&=
-\rhostar b_c\phi(T-t),
\label{eq:q0-sol}\\
A_H^{(\infty)}(t)
&=
2a_c^2(T-t),
\label{eq:A20-sol}\\
B_H^{(\infty)}(t)
&=
-2a_c(\rhostar-b_c\phi)(T-t),
\label{eq:B20-sol}\\
C_H^{(\infty)}(t)
&=
\left(
b_c^2\phi^2
-2\rhostar b_c\phi
\right)(T-t),
\label{eq:C20-sol}
\end{align}
and write
$U^{(\infty)}:=(p_G^{(\infty)},q_G^{(\infty)},A_H^{(\infty)},B_H^{(\infty)},C_H^{(\infty)})$.

Let
$M:=\overline M_c(m,\phi)=a_c m+b_c\phi$. By
\eqref{eq:p0-sol}--\eqref{eq:C20-sol}, the normalized first and second moment factors corresponding to the limiting system are
\begin{align*}
G^{(\infty)}(m,\phi)
=
\exp\{-\rhostar MT\},\quad
H^{(\infty)}(m,\phi)
=
\exp\{M^2T-2\rhostar MT\}.
\end{align*}
At the oracle point $(m,\phi)=(\rhostar,\rhostar)$, we have $M=\rhostar$, and hence
\[
H^{(\infty)}-(G^{(\infty)})^2
=
e^{-2(\rhostar)^2T}
\left(e^{(\rhostar)^2T}-1\right)>0.
\]
Since
$H^{(\infty)}(m,\phi)-G^{(\infty)}(m,\phi)^2$
is continuous in $(m,\phi)$, one may choose
$\delta_{\mathrm{ODE}}>0$ and $d_0>0$ such that, whenever
$|m-\rhostar|\vee|\phi-\rhostar|\le\delta_{\mathrm{ODE}}$,
$H^{(\infty)}(m,\phi)-G^{(\infty)}(m,\phi)^2\ge2d_0$.
In what follows, fix the neighborhood selected in this way:
$|m-\rhostar|\vee|\phi-\rhostar|\le\delta_{\mathrm{ODE}}$. By continuity of the explicit solution above with respect to $(t,\phi)$, the family of trajectories $\Gamma_\infty
:=
\left\{
U^{(\infty)}(t;\phi):
0\le t\le T,
\ |\phi-\rhostar|\le\delta_{\mathrm{ODE}}
\right\}$
is compact. Choose $\varepsilon_{\mathcal K}>0$ and define the closed tubular neighborhood
$\mathcal K
:=
\left\{
U\in\mathbb R^5:
\operatorname{dist}(U,\Gamma_\infty)
\le\varepsilon_{\mathcal K}
\right\}$. As $\Gamma_\infty$ is compact, $\mathcal K$ is also compact. Combining the uniform bounds on $R,V,a,b$ with boundedness of $\mathcal K$, the polynomial vector field in $U$ has a uniformly bounded Jacobian. Therefore there exists $L<\infty$, independent of $n$, $t$, and $\phi$ in the neighborhood, such that for $U,\widetilde U\in\mathcal K$,
\[
|F^{(n)}(t,U;\phi)-F^{(n)}(t,\widetilde U;\phi)|
\le L|U-\widetilde U|.
\]
On the other hand, using \eqref{eq:decay-R-path}, \eqref{eq:V-bounds}, and \eqref{eq:main-policy-coefficient-rate} term by term, for $U\in\mathcal K$ we have
\begin{align*}
|F_{p_G}^{(n)}-F_{p_G}^{(\infty)}|
&\le C\left(R|p_G|+|a-a_c|\right),\\
|F_{q_G}^{(n)}-F_{q_G}^{(\infty)}|
&\le C\left(R^2|p_G|^2+R(1+|p_G|)
+|b-b_c|\right),\\
|F_{A_H}^{(n)}-F_{A_H}^{(\infty)}|
&\le C\left(R^2|A_H|^2+R|A_H|
+|a-a_c|\right),\\
|F_{B_H}^{(n)}-F_{B_H}^{(\infty)}|
&\le C\left(R^2|A_HB_H|+R(1+|A_H|+|B_H|)
+|a-a_c|
+|b-b_c|\right),\\
|F_{C_H}^{(n)}-F_{C_H}^{(\infty)}|
&\le C\left(R^2(1+|A_H|+|B_H|^2)
+R(1+|B_H|)+V
+|b-b_c|\right).
\end{align*}
Since $\mathcal K$ is bounded, the five estimates above jointly give
\[
\sup_{\substack{0\le t\le T,\ U\in\mathcal K\\
|\phi-\rhostar|\le\delta_{\mathrm{ODE}}}}
|F^{(n)}(t,U;\phi)-F^{(\infty)}(t,U;\phi)|
\le\frac{C}{(n+1)^\alpha}.
\]
Let $[t_n^\ast,T]$ be the maximal backward interval from $T$ on which $U^{(n)}$ remains in $\mathcal K$. Since the two solutions have the same terminal values, on this interval
\begin{align*}
|U_t^{(n)}-U_t^{(\infty)}|
&\le
\int_t^T
|F^{(n)}(s,U_s^{(n)};\phi)
-F^{(n)}(s,U_s^{(\infty)};\phi)|ds\\
&\quad+
\int_t^T
|F^{(n)}(s,U_s^{(\infty)};\phi)
-F^{(\infty)}(s,U_s^{(\infty)};\phi)|ds\\
&\le
L\int_t^T|U_s^{(n)}-U_s^{(\infty)}|ds
+\frac{C}{(n+1)^\alpha}.
\end{align*}
Reverse time in the terminal-value problem by setting $s=T-t$. Then the generalized Gronwall inequality
\citep[Lemma~2.7]{teschl2012ode} gives
\[
\sup_{t_n^\ast\le t\le T}|U_t^{(n)}-U_t^{(\infty)}|
\le Ce^{LT}(n+1)^{-\alpha}.
\]
Choose $n_{1,\mathrm{SR}}$ so that the right-hand side is at most
$\varepsilon_{\mathcal K}/2$. Since
$U_t^{(\infty)}\in\Gamma_\infty$, it follows that
$\operatorname{dist}(U_t^{(n)},\Gamma_\infty)\le
\varepsilon_{\mathcal K}/2$, so $U^{(n)}$ cannot reach
$\partial\mathcal K$. The ODE continuation criterion
\citep[Lemma~2.14 and Corollary~2.15]{teschl2012ode} then implies that $t_n^\ast=0$. Absorbing $e^{LT}$ gives
\begin{equation}
\sup_{0\le t\le T}|U_t^{(n)}-U_t^{(\infty)}|
\le\frac{C}{(n+1)^\alpha}.
\label{eq:ODE-Gronwall}
\end{equation}
This also gives existence and uniform boundedness of the perturbed ODE solution. For
$k\in\mathbb N_0$, define the normalized first and second moment factors
\begin{align*}
G^{(k)}(m,\phi)
&:=
\exp\{p_G^{(k)}(0)m+q_G^{(k)}(0)\},
\\
H^{(k)}(m,\phi)
&:=
\exp\left\{
\frac12A_H^{(k)}(0)m^2+B_H^{(k)}(0)m+C_H^{(k)}(0)
\right\}.
\end{align*}
For $k=\infty$, retain the definitions of
$G^{(\infty)}$ and $H^{(\infty)}$ given above. On the chosen compact neighborhood of $(m,\phi)$, \eqref{eq:ODE-Gronwall} gives
\[
|\log G^{(n)}-\log G^{(\infty)}|
+|\log H^{(n)}-\log H^{(\infty)}|
\le\frac{C}{(n+1)^\alpha}.
\]
These exponents are uniformly bounded on the compact neighborhood, and the exponential function is Lipschitz on bounded intervals. Hence
\[
|G^{(n)}-G^{(\infty)}|+|H^{(n)}-H^{(\infty)}|
\le\frac{C}{(n+1)^\alpha}.
\]
Because $w>x_0$, the common factor $w-x_0$ in the Sharpe ratio cancels exactly, and therefore
\begin{equation*}
\mathcal S(m,\phi;R^{(n)})
=
\frac{1-G^{(n)}(m,\phi)}
{\sqrt{H^{(n)}(m,\phi)-G^{(n)}(m,\phi)^2}}.
\end{equation*}
By the explicit expressions for $G^{(\infty)}$ and $H^{(\infty)}$ given above,
\[
\frac{1-G^{(\infty)}(m,\phi)}
{\sqrt{H^{(\infty)}(m,\phi)-G^{(\infty)}(m,\phi)^2}}
=
\frac{e^{\rhostar MT}-1}
{\sqrt{e^{M^2T}-1}}
=SR_\infty(M).
\]
Moreover, by the choice of $\delta_{\mathrm{ODE}}$ and $d_0$ above,
$H^{(\infty)}-(G^{(\infty)})^2
\ge2d_0$. Since
\begin{align*}
\left|
(H^{(n)}-(G^{(n)})^2)
-(H^{(\infty)}-(G^{(\infty)})^2)
\right|
&\le
|H^{(n)}-H^{(\infty)}|
+|(G^{(n)})^2-(G^{(\infty)})^2|\\
&\le\frac{C}{(n+1)^\alpha},
\end{align*}
take $n_{2,\mathrm{SR}}$ such that, for all $n\ge n_{2,\mathrm{SR}}$,
$
{C}/{(n+1)^\alpha}\le d_0$. Then for all $n\ge \max\{n_{1,\mathrm{SR}},n_{2,\mathrm{SR}}\}$, we have $H^{(n)}-(G^{(n)})^2\ge d_0$. Hence the map $\Psi(G,H):=\frac{1-G}{\sqrt{H-G^2}}$
is continuously differentiable on the compact set under consideration, and its gradient is uniformly bounded. The mean value theorem gives
\[
\left|
\Psi(G^{(n)},H^{(n)})-\Psi(G^{(\infty)},H^{(\infty)})
\right|
\le\frac{C}{(n+1)^\alpha}.
\]
Combining this with the limiting identity above gives
the convergence rate.
\end{proof}

\subsection*{Proof of Theorem~\ref{thm:main-sr}}
\begin{proof}
Take $\delta:=\min\left\{{|\rhostar|}/{4},\,\delta_{\mathrm{ODE}}\right\}$, and define
$\mathcal E_n(\delta)
:=
\left\{
|m_0^{(n)}-\rhostar|\le\delta,
\ |\phi_n-\rhostar|\le\delta
\right\}$.
Denote $r_{\phi,n}:=
\sqrt{
L_n\left[
(n+1)^{-2\alpha}+(n+1)^{-\gamma_2}
\right]
}$. By Lemma~\ref{lem:main-filter} and \eqref{eq:phi-Op-rate}, we have $e_0^{(n)}=O_p((n+1)^{-\alpha/2})$ and
$\phi_n-\rhostar=O_p(r_{\phi,n})$. Therefore
\[
\overline M_{c,n}-\rhostar
=a_c e_0^{(n)}+b_c(\phi_n-\rhostar)
=O_p((n+1)^{-\alpha/2}+r_{\phi,n}).
\]
On the event $\mathcal E_n(\delta)$, by $a_c+b_c=1$,
\[
|\overline M_{c,n}-\rhostar|
\le a_c|e_0^{(n)}|+b_c|\phi_n-\rhostar|
\le\delta\le\frac{|\rhostar|}{4}.
\]
Thus, by the Taylor bound in Lemma \ref{lem:main-sr-local},
\[
\mathbf1_{\mathcal E_n(\delta)}
|SR_\infty(\overline M_{c,n})-SR_{\mathrm{or}}|
=O_p\left((n+1)^{-\alpha}+r_{\phi,n}^2\right).
\]
On the other hand, since $\delta\le\delta_{\mathrm{ODE}}$, \eqref{eq:main-sr-ode-stability} gives
\[
\mathbf1_{\mathcal E_n(\delta)}
|SR_{n}-SR_\infty(\overline M_{c,n})|
\le C(n+1)^{-\alpha}.
\]
Combining the two displays and absorbing the term controlled by $(n+1)^{-2\alpha}$, we obtain
\[
\mathbf1_{\mathcal E_n(\delta)}
|SR_n-SR_{\mathrm{or}}|
=
O_p\left(
\frac{L_n}{(n+1)^\alpha}
+
\frac{L_n}{(n+1)^{\gamma_2}}
\right).
\]
Write $r_{SR,n}:=
\frac{L_n}{(n+1)^\alpha}
+
\frac{L_n}{(n+1)^{\gamma_2}}$.
For any $\varepsilon>0$, take $K<\infty$ such that
\[
\limsup_{n\to\infty}
\mathbb P\left(
\mathbf1_{\mathcal E_n(\delta)}
|SR_n-SR_{\mathrm{or}}|>Kr_{SR,n}
\right)
\le\varepsilon.
\]
Also, by Lemma \ref{lem:main-filter} and Theorem \ref{thm:main-actor},
\[
\mathbb P(\mathcal E_n(\delta)^c)
\le
\mathbb P(|e_0^{(n)}|>\delta)
+
\mathbb P(|\phi_n-\rhostar|>\delta)
\longrightarrow0.
\]
Therefore
\[
\begin{aligned}
\limsup_{n\to\infty}
\mathbb P\left(
|SR_n-SR_{\mathrm{or}}|>Kr_{SR,n}
\right)
&\le
\limsup_{n\to\infty}
\mathbb P\left(
\mathbf1_{\mathcal E_n(\delta)}
|SR_n-SR_{\mathrm{or}}|>Kr_{SR,n}
\right)\\
&\quad+
\limsup_{n\to\infty}
\mathbb P(\mathcal E_n(\delta)^c)
\le\varepsilon.
\end{aligned} 
\]
This is exactly \eqref{eq:main-sr-rate}.
\end{proof}

\subsection*{Proof of Theorem~\ref{thm:main-w}}

We first have the next preparation. 

\begin{lemma}
\label{lem:w-moment-comparison}
Let $\widetilde G$ and $\widetilde H$ be the sampled first- and second-moment factors of $(X_T-w)/(x_0-w)$ under \eqref{micro wealth}. Episode and batch indices are suppressed in $\widetilde G$ and $\widetilde H$. There exist constants $R_*>0$ and $C<\infty$, depending on the fixed grid and model coefficients, such that all four moments are finite and positive for $0<R\le R_*$ and
\begin{align}
\left|\log\frac{\widetilde G}{G(0,m)}\right|
&\le CR(1+|m|+|\phi|),
\label{eq:w-first-moment-comparison}\\
\left|\log\frac{\widetilde H}{H(0,m)}\right|
&\le CR(1+m^2+\phi^2).
\label{eq:w-second-moment-comparison}
\end{align}
The constants are independent of $m$ and $\phi$.
\end{lemma}

\begin{proof}
By \eqref{micro wealth}, It\^o's formula gives $\frac{Y_T}{Y_0}
=\exp\left\{-\int_0^T\left(\rhostar\rho_t+\frac12\rho_t^2\right)dt-\int_0^T\rho_t\,dW_t\right\}$. Condition throughout on the episode-start information $\F_n$. The conditional Gaussian sampling law gives
\[
\E\left[\int_0^T(\rho_t-M_t^\phi)^2dt\middle|\F_n\right]
=\int_0^T V_tdt\le CR_0.
\]
Thus the centered action is square-integrable in time almost surely. Since $V_t$ is deterministic given $\F_n$, the independent sampling construction in Algorithm~\ref{alg:offline_learning} allows us to average first over the market noise:
\begin{align*}
\log\E\left[\frac{Y_T}{Y_0}\middle|\F_n,\{\rho_s-M_s^\phi:0\le s\le T\}\right]
&=\widetilde p_G(0)m+\widetilde q_G(0),\\
\log\E\left[\left(\frac{Y_T}{Y_0}\right)^2\middle|\F_n,\{\rho_s-M_s^\phi:0\le s\le T\}\right]
&=\tfrac12\widetilde A_H(0)m^2+\widetilde B_H(0)m+\widetilde C_H(0).
\end{align*}
The tilded coefficients depend on the centered sampling randomness.
For the macroscopic coefficients, the $A_H$ equation and terminal condition are independent of $\phi$. Thus the stability argument for \eqref{eq:ODE-Gronwall}, with limit \eqref{eq:A20-sol}, bounds $A_H$ uniformly for small $R_0$ and all $\phi$. Together with \eqref{eq:ab-uniform}, the integral formulas $A_H$ and $B_H$, this gives
\[
\sup_{t\le T}(|p_G(t)|+|A_H(t)|)\le C,
\quad \sup_{t\le T}|B_H(t)|\le C(1+|\phi|).
\]
It\^o verification is justified by Novikov's condition: the integrands are affine in the Gaussian filter state plus deterministic $L^2$ terms. Its conditional variance is $O(R_0^2)$ by \eqref{eq:appendix-filter-explicit}, so Jensen's inequality and Gaussian exponential moments give the required integrability for small $R_0$.
Matching the coefficients of $m$ in the first-moment equations and $m^2$ in the second-moment equations gives the same terminal-value problems; uniqueness yields $\widetilde p_G(t)=p_G(t)$ and $\widetilde A_H(t)=A_H(t)$. Variation of constants in the remaining coefficient equations gives
\[
|\widetilde q_G(0)-q_G(0)|
+\sup_{t\le T}|\widetilde B_H(t)-B_H(t)|
\le C\int_0^T|\rho_t-M_t^\phi|dt.
\]
Both differences are centered linear functionals of $\rho-M^\phi$ with uniformly bounded deterministic kernels. Moreover, $\widetilde C_H(0)-C_H(0)+\int_0^T V_tdt$ is the sum of a centered linear functional with deterministic kernel bounded by $C(1+|\phi|)$ and a quadratic functional bounded in absolute value by $C\int_0^T(\rho_t-M_t^\phi)^2dt$.
For fixed $(m,\phi)$, the tower property gives
\begin{align}
\frac{\widetilde G(0,m)}{G(0,m)}
&=\E\left[\exp\{\widetilde q_G(0)-q_G(0)\}\middle|\F_n\right],
\label{eq:w-conditioned-first-ratio}\\
\frac{\widetilde H(0,m)}{H(0,m)}
&=\E\left[\exp\{m(\widetilde B_H(0)-B_H(0))
+\widetilde C_H(0)-C_H(0)\}\middle|\F_n\right].
\label{eq:w-conditioned-second-ratio}
\end{align}
For every integrable $f$ that is deterministic conditional on $\F_n$, weighted Jensen's inequality and the conditional Gaussian law \eqref{eq:actor-draw-centered} yield
\[
\E\left[\exp\left\{\int_0^T f_t(\rho_t-M_t^\phi)dt\right\}\middle|\F_n\right]
\le\exp\left\{\frac12\left(\int_0^T|f_t|\sqrt{V_t}\,dt\right)^2\right\}.
\]
Since the first exponent is a centered linear functional with a bounded coefficient and $V_t\le CR_0$,
\[
1\le\E\left[e^{\widetilde q_G(0)-q_G(0)}\middle|\F_n\right]\le e^{CR_0}.
\]
The second exponent has a centered linear part with kernel bounded by $C(1+|m|+|\phi|)$. Its quadratic part satisfies the preceding bound. For small $R_0$,
\[
\begin{aligned}
\E\left[\exp\left\{C\int_0^T(\rho_t-M_t^\phi)^2dt\right\}\middle|\F_n\right]
&\le\frac1T\int_0^T
\E\left[e^{CT(\rho_t-M_t^\phi)^2}\middle|\F_n\right]dt\\
&=\frac1T\int_0^T(1-2CTV_t)^{-1/2}dt\le e^{CR_0}.
\end{aligned}
\]
Thus Jensen's inequality for the lower bound and Cauchy--Schwarz for the upper bound give
\[
0\le\log\frac{\widetilde G(0,m)}{G(0,m)}\le CR_0,
\quad
-CR_0\le\log\frac{\widetilde H(0,m)}{H(0,m)}\le CR_0(1+m^2+\phi^2).
\]
These estimates ensure finiteness for sufficiently small $R_*$; positivity follows from the stochastic exponential representations.
At $R_0=0$, the action is constant, $\rho_t=a_cm+b_c\phi$, and \eqref{eq:p0-sol}--\eqref{eq:C20-sol} give
\[
\left.\widetilde G(0,m)\right|_{R_0=0}=\left.G(0,m)\right|_{R_0=0}=G^{(\infty)}(m,\phi),
\quad
\left.\widetilde H(0,m)\right|_{R_0=0}=\left.H(0,m)\right|_{R_0=0}=H^{(\infty)}(m,\phi).
\]
The preceding logarithmic bounds prove \eqref{eq:w-first-moment-comparison}--\eqref{eq:w-second-moment-comparison}.
\end{proof}

In the multiplier update, \(G,H\) and their ODE coefficients correspond to the macroscopic model with actor parameter \(\phi_\ell\) and covariance path \(R^{(\ell)}\). Let $\F_r^w$ denote the information before the start of the batch used for the $(r+1)$st multiplier update.
Within this batch, the multiplier remains the $\F_r^w$-measurable variable $w_r$, and is updated after the episodes $\mathcal I_{r+1}$ end. For $\ell\in\mathcal I_{r+1}$,
let $\F_\ell$ denote the information before episode $\ell$ starts; therefore
$\F_r^w\subseteq\F_\ell$. First define the episode-start normalized first-moment factor
\begin{equation*}
\widetilde G_w^{(\ell,r)}
:=
\E\left[
\frac{X_T^{(\ell)}-w_r}{x_0-w_r}
\middle|\F_\ell
\right].
\end{equation*}
Next define the factor truly conditioned on the batch-start information:
\begin{equation*}
G_w^{(\ell,r)}
:=
\E\left[
\frac{X_T^{(\ell)}-w_r}{x_0-w_r}
\middle|\F_r^w
\right]
=
\E\left[
\widetilde G_w^{(\ell,r)}
\middle|\F_r^w
\right],
\quad
\overline G_{w,r}
:=
\frac1{N_w}\sum_{\ell\in\mathcal I_{r+1}}G_w^{(\ell,r)}.
\end{equation*}
The second equality follows from the tower property, and $G_w^{(\ell,r)}$ and
$\overline G_{w,r}$ are both $\F_r^w$-measurable.
As long as $x_0\ne w_r$, the terminal deviation $X_T^{(\ell)}-w_r$ is homogeneous in the initial deviation $x_0-w_r$, and the normalized factor does not depend on $w_r$.
Furthermore,
\begin{equation}
\E[\overline X_{T,r+1}-z\mid\F_r^w]
=
(1-\overline G_{w,r})w_r+\overline G_{w,r}x_0-z.
\label{eq:w-estimator-bias}
\end{equation}
In the oracle limit,
\begin{equation}
G_w^\ast
:=
\exp\left\{-(\rhostar)^2T\right\},
\quad
w^\ast
:=
\frac{z-G_w^\ast x_0}{1-G_w^\ast}
=
\frac{ze^{(\rhostar)^2T}-x_0}
{e^{(\rhostar)^2T}-1}.
\label{eq:w-oracle}
\end{equation}
For any terminal wealth in the window generated by episode-start parameters and using the same batch multiplier $w_r$, define the sampled second-moment factor and recall that
\begin{align*}
\widetilde H^{(\ell,r)}
:=
\E\left[
\left.
\left(\frac{X_T^{(\ell)}-w_r}{x_0-w_r}\right)^2
\right|\F_\ell
\right],\quad
H(0,m_0^{(\ell)})
=
\exp\left\{
\frac12 A_{H\ell,r}(0)(m_0^{(\ell)})^2
+ B_{H\ell,r}(0)m_0^{(\ell)}
+ C_{H\ell,r}(0)
\right\}.
\end{align*}
It does not depend on $w_r$. Conditioning on the batch-start information, it follows from the tower property that
\begin{equation}
H_w^{(\ell,r)}
:=
\E\left[
\widetilde H^{(\ell,r)}\mid\F_r^w
\right],
\quad
\E\left[
\left(X_T^{(\ell)}-w_r\right)^2
\middle|\F_r^w
\right]
=
(x_0-w_r)^2H_w^{(\ell,r)}.
\label{eq:w-second-moment-formula}
\end{equation}

Choose $\ell_0$ so that the comparison lemma applies to all $\ell\ge\ell_0$. For the Riccati coefficient $A_{H\ell,r}$, since the ODE for $A_{H\ell,r}$ does not contain $\phi_\ell$, Lemma~\ref{lem:main-sr-local} gives the $A_H$-component estimate \eqref{eq:ODE-Gronwall} without requiring a local restriction on $\phi$. Therefore, when $\ell$ is sufficiently large, $A_{H\ell,r}$ is uniformly close on $[0,T]$ to $A_H^{(\infty)}(t)=2a_c^2(T-t)$, and hence is uniformly bounded. Thus there exists a deterministic constant $C_A<\infty$ such that
\begin{equation}
\sup_{r,\, \ell\in\mathcal I_{r+1},\,\ell\ge\ell_0}
\sup_{0\le t\le T}|A_{H\ell,r}(t)|\le C_A.
\label{eq:w-A-global-bound}
\end{equation}

Based on the formula above, the ODE for $B_{H\ell,r}$ is linear in $B_{H\ell,r}$. Using \eqref{eq:ab-uniform}, $|\phi_\ell|\le c_\ell$, $R_t^{(\ell)}\le R_0$, and the terminal condition $B_{H\ell,r}(T)=0$, its backward integral form gives
\[
|B_{H\ell,r}(t)|
\le
C\int_t^T|B_{H\ell,r}(s)|\,ds
+C(1+c_\ell),
\quad 0\le t\le T.
\]
The backward Gronwall inequality therefore gives
\begin{equation}
\sup_{0\le t\le T}|B_{H\ell,r}(t)|
\le C(1+c_\ell).
\label{eq:w-B-global-bound}
\end{equation}
Finally, substituting \eqref{eq:w-A-global-bound}--\eqref{eq:w-B-global-bound} into the ODE for $C_{H\ell,r}$ and using $C_{H\ell,r}(T)=0$ gives
\begin{align*}
|C_{H\ell,r}(t)|
\le
C\int_t^T
\left[
1+|B_{H\ell,r}(s)|^2
+c_\ell|B_{H\ell,r}(s)|+c_\ell^2
\right]ds\le C(1+c_\ell^2).
\end{align*}
Consequently, by ODE bounds and \eqref{eq:w-second-moment-comparison}, it holds that
\begin{equation}
\widetilde H^{(\ell,r)}
\le H(0,m_0^{(\ell)})
\exp\left\{CR_0^{(\ell)}\left(1+|m_0^{(\ell)}|^2+\phi_\ell^2\right)\right\}
\le
\exp\left\{C\left(1+|m_0^{(\ell)}|^2+c_\ell^2\right)\right\}.
\label{eq:H-factor-global-growth}
\end{equation}
Conditional on $\F_r^w$, the $m_0^{(\ell)}$ in the window are Gaussian random variables whose means are eventually bounded and whose variances are $O((r+1)^{-2\alpha})$. Since each batch contains only finitely many episodes, the Gaussian quadratic exponential moment formula and \eqref{eq:H-factor-global-growth} show that there exists an a.s. finite random constant $C_H$ such that eventually
\begin{equation}
Q_r
:=
1+\max_{\ell\in\mathcal I_{r+1}}H_w^{(\ell,r)}
\le C_H\exp\{Cc_{(r+1)N_w}^2\}
\quad\text{a.s.}
\label{eq:H-factor-conditional-bound}
\end{equation}
By $w_r\in\mathcal W_r$ and \eqref{eq:w-second-moment-formula}, we obtain
\[
\E\left[
\left|X_T^{(\ell)}\right|^2
\middle|\F_r^w
\right]
\le C(1+c_{w,r}^2)Q_r.
\]
By \eqref{eq:terminal-wealth-batch-statistic} and linearity of conditional expectation, applying Jensen's inequality to the convex function $x\mapsto |x|^2$, we have
\begin{align*}
\left|
\overline X_{T,r+1}
-
\E[\overline X_{T,r+1}\mid\F_r^w]
\right|^2
&\le
\frac1{N_w}
\sum_{\ell\in\mathcal I_{r+1}}
\left|
X_T^{(\ell)}
-
\E[X_T^{(\ell)}\mid\F_r^w]
\right|^2.
\end{align*}
Taking expectation conditional on $\F_r^w$ and using the conditional variance identity,
\begin{align}
\E\left[
\left|
\overline X_{T,r+1}
-
\E[\overline X_{T,r+1}\mid\F_r^w]
\right|^2
\middle|\F_r^w
\right]
&\le
\frac1{N_w}
\sum_{\ell\in\mathcal I_{r+1}}
\E\left[
\left|
X_T^{(\ell)}
-
\E[X_T^{(\ell)}\mid\F_r^w]
\right|^2
\middle|\F_r^w
\right]\notag\\
&\le
\frac1{N_w}
\sum_{\ell\in\mathcal I_{r+1}}
\E\left[
\left|X_T^{(\ell)}\right|^2
\middle|\F_r^w
\right]\notag\\
&\le
C(1+c_{w,r}^2)Q_r.
\label{eq:w-estimator-noise}
\end{align}
Choose a deterministic batch index $r_0$ so that all episodes in batches $r\ge r_0$ have $\ell\ge\ell_0$. Enlarge $\ell_0$ if necessary to make the within-batch Gaussian quadratic exponential moments finite, using their variances $O((r+1)^{-2\alpha})$. Since $\gamma_3>1/2$, using \eqref{eq:H-factor-conditional-bound} gives
\begin{equation}
\sum_{r=r_0}^\infty
(\alpha_w(r))^2(1+c_{w,r}^2)Q_r<\infty
\quad\text{a.s.}
\label{eq:w-noise-summability}
\end{equation}
We next estimate the batch first-moment factor. Similarly, Jensen's inequality gives
\begin{align}
\E\left[
|\overline G_{w,r}-G_w^\ast|^2
\right]
&\le
\frac1{N_w}
\sum_{\ell\in\mathcal I_{r+1}}
\E\left[
\left|G_w^{(\ell,r)}-G_w^\ast\right|^2
\right].
\label{eq:G-factor-Jensen}
\end{align}
By the conditional Jensen inequality and the tower property,
\begin{equation*}
\E\left[
|G_w^{(\ell,r)}-G_w^\ast|^2
\right]
\le
\E\left[
|\widetilde G_w^{(\ell,r)}-G_w^\ast|^2
\right].
\end{equation*}
For the episode-start factor, the first-moment formula in the paper gives
\begin{equation*}
G(0,m_0^{(\ell)})
=
\exp\left\{
p_{G\ell,r}(0)m_0^{(\ell)}
+q_{G\ell,r}(0)
\right\}.
\end{equation*}
The first-moment coefficient ODE listed in the proof of Lemma \ref{lem:main-sr-local} is triangular: the equation for $p_G$ is linear, and after $p_G$ is obtained, $q_G$ is determined by backward integration. By
\eqref{eq:decay-R-path}, \eqref{eq:ab-uniform}, $|\phi_\ell|\le c_\ell$
and the backward Gronwall inequality, it holds that $\sup_{0\le t\le T}|p_{G\ell,r}(t)|\le C$ and
$\sup_{0\le t\le T}|q_{G\ell,r}(t)|\le C(1+c_\ell)$.

Next compare these two ODEs term by term with the $R=0$ limit in \eqref{eq:p0-sol}--\eqref{eq:q0-sol} using the same $\phi_\ell$. By \eqref{eq:decay-R-path} and \eqref{eq:main-policy-coefficient-rate}, the same backward integration estimate gives
\[
\sup_{0\le t\le T}|p_{G\ell,r}(t)-p_G^{(\infty)}(t)|
\le C R_0^{(\ell)},
\quad
\sup_{0\le t\le T}|q_{G\ell,r}(t)-q_G^{(\infty)}(t)|
\le C R_0^{(\ell)}(1+c_\ell).
\]
Also, by $a_c+b_c=1$,
\[
p_G^{(\infty)}(0)m_0^{(\ell)}+q_G^{(\infty)}(0)+(\rhostar)^2T
=-\rhostar T\left[
a_c(m_0^{(\ell)}-\rhostar)
+b_c(\phi_\ell-\rhostar)
\right].
\]
Therefore
\begin{align*}
\left|
p_{G\ell,r}(0)m_0^{(\ell)}+q_{G\ell,r}(0)
+(\rhostar)^2T
\right|
\le
C\left[
|m_0^{(\ell)}-\rhostar|
+|\phi_\ell-\rhostar|
+R_0^{(\ell)}\left(1+|m_0^{(\ell)}|+c_\ell\right)
\right].
\end{align*}
By \eqref{eq:w-first-moment-comparison}, the preceding exponent bound also holds for $|\log\widetilde G_w^{(\ell,r)}+(\rhostar)^2T|$, after increasing $C$. The mean value theorem for the exponential map and the preceding display give
\begin{align}
|\widetilde G_w^{(\ell,r)}-G_w^\ast|^2
\le{}&
C\exp\left\{C\left(|m_0^{(\ell)}-\rhostar|+c_\ell\right)\right\}
\notag\\
&\times\left[
|m_0^{(\ell)}-\rhostar|^2
+|\phi_\ell-\rhostar|^2
+(R_0^{(\ell)})^2
\left(1+|m_0^{(\ell)}-\rhostar|^2+c_\ell^2\right)
\right].
\label{eq:G-factor-global-growth}
\end{align}
The same estimates hold with $c_\ell$ replaced by $|\phi_\ell|$, giving the uniform local bound used below. By the explicit recursion for filtering and the conditional second-moment bound \eqref{eq:correct-Fhat-bound} for the actor, for every realized sample path satisfying $m_0^{(rN_w+1)}\to\rhostar$ and $\phi_{rN_w+1}\to\rhostar$, we have
\begin{equation*}
\Delta_r
:=\max_{\ell\in\mathcal I_{r+1}}
\E\left[
|m_0^{(\ell)}-\rhostar|^2
+|\phi_\ell-\rhostar|^2
+(R_0^{(\ell)})^2
\middle|\F_r^w
\right]
\longrightarrow0.
\end{equation*}
Fix $\delta>0$, and let $\mathcal A_r(\delta)
:=
\left\{
\max_{\ell\in\mathcal I_{r+1}}
\left(
|m_0^{(\ell)}-\rhostar|
\vee|\phi_\ell-\rhostar|
\right)
\le\delta
\right\}$.
Using the local Lipschitz bound on $\mathcal A_r(\delta)$ gives
\[
\max_{\ell\in\mathcal I_{r+1}}
\E\left[
|\widetilde G_w^{(\ell,r)}-G_w^\ast|^2
\mathbf 1_{\mathcal A_r(\delta)}
\middle|\F_r^w
\right]
\le C\Delta_r
\longrightarrow0
\quad\text{a.s.}
\]
On the other hand, batch-start convergence ensures that eventually a.s. every batch starts inside the $\delta/4$-neighborhood:
\[
|m_0^{(rN_w+1)}-\rhostar|
\vee|\phi_{rN_w+1}-\rhostar|\le\delta/4.
\]
Thus, $\mathcal A_r(\delta)^c$ can occur only if a filtering fluctuation or a cumulative actor increment within the window has magnitude at least $3\delta/4$. We control these two cases separately.
By the explicit recursion in the proof of Lemma~\ref{lem:main-filter}, conditional on $\F_r^w$, the variables
$m_0^{(\ell)}-\rhostar$ in the window are Gaussian. Eventually a.s.,
\[
\max_{\ell\in\mathcal I_{r+1}}
|\E[m_0^{(\ell)}-\rhostar\mid\F_r^w]|
\le\delta/4,
\quad
\max_{\ell\in\mathcal I_{r+1}}
\Var(m_0^{(\ell)}-\rhostar\mid\F_r^w)
\le C(r+1)^{-2\alpha}.
\]
Therefore, for every $\ell\in\mathcal I_{r+1}$,
\begin{align*}
&\mathbb P\left(
|m_0^{(\ell)}-\rhostar|>\delta
\middle|\F_r^w
\right)\le
\mathbb P\left(
\left|m_0^{(\ell)}-\rhostar-\E[m_0^{(\ell)}-\rhostar\mid\F_r^w]\right|
>\frac{3\delta}{4}
\middle|\F_r^w
\right).
\end{align*}
The Gaussian tail bound then gives, for some $\kappa_1>0$,
\[
\mathbb P\left(
\max_{\ell\in\mathcal I_{r+1}}|m_0^{(\ell)}-\rhostar|>\delta
\middle|\F_r^w\right)
\le Ce^{-2\kappa_1(r+1)^{2\alpha}}.
\]
On the other hand, since $\phi_k\in\Phi_k\subseteq\Phi_{k+1}$, we have $|\phi_{k+1}-\phi_k|
\le
\alpha_\phi(k)|\widehat F_{k+1}^c(\phi_k,\theta_k)|$.
Together with the batch-start bound, this yields
\[
\left\{
\max_{\ell\in\mathcal I_{r+1}}|\phi_\ell-\rhostar|>\delta
\right\}
\subseteq
\left\{
\sum_{k=rN_w+1}^{(r+1)N_w-1}
\alpha_\phi(k)|\widehat F^{\mathrm c}_{k+1}(\phi_k,\theta_k)|
>\frac{3\delta}{4}
\right\}.
\]
Therefore, the conditional Markov inequality and the Cauchy--Schwarz inequality give
\begin{align*}
\mathbb P\left(
\max_{\ell\in\mathcal I_{r+1}}|\phi_\ell-\rhostar|>\delta
\middle|\F_r^w\right)
&\le \frac{16}{9\delta^2}
\E\left[
\left(
\sum_{k=rN_w+1}^{(r+1)N_w-1}
\alpha_\phi(k)|\widehat F^{\mathrm c}_{k+1}(\phi_k,\theta_k)|
\right)^2
\middle|\F_r^w
\right]\\
&\le \frac{16(N_w-1)}{9\delta^2}
\sum_{k=rN_w+1}^{(r+1)N_w-1}
(\alpha_\phi(k))^2\E\left[
|\widehat F^{\mathrm c}_{k+1}(\phi_k,\theta_k)|^2\mid\F_r^w
\right]\\
&\le C(r+1)^{-2\gamma_2}
\left(1+c_{(r+1)N_w}^4\right)
=O((r+1)^{-2\kappa_2})
\quad\text{a.s.}
\end{align*}
where $0<\kappa_2<\gamma_2$ is arbitrary; the last inequality uses
\eqref{eq:correct-Fhat-bound} and the tower property.
In addition, \eqref{eq:G-factor-global-growth} and the conditional Gaussian exponential-moment bound give eventually a.s.
\[
|\widetilde G_w^{(\ell,r)}-G_w^\ast|^2
\le Ce^{c_{(r+1)N_w}+|m_0^{(\ell)}-\rhostar|}
\left(1+c_{(r+1)N_w}^2\right)
\left(1+|m_0^{(\ell)}-\rhostar|^2\right),
\]
\[
\max_{\ell\in\mathcal I_{r+1}}
\E\left[
e^{2C|m_0^{(\ell)}-\rhostar|}
\left(1+|m_0^{(\ell)}-\rhostar|^2\right)^2
\middle|\F_r^w\right]
\le C.
\]
Thus the conditional Cauchy--Schwarz inequality gives
\begin{align*}
&\max_{\ell\in\mathcal I_{r+1}}
\E\left[
|\widetilde G_w^{(\ell,r)}-G_w^\ast|^2
\mathbf 1_{\mathcal A_r(\delta)^c}
\middle|\F_r^w
\right]\le
\exp\{Cc_{(r+1)N_w}\}
\left[
e^{-\kappa_1(r+1)^{2\alpha}}
+O((r+1)^{-\kappa_2})
\right],
\end{align*}
The last display is clearly $o(1)$. 
Combining the two parts above and using the conditional Jensen inequality,
\begin{align*}
\max_{\ell\in\mathcal I_{r+1}}
|G_w^{(\ell,r)}-G_w^\ast|^2
\le
\max_{\ell\in\mathcal I_{r+1}}
\E\left[
|\widetilde G_w^{(\ell,r)}-G_w^\ast|^2
\middle|\F_r^w
\right]\le C\Delta_r+o(1)
\longrightarrow0
\quad\text{a.s.}
\end{align*}
The finite average therefore satisfies $\overline G_{w,r}\to G_w^\ast$ a.s. Thus we obtain
\begin{equation}
\overline G_{w,r}\longrightarrow G_w^\ast
\quad\text{a.s.}
\label{eq:G-factor-tracking}
\end{equation}

\begin{proof}
It suffices to consider the tail $r\ge r_0$, where the sampled-wealth conditional moments are finite.
Write $U_{w,r}:=w_r-w^\ast$, $\xi_{w,r+1}
:=
\overline X_{T,r+1}
-
\E[\overline X_{T,r+1}\mid\F_r^w]$. By Theorem~\ref{thm:main-w} and \eqref{eq:w-oracle},
\[
(1-G_w^\ast)w^\ast
+G_w^\ast x_0-z=0.
\]
So the conditional mean in \eqref{eq:w-estimator-bias} can be written as
\begin{equation}
\E[\overline X_{T,r+1}-z\mid\F_r^w]
=
(1-\overline G_{w,r})U_{w,r}+\delta_{w,r},
\quad
\delta_{w,r}
:=
(\overline G_{w,r}-G_w^\ast)(x_0-w^\ast).
\label{eq:w-field-decomposition}
\end{equation}
Since $c_{w,r}\uparrow\infty$, there exists a deterministic finite index $r_w$ such that $w^\ast\in\mathcal W_r$ for all $r\ge r_w$.
Therefore, for all $r\ge r_w$, projection non-expansiveness gives
\begin{align*}
|U_{w,r+1}|^2
&\le
\left|
U_{w,r}
-
\alpha_w(r)
\left[
	(1-\overline G_{w,r})U_{w,r}
+\delta_{w,r}
+\xi_{w,r+1}
\right]
\right|^2.
\end{align*}
Since $\overline G_{w,r}\to G_w^\ast<1$, there exists a random but finite $r_w^{\mathrm{stab}}\ge r_w$ such that, for all $r\ge r_w^{\mathrm{stab}}$,
\[
\frac{1-G_w^\ast}{2}
\le
1-\overline G_{w,r}
\le
1+G_w^\ast
:=C_G<\infty.
\]
We again use the squared distance to a fixed neighborhood. Fix
$0<\varepsilon<{(w^\ast-\underline w)}/{2}$, and define
\begin{equation*}
D_\varepsilon^w:=[-\varepsilon,\varepsilon],
\quad
V_{w,r}^{(\varepsilon)}
:=|U_{w,r}-\Pi_{D_\varepsilon^w}(U_{w,r})|^2.
\end{equation*}
Since $c_{w,r}\uparrow\infty$ and $\underline w<w^\ast$, eventually
$D_\varepsilon^w\subseteq\mathcal W_{r+1}-w^\ast$. Applying
\eqref{eq:nested-projection-distance}--\eqref{eq:squared-distance-one-step}
and taking conditional expectation given $\F_r^w$ yields
\begin{align}
\E[V_{w,r+1}^{(\varepsilon)}\mid\F_r^w]
\le{}&
V_{w,r}^{(\varepsilon)}
-2\alpha_w(r)\left(U_{w,r}-\Pi_{D_\varepsilon^w}(U_{w,r})\right)
\left[
(1-\overline G_{w,r})U_{w,r}+\delta_{w,r}
\right]
\notag\\
&+
(\alpha_w(r))^2
\E\left[
\left|
(1-\overline G_{w,r})U_{w,r}
+\delta_{w,r}+\xi_{w,r+1}
\right|^2
\middle|\F_r^w
\right].
\label{eq:w-epsilon-RS-pre}
\end{align}
Let $\kappa_w:=(1-G_w^\ast)/2>0$. By
\eqref{eq:G-factor-tracking} and \eqref{eq:w-field-decomposition}, eventually
we have $1-\overline G_{w,r}\ge\kappa_w$ and
$|\delta_{w,r}|\le{\kappa_w\varepsilon}/{2}$.
As in the actor case, $U_{w,r}-\Pi_{D_\varepsilon^w}(U_{w,r})$ has the same sign as
$U_{w,r}$ when $|U_{w,r}|>\varepsilon$, and hence eventually
\begin{equation}
\left(U_{w,r}-\Pi_{D_\varepsilon^w}(U_{w,r})\right)
\left[
(1-\overline G_{w,r})U_{w,r}+\delta_{w,r}
\right]
\ge
\kappa_wV_{w,r}^{(\varepsilon)}.
\label{eq:w-epsilon-drift}
\end{equation}
Also, by $w_r\in\mathcal W_r$, the eventual bound
$|1-\overline G_{w,r}|\le 1+G_w^\ast$, and
\eqref{eq:w-estimator-noise},
\[
\E\left[
\left|
(1-\overline G_{w,r})U_{w,r}
+\delta_{w,r}+\xi_{w,r+1}
\right|^2
\middle|\F_r^w
\right]
\le
C\left((1+c_{w,r}^2)Q_r+|\delta_{w,r}|^2\right).
\]
Finally, by \eqref{eq:w-noise-summability}, $\delta_{w,r}\to0$, and
$\sum_r(\alpha_w(r))^2<\infty$,
\begin{equation*}
\sum_rC(\alpha_w(r))^2
\left((1+c_{w,r}^2)Q_r+|\delta_{w,r}|^2\right)<\infty
\quad\text{a.s.}
\end{equation*}
Therefore \eqref{eq:w-epsilon-RS-pre}--\eqref{eq:w-epsilon-drift} gives
\begin{equation*}
\E[V_{w,r+1}^{(\varepsilon)}\mid\F_r^w]
\le
V_{w,r}^{(\varepsilon)}
-2\kappa_w\alpha_w(r)V_{w,r}^{(\varepsilon)}
+C(\alpha_w(r))^2
\left((1+c_{w,r}^2)Q_r+|\delta_{w,r}|^2\right).
\end{equation*}
The Robbins--Siegmund theorem \citep{robbins1971convergence} and $\sum_r\alpha_w(r)=\infty$ give
$V_{w,r}^{(\varepsilon)}\to0$ a.s. Finally, take
\[
\varepsilon_m
:=
\min\left\{
\frac1m,\frac{w^\ast-\underline w}{4}
\right\},
\quad m\ge1,
\]
and intersect the countably many probability-one events to obtain
$U_{w,r}\to0$ a.s., that is, $w_r\to w^\ast$ a.s.
\end{proof}

\end{document}